\ifx \labelsONmargin\undefined

  \ifx\RequirePackage\undefined
    \expandafter\ifx\csname amsart.sty\endcsname\relax
        \documentstyle[12pt,amssymb,amscd]{amsart}
    \fi

    \def\atSign{@@}

    \def\mathbb{\Bbb}
    \def\mathfrak{\frak}
    \def\mathbf{\bold}
    \ifx\boldsymbol\undefined	
      \def\boldsymbol#1{{\bold #1}}
    \fi

    \def\mathbit{\boldsymbol}

    \makeatletter
    \newenvironment{proof}{%
         \@ifnextchar[{%
                       \expandafter\let\expandafter\end@proof
                         \csname endpf*\endcsname
                         \my@proof
                      }{\let\end@proof\endpf\pf}%
        }{\end@proof}
    \def\my@proof[#1]{\@nameuse{pf*}{#1}}
    \makeatother

    \def\xrightarrow[#1]#2{@>{#2}>{#1}>}
    \def\xleftarrow[#1]#2{@<{#2}<{#1}<}

    \def\providecommand#1{\def#1}
    \def\emph#1{{\em #1}}
    \def\textbf#1{{\bf #1}}

    \def\mathring{\overset{\,\,{}_\circ}}
  \else
      \ifx\usepackage\RequirePackage
      \else
        \documentclass[12pt]{amsart}
        \usepackage{amssymb,amscd}
      \fi
      \ifx \mathring\undefined		
        \DeclareMathAccent{\mathring}{\mathalpha}{operators}{"17}
      \fi

      \long\def\FAKEendPROOF{\endtrivlist}
      \ifx\endproof\FAKEendPROOF
	  \def\endproof{\qed\endtrivlist}
      \fi

      \ifx\mathbit\undefined	
        \DeclareMathAlphabet{\mathbit}{OML}{cmm}{b}{it}
      \fi

      \def\atSign{@}

      \advance\oddsidemargin -0.1\textwidth
      \evensidemargin\oddsidemargin
      \def\Sb#1\endSb{_{\substack{#1}}}
      \def\Sp#1\endSp{^{\substack{#1}}}

  \fi
\makeatletter
\@ifundefined{DeclareFontShape} 
     {\@ifundefined{selectfont}
             {
                \message{We need AmS-LaTeX, this version of LaTeX does not 
                        support it}\@@end
             }{
                \def\mathcal{\cal}
                
                \def\pcyr{%
                        \def\default@family{UWCyr}%
                        \let\oldSl@\sl
                        \def\sl{\def\default@shape{it}\oldSl@}%
                        \cyracc
                        \language\Russian\family{UWCyr}\selectfont
                }
             }}{
                \DeclareFontEncoding{OT2}{}{\noaccents@}
                \DeclareFontSubstitution{OT2}{cmr}{m}{n}
                \DeclareFontFamily{OT2}{cmr}{\hyphenchar\font45 }
                \DeclareFontShape{OT2}{cmr}{m}{n}{%
                     <5><6><7><8><9><10>gen*wncyr %
                     <10.95><12><14.4><17.28><20.74><24.88> wncyr10 %
                }{}
                \DeclareFontShape{OT2}{cmr}{m}{it}{%
                     <5><6><7><8><9><10> gen * wncyi%
                     <10.95><12><14.4><17.28><20.74><24.88> wncyi10%
                }{}
                \DeclareFontShape{OT2}{cmr}{bx}{n}{%
                     <5><6><7><8><9><10> gen * wncyb%
                     <10.95><12><14.4><17.28><20.74><24.88> wncyb10%
                }{}
                \DeclareFontShape{OT2}{cmr}{m}{sl}{%
                     <-> ssub * cmr/m/it%
                }{}
                \DeclareFontShape{OT2}{cmr}{m}{sc}{%
                     <5><6><7><8><9><10>%
                     <10.95><12><14.4><17.28><20.74><24.88> wncysc10%
                }{}
                \DeclareFontFamily{OT2}{cmss}{\hyphenchar\font45 }
                \DeclareFontShape{OT2}{cmss}{m}{n}{%
                     <8><9><10> gen * wncyss%
                     <10.95><12><14.4><17.28><20.74><24.88> wncyss10%
                }{}
                
                \def\cyrencodingdefault{OT2}
                \def\pcyr{%
                        \cyracc
                        \let\encodingdefault\cyrencodingdefault
                        \language\Russian\fontencoding{OT2}\selectfont
                }
             }   

\@ifundefined{theorembodyfont} 
     {
        \def\theorembodyfont#1{\relax}
        \makeatletter
          \let\@@th@plain\th@plain
          \def\th@plain{ \@@th@plain \slshape }
        \makeatother
        \let\normalshape\relax
     }{}   

\ifx\cprime\undefined
     \def\cprime{$'$}
\fi
\makeatletter
  \def\@sect@my#1#2#3#4#5#6[#7]#8{%
\ifnum #2>\c@secnumdepth
   \let\@svsec\@empty
 \else
   \refstepcounter{#1}%
\edef\@svsec{\ifnum#2<\@m
             \@ifundefined{#1name}{}{\csname #1name\endcsname\ }\fi
\noexpand\rom{\csname the#1\endcsname.}\enspace}\fi
 \@tempskipa #5\relax
 \ifdim \@tempskipa>\z@ 
   \begingroup #6\relax
   \@hangfrom{\hskip #3\relax\@svsec}{\interlinepenalty\@M #8\par}%
   \endgroup
   \if@article\else\csname #1mark\endcsname{%
        \ifnum \c@secnumdepth >#2\relax\csname the#1\endcsname. \fi#7}\fi
\ifnum#2>\@m \else
       \let\@tempf\\ \def\\{\protect\\}\addcontentsline{toc}{#1}%
{\ifnum #2>\c@secnumdepth \else
             \protect\numberline{%
               \ifnum#2<\@m
               \@ifundefined{#1name}{}{\csname #1name\endcsname\ }\fi
               \csname the#1\endcsname.}\fi
           #8}\let\\\@tempf
     \fi
 \else
  \def\@svsechd{#6\hskip #3\@svsec
    \@ifnotempty{#8}{\ignorespaces#8\unskip
       \ifnum\spacefactor<1001.\fi}%
        \ifnum#2>\@m \else
          \let\@tempf\\ \def\\{\protect\\}\addcontentsline{toc}{#1}%
            {\ifnum #2>\c@secnumdepth \else
              \protect\numberline{%
                \ifnum#2<\@m
                \@ifundefined{#1name}{}{\csname #1name\endcsname\ }\fi
                \csname the#1\endcsname.}\fi
             #8}\let\\\@tempf\fi}%
 \fi
\@xsect{#5}}
  \ifx\RequirePackage\undefined
  \let\@sect\@sect@my             
  \fi
  \ifx\RequirePackage\undefined	
  \def\th@remark@my{\theorempreskipamount6\p@\@plus6\p@
    \theorempostskipamount\theorempreskipamount
    \def\theorem@headerfont{\it}\normalshape}
    \let\th@remark\th@remark@my
  \else				
    \let\o@@remark\th@remark

    \ifx\theorempostskipamount\undefined		
    \else						
      \def\th@remark{\o@@remark
	\ifdim\theorempostskipamount < 2pt\relax
	  \theorempostskipamount\theorempreskipamount
	     \multiply\theorempostskipamount\tw@
	     \divide\theorempostskipamount\thr@@
	\fi
      }
    \fi
  \fi
\let\myLabel\@gobble
\def\labelsONmargin{\@mparswitchfalse\def\myLabel##1{\@bsphack\marginpar
                                  {\normalshape\tiny\rm Label ##1}\@esphack}}
\ifx \url\undefined
  \def\url#1{{\tt #1}}%
\fi
\def\cyracc{\def\u##1{
                \if \i##1\char"1A%
                \else \if I##1\char"12%
                \else \accent"24 ##1\fi\fi }%
\def\"##1{\if e##1{\char"1B}%
                \else \if E##1{\char"13}%
                \else \accent"7F ##1\fi\fi }%
\def\9##1{\if##1z\char"19 
\else\if##1Z\char"11 
\else\if##1E\char"03 
\else\if##1e\char"0B 
\else\if##1u\char"18 
\else\if##1U\char"10 
\else\if##1A\char"17 
\else\if##1a\char"1F 
\else\if##1p\char"7E 
\else\if##1P\char"5E 
\else\if##1Q\char"5F 
\else\if##1q\char"7F 
\else\if##1i\char"1A 
\else\if##1I\char"12 
\else\if##1N\char"7D 
\fi
\fi
\fi
\fi
\fi
\fi
\fi
\fi
\fi
\fi
\fi
\fi
\fi
\fi
\fi
}%
\def\cydot{{\kern0pt}}}%
\def\cydot{$\cdot$}

\ifx\Russian\undefined
        \def\Russian{0\relax
    \message{Don't know the hyphenation rules for Russian^^J
                        Please do INITeX with `input  russhyph' in the 
                        command line}%
                \gdef\Russian{0\relax}%
        }
\fi
\ifx\RequirePackage\undefined			
  \def\@putname#1#2#3#4{\def\@@ref{#3}\let\old@bf\bf
        \def\bf##1{\old@bf\if?\noexpand##1?{#4}\else##1\fi}%
	#1{#2}%
        \let\bf\old@bf}
\else						
  \def\@putname#1#2#3#4{\def\@@ref{#3}\let\old@bf\bf	
	\let\old@reset@font\reset@font			
        \def\bf##1{\old@bf\if?\noexpand##1?{#4}\else##1\fi}%
	\def\reset@font##1##2{\old@reset@font##1\if?\noexpand##2?{#4}\else##2\fi}#1{#2}%
        \let\bf\old@bf\let\reset@font\old@reset@font}
\fi
\let\my@ref=\ref
\def\ref#1{\@putname\my@ref{#1}{#1}{\tiny\rm\@@ref}}
\let\my@pageref=\pageref
\def\pageref#1{\@putname\my@pageref{#1}{#1}{\tiny\rm\@@ref}}
\let\my@cite=\cite
\def\cite#1{\@putname\my@cite{#1}{\@citeb}{\tiny\rm\@@ref}}
\makeatother
\fi
\relax 
\theoremstyle{plain} 
\theorembodyfont{\sl}

\newtheorem{nwthrmi}{Fact } 
\newtheorem{nwthrmii}{Fact } 
\newtheorem{nwthrmiii}{Heuristic }

\numberwithin{equation}{section}
\theorembodyfont{\rm}
\theoremstyle{definition}

\newtheorem{definition}{Definition}[section]

\newtheorem{example}[definition]{Example}

\theoremstyle{remark}

\newtheorem{remark}[definition]{Remark} 

\theoremstyle{plain} 
\theorembodyfont{\sl}

\newtheorem{theorem}[definition]{Theorem}
\newtheorem{lemma}[definition]{Lemma}
\newtheorem{corollary}[definition]{Corollary}
\newtheorem{proposition}[definition]{Proposition}
\newtheorem{amplification}[definition]{Amplification}

\begin{document}
\bibliographystyle{amsplain}
\relax 

\title{ Geometric vertex operators }

\author{ Ilya Zakharevich }

\email {math\atSign{}ilyaz.org}

\date{\quad May 1998 (Revision II: August 2002), archived at
\url{arXiv:math.AG/980509}.\quad Printed: \today }

\setcounter{section}{-1}

\maketitle
\begin{abstract}
We investigate relationships between (infinite-dimensional)
algebraic geometry of loop spaces, smooth families of partial linear
operators, and classical vertex operators. Vertex operators are shown to
be actions of birational transformations of infinite-dimensional
algebraic ``varieties'' $ M $ on appropriate line bundles on $ M $.

These vertex operators act in the vector space $ {\mathcal M}\left(M\right) $ of meromorphic
functions on $ M $ as {\em partial operators\/}: they are defined on a subspace (in an
appropriate {\em lattice of subspaces\/} of $ {\mathcal M}\left(M\right) $), and send smooth families in
such a subspace to smooth families. Axiomatizing this, we define
{\em conformal fields\/} as arbitrary families of partial operators in $ {\mathcal M}\left(M\right) $ which
satisfy both these properties.

The ``variety'' $ M $ related to standard vertex operators is formed by
rational functions of one variable $ z\in Z={\mathbb P}^{1} $, changing the variety $ Z $ one
obtains different examples of $ M $, and multidimensional analogues of vertex
operators. One can cover $ {\mathcal M}\left(M\right) $ by ``big smooth subsets''; these subsets are
parameterized by appropriate projective bundles over Hilbert schemes of
points on $ Z $. We deduce conformal associativity relation for conformal
fields, and conformal commutation relations for Laurent coefficients of
commuting conformal fields from geometric properties of the Hilbert
schemes.

One of the main stumbling blocks for a mathematician when dealing
with conformal fields is the very first sentence they see: ``Consider two
holomorphic functions $ a\left(s\right) $ and $ b\left(t\right) $ such that $ a\left(s\right)b\left(t\right) $ has a pole when
$ s=t. $'' We {\em start\/} with examples of such functions (with values in partially
defined operators in spaces of meromorphic functions), thus make this
problem nonextant.

We also provide geometric description of boson-fermion
correspondence, and relations of this correspondence to geometry of the
set of meromorphic functions.
\end{abstract}
\tableofcontents

\section{Introduction }\label{h0.2}\myLabel{h0.2}\relax 

All the geometry and linear algebra considered in this paper happens
over an arbitrary algebraically closed field $ {\mathbb K} $. However, the main
application is when $ {\mathbb K}={\mathbb C} $, thus one can feel free to substitute $ {\mathbb C} $ instead
of $ {\mathbb K} $.

Note however that in a case of a topological field (such as $ {\mathbb C} $) our
weak axioms will not work in the case of an analytic manifold, to
implement field theory over an analytic manifold one may need much more
rigid axioms. We confine ourself to the case of algebraic manifolds and
algebraic functions, so one does not need to worry about significant
singularities and non-dense open sets.

\subsection{Loops spaces and OPEs }\label{s0.01}\myLabel{s0.01}\relax  Fix a variety $ Z $. Given a
variety $ X $, a {\em loop in\/} $ X $ is a mapping $ \varphi\colon Z \to X $ (smooth, or regular, or
formal, or rational, etc., as specified by context). Denote the set of
loops on $ X $ by $ {\mathcal L}X $ or $ {\mathcal L}_{Z}X $; later we equip this set with additional
structures, and call it the {\em loop space}. Assume that these additional
structures on $ {\mathcal L}X $ lead to a particular algebra $ \operatorname{Func}\left({\mathcal L}X\right) $ of functions on $ {\mathcal L}X $
(continuous, or regular, or rational, etc., as needed). To avoid
confusion with the loops (which are functions themselves), call elements
of $ \operatorname{Func}\left({\mathcal L}X\right) $ {\em functionals on\/} $ {\mathcal L}X $.

Loops spaces are the first examples of infinite-dimensional
varieties with properties dramatically different from those familiar
from finite-dimensional algebraic geometry. {\em Vertex\/} flavors of
algebro-geometric constructions are convenient tools helping handling
objects from this still-unchartered territory.

The simplest manifestation of this new geometry is the appearance of
{\em operator product expansions (OPEs)\/} in the rule of composition of families
of geometric linear operators acting on $ \operatorname{Func}\left({\mathcal L}X\right) $. Recall that given a
manifold $ M $, a {\em geometric operator\/} in $ \operatorname{Func}\left(M\right) $ is a composition of an
operator $ {\mathbit F} $ of multiplication by $ F\in\operatorname{Func}\left(M\right) $ and of the operator $ h^{*} $ for a
diffeomorphism $ h\in\operatorname{Diff}\left(M\right) $ of $ M $. Recall also that given two smooth families
$ A_{t} $, $ t\in T $, and $ B_{s} $, $ s\in S $, of linear operators acting in a vector space $ V $, the
product family $ A_{t}B_{s} $ is {\em not necessarily\/}\footnote{It would be smooth if the multiplication mapping $ \operatorname{End}\left(V\right)^{2} \to \operatorname{End}\left(V\right) $ were
continuous (in an appropriate sense). However, it is not, with a few
exceptions, such as $ \dim  V<\infty $, or the case of Banach $ V $. Later, in Section~%
\ref{s0.55}, we discuss the actual mechanism of appearence of singularities;
however, right now it is enough to understand that there is no {\em a priori\/}
reason to assume $ A_{t}B_{s} $ to be smooth.} smooth on the whole parameter
space $ T\times S $. We say that these families have an {\em operator product expansions\/}
if the product is smooth outside of a hypersurface $ \Sigma\subset T\times S $; moreover, $ A_{t}B_{s} $
has {\em at most a pole\/} at $ \Sigma $; in other words, there is a number $ k\in{\mathbb Z} $
such that the product $ \Phi^{k}\left(t,s\right)A_{t}B_{s} $ allows a smooth continuation from $ T\times S\smallsetminus\Sigma $
to $ T\times S $; here $ \Phi $ is defined by $ \Sigma=\left\{\Phi\left(t,s\right)=0\right\} $.

\subsection{OPEs and linear algebra } The aim of this introduction is to analyze
the relationship of OPEs with infinite-dimensional geometry of loop
spaces in more
detail. The bulk of this paper develops a linear algebra framework in
which ``good'' families of linear operators in a vector space $ V $
{\em automatically\/} satisfy OPEs. It turns out that the notion of ``goodness''
can be expressed in terms of two compatible structures on $ V $: on the
notion of {\em smoothness\/} on $ V $ (i.e., a description of which mappings $ M \to V $
are ``smooth''; here $ M $ is an arbitrary manifold), and on a {\em lattice\/}\footnote{I.e., a collection of vector subspaces of $ V $ closed w.r.t.~finite sums and
intersections.} of
vector subspaces in $ V $ parameterized by a union of algebraic varieties.

The principal results of this paper claim that
\begin{enumerate}
\item
it makes sense to consider the set of functionals on (an appropriate
flavor of) the loop space; this set is a vector space (even an algebra)
$ V $;
\item
this vector space $ V $ carries a smoothness structure and a compatible
lattice of vector subspaces;
\item
it makes sense to consider smooth families of (partial\footnote{In other words, linear operators $ V_{0} \to V $, here $ V_{0}\subset V $.}) linear
operators in $ V $ which are compatible with two structures on $ V $ introduced
above;
\item
any two such families satisfy an OPE;
\item
the usual ``physical parlance'' lexicon and toolset of conformal
fields is {\em literally\/} applicable to the families of the above form;
\item
the ``standard'' families of linear operators appearing in the theory
of vertex operators act in a vector subspace of $ V $; they can be
``continuously'' extended to the families of partial operators in $ V $
compatible with two structures on $ V $;
\item
moreover, these ``standard vertex operators'' have a geometric nature,
corresponding to the action of certain groups in certain line bundles on
the loop space.
\end{enumerate}

The corresponding flavor of the loop space is the space $ {\mathcal L}^{r}X $ of {\em rational
loops}, i.e., of rational mappings $ Z \to X $. (Here we assume that $ Z $ is a
projective manifold.)

This introduction serves a special purpose: in the principal part of
this paper we mostly ignore the last part of the program above, and
introduce the vector space $ V $ in purely algebraic terms. Here, in the
introduction, we discuss in details the geometric realization of $ V $. Given
that the corresponding infinite-dimensional geometry is quite hairy,
sometimes the discussion becomes superficial, e.g., providing only
circumstantial evidence, parables and analogies. We are going to discuss
the the geometry in details in \cite{Zakh98Mer}.

\subsection{Finite dimensional algebraic geometry }\label{s00.10}\myLabel{s00.10}\relax  The geometry of an affine
finite-dimensional algebraic variety $ M $ (possibly equipped with an action
of a group $ G $) is encoded in commutator rules for geometric linear
operators (defined in Section~\ref{s0.01}) in $ \operatorname{Func}\left(M\right) $. First of all, the
commutative algebra $ A=\operatorname{Func}\left(M\right) $ of functions on $ M $ (the corresponding
commutator rule is $ aa'=a'a $ for $ a,a'\in A $) determines $ M $. Moreover, given a
group $ G $ of symmetries of $ M $, the algebra $ \widetilde{A}=k\left[G\right]\ltimes A $ acts on $ \operatorname{Func}\left(M\right) $
by geometric operators. The structure of the group $ G $ is encoded into the
generators and relations of $ k\left[G\right] $; the structure of the action is encoded
in relationships of the form $ ga=a'g $; here $ a\in A $ and $ g\in G $ are given, and $ a' $
is an appropriate element of $ A $. In other words, given a pair $ \widetilde{A}\supset A $ of
algebras with commutation relations as described above, one can
reconstruct $ M $ and $ G $.

In a non-affine situation an extra ingredient is added; the {\em locality\/}
along $ M $ enters the play via consideration of locally-defined functions.
The action of $ G $ mixes the domains of definition; the situation may be
simplified again in the quasiprojective case, when one can replace $ A=\operatorname{Func}\left(M\right) $
with the algebra $ A=\operatorname{MerFunc}\left(M\right) $ of meromorphic (alternatively named
{\em rational\/}) functions on $ M $. After such a translation $ G $ acts again on
the algebra $ A $.

Suppose that the action $ \pi $ of $ G $ on $ M $ is extended to an action $ \Pi $ on a
line bundle $ {\mathcal L} $ over $ M $. Trivialize $ {\mathcal L} $ over a Zariski open subset $ U\subset M $; this
associates a meromorphic function $ a_{\sigma} $ on $ U $ to any locally defined section
$ \sigma $ of $ {\mathcal L} $. Now the action $ \Pi $ of $ G $ on locally defined (or meromorphic)
sections of $ {\mathcal L} $ is a geometric operator, $ \Pi_{g}a=\mu_{g}\cdot\pi_{g}^{*}a $; here $ \mu_{g} $ is a
meromorphic function on $ M $ (a {\em Schur multiplier\/}) describing how the action
of $ G $ on $ {\mathcal L} $ changes the trivialization of $ {\mathcal L} $.

In other words, the action of $ g\in G $ on the vector space of sections of
$ {\mathcal L} $ can be expressed in terms of operators of two types: the action of
diffeomorphisms of $ M $ on (meromorphic) functions of $ M $, and the operator of
multiplication by a (meromorphic) function on $ M $. In the introduction we
discuss analogues of these operators when $ M $ is a loop space.

\subsection{Schemes and loop spaces }\label{s0.04}\myLabel{s0.04}\relax  Note that the appearance of meromorphic
functions on $ M $ is an artifact of our usage of a trivialization of $ {\mathcal L} $.
However, in the case of loop spaces this is going to be an immanent
feature: it turns out that the ``manifolds'' $ M={\mathcal L}^{r}X $ (consisting of {\em rational\/}
loops on $ X $) has no regular functions, even locally defined regular
functions:

\begin{example} \label{ex00.40}\myLabel{ex00.40}\relax  Suppose $ Z={\mathbb P}^{1} $, $ X={\mathbb A}^{1} $. Let $ R_{k}\subset{\mathcal L}^{r}X $ consists of rational
loops of degree $ k $. One can parameterize $ R_{1} $ via $ \varphi_{P}\left(z\right)=\frac{az+b}{cz+d} $; here
$ P=\left(a:b:c:d\right)\in{\mathbb P}^{3} $. On the other hand, the loops $ \varphi_{P} $ with $ P= \left(\alpha\beta:\alpha'\beta:\alpha\beta':\alpha'\beta'\right) $
are just constant loops (from $ R_{0} $) with the value $ \left(\beta:\beta'\right)\in{\mathbb P}^{1} $. In other
words, $ R_{1} $ is the image of $ {\mathbb P}^{3} $ with all the projective lines
$ \left(\alpha B:\alpha'B:\alpha B':\alpha'B'\right) $ shrunk to a point (here $ B,B' $ are given numbers, and
$ \left(\alpha:\alpha'\right)\in{\mathbb P}^{1} $ is a parameter on the line).

On the other hand, it is obvious that $ {\mathbb P}^{3} $ with a line $ l $ shrunk to a
point has few locally defined functions; indeed, let $ L\in Q $ be the image of
the line $ l $ on the quotient space $ Q $; consider a neighborhood $ U $ of $ L $. Its
preimage $ U' $ on $ {\mathbb P}^{3} $ contains $ l $, thus it contains any projective line $ l' $
which is close to $ l $; in particular, any regular function on $ U' $ is
constant, since its pull back must be constant on any (projective!) line
$ l{}' $.

Similarly, $ R_{1} $ can have no regular functions in the neighborhood of
$ R_{0}\subset R_{1} $; same for a neighborhood of $ R_{k} $ in $ R_{k+1} $. Since any point of $ {\mathcal L}^{r}X $ is
in $ R_{k} $ for a large enough $ k $, this shows that $ {\mathcal L}^{r}X $ can have no regular
function. \end{example}

\begin{remark} \label{rem00.44}\myLabel{rem00.44}\relax  It turns out that the principal reason for this
phenomenon is that in the example we {\em allow a ``moving singularity''\/}: the
position of the pole of the loop $ \varphi_{P}\left(z\right) $ depends on the parameter $ P $;
however, we assume that the family $ \left\{\varphi_{P}\right\}_{P\in{\mathbb P}^{3}} $ is a ``smooth family'' in $ {\mathcal L}^{r}X $.
In turn, the reason for such an assumption (and the reason to consider
$ {\mathcal L}^{r}X! $) is the observation (see Sections~\ref{s0.05} and~\ref{s0.06}) that the
``classical vertex operators'' $ \exp \sum_{n>0}\frac{s^{-n}}{n}\partial_{n} $ and $ \sum_{n\geq0}s^{-n-1}\partial_{n} $ ``create'' a
pole of a loop at a moving point (actually at the point $ s\in Z={\mathbb P}^{1} $). \end{remark}

The example above shows that it does not make a lot of sense to
describe the ``space'' $ {\mathcal L}^{r}X $ (as well as the action of vertex operators)
using the language of schemes; indeed, this language is based on the idea
that the locally defined functions separate points. This may explain why
the scheme-theoretic approaches of \cite{BeiDriChiral} and \cite{KapVas01Ver} to
vertex algebras and loop spaces lead to so complicated constructions.

\subsection{Geometry of classical vertex operators and moving poles}\label{s0.05}\myLabel{s0.05}\relax  In a lot of
applications of families of linear operators with OPEs, the families are
given by explicit formulae involving the following series depending on
parameters $ t $ and $ s $:
\begin{gather} \exp \sum_{n\geq0}t^{n}{\mathbit y}_{n},\qquad \exp \sum_{n>0}\frac{s^{-n}}{n}\partial_{n}
\notag\\
\sum_{n\geq0}t^{n}{\mathbit x}_{n},\qquad \sum_{n\in{\mathbb Z}}t^{n}{\mathbit x}_{n},\qquad \sum_{n\geq0}s^{-n-1}\partial_{n},\qquad \sum_{n\in{\mathbb Z}}s^{-n-1}\partial_{n} 
\notag\end{gather}
Here $ x_{n} $ (or $ y_{n} $) are free generators of $ A={\mathbb C}\left[x_{n}\right]_{n\in{\mathbb Z}} $ or $ A={\mathbb C}\left[y_{n}\right]_{n\geq0} $
correspondingly, and $ {\mathbit x}_{n} $ (or $ {\mathbit y}_{n} $) are operators in $ A $ of multiplication by
the corresponding symbol $ x_{n} $ (or $ y_{n} $); $ \partial_{n} $ is the operator of taking the
derivative w.r.t.~the corresponding symbol $ x_{n} $ or $ y_{n} $. To show how these
operators are related to loop spaces and moving singularities of
meromorphic loops, consider $ A $ as an algebra of functions on a particular
loop space of the additive (for operators with $ x_{n} $) or multiplicative (for
operators with $ y_{n} $) group.

The flavor of a loop space on which this relation is simplest to
illustrate corresponds to $ Z $ being a punctured disk,
$ Z=\operatorname{PDisk}\left(R\right)\buildrel{\text{def}}\over{=}\left\{0<|z|<R\right\}\subset{\mathbb C} $. First, consider the additive group $ X={\mathbb A}^{1} $. Then
any loop $ \varphi\in{\mathcal L}X $ can be written as $ \varphi\left(z\right)=\sum_{n\in{\mathbb Z}}x_{n}z^{n} $; this gives a coordinate
system $ \left(x_{n}\right) $ on $ {\mathcal L}X $. Consider any ``reasonable'' definition of $ \operatorname{Func}\left({\mathcal L}X\right) $; then
the evaluation functional $ \operatorname{Eval}_{t}\colon \varphi \mapsto \varphi\left(t\right) $ must be an element of $ \operatorname{Func}\left({\mathcal L}X\right) $
for $ |t|<R $. Denote by $ {\mathbit E}_{t} $ the operator in $ \operatorname{Func}\left({\mathcal L}X\right) $ of multiplication by
$ \operatorname{Eval}_{t} $. Similarly, consider the mapping
\begin{equation}
G_{s}^{\left(\varepsilon\right)}\colon {\mathcal L}{\mathbb A}^{1} \to {\mathcal L}{\mathbb A}^{1}\colon \varphi \to \varphi',\qquad \varphi'\left(z\right)=\varphi\left(z\right)+\frac{\varepsilon}{z-s}.
\notag\end{equation}
These mappings form a well-defined one-parameter group of automorphisms
of $ {\mathcal L}{\mathbb A}^{1} $ for $ |s|>R $; denote by $ {\mathbit d}_{s} $ the corresponding vector field on $ {\mathcal L}{\mathbb A}^{1} $.
Again, for any reasonable definition of $ \operatorname{Func}\left({\mathcal L}{\mathbb A}^{1}\right) $, this vector field
induces a linear operator on $ \operatorname{Func}\left({\mathcal L}{\mathbb A}^{1}\right) $; denote it by the same symbol.

In the coordinates $ x_{n} $ one can write the function $ \operatorname{Eval}_{t} $ on $ {\mathcal L}{\mathbb A}^{1} $ as
$ \sum_{n\in{\mathbb Z}}t^{n}x_{n} $. Moreover, since the vector field $ {\mathbit d}_{s} $ is translation-invariant on
$ {\mathcal L}{\mathbb A}^{1} $, it must be a linear combination of the operators $ \partial_{n} $. Writing the
Taylor coefficients of $ 1/\left(z-s\right) $, one obtains

\begin{theorem} \label{th02.10}\myLabel{th02.10}\relax  The operators $ {\mathbit E}_{t} $ and $ {\mathbit d}_{s} $ can be written as
\begin{equation}
{\mathbit E}_{t}= \sum_{n\in{\mathbb Z}}t^{n}{\mathbit x}_{n},\qquad {\mathbit d}_{s}=-\sum_{n\geq0}s^{-n-1}\partial_{n}.
\notag\end{equation}
here $ {\mathbit E}_{t} $ is a well-defined\footnote{Since these sums are infinite, one needs to describe the topology on
$ \operatorname{Func}\left({\mathcal L}_{\operatorname{PDisk}\left(R\right)}{\mathbb A}^{1}\right) $ first. However, this statement holds with all the
natural choices of topology on this space, so we skip this issue now.} linear operator in $ \operatorname{Func}\left({\mathcal L}_{\operatorname{PDisk}\left(R\right)}{\mathbb A}^{1}\right) $ for
$ 0<|t|<R $; $ {\mathbit d}_{s} $ is a well-defined linear operator in $ \operatorname{Func}\left({\mathcal L}_{\operatorname{PDisk}\left(R\right)}{\mathbb A}^{1}\right) $ for
$ |s|\geq R $. \end{theorem}

One can also consider a non-punctured disk $ \operatorname{Disk}\left(R\right) $ in similar terms.
In this case the coordinate systems is $ \left(x_{n}\right)_{n\geq0} $, and the sum for the
operator $ {\mathbit E}_{t} $ has only the terms for $ n\geq0 $.

Now consider a similar loop space $ \left\{\varphi\colon \operatorname{Disk}\left(R\right) \to {\mathbb C}^{*}\right\} $ of the
multiplicative group. Consider a coordinate system $ \left(y_{n}\right)_{n\geq0} $ on $ {\mathcal L}{\mathbb C}^{*} $ by
writing $ \varphi\left(z\right)=y_{0}\exp \sum_{n>0}y_{n}z^{n} $. Abusing notations, denote $ \partial_{y_{n}} $ as $ \partial_{n} $.

As above, consider the operator $ {\mathbit E}_{t} $; consider also the mapping
\begin{equation}
H_{s}\colon {\mathcal L}{\mathbb C}^{*} \to {\mathcal L}{\mathbb C}^{*}\colon \varphi \to \varphi',\qquad \varphi'\left(z\right)=\frac{s\varphi\left(z\right)}{s-z}.
\notag\end{equation}
This mapping is well-defined for $ |s|>R $; denote by $ {\mathbit H}_{s} $ the corresponding
operator on $ \operatorname{Func}\left({\mathcal L}{\mathbb C}^{*}\right) $. Again, since $ H_{s} $ preserves $ y_{0} $ and is a parallel
translation in coordinates $ y_{n} $, $ n>0 $, one can write $ {\mathbit H}_{s} $ as $ \exp \left(v\right) $; here $ v $ is
a translation-invariant vector field in coordinates $ y_{n} $. The calculation
of the Taylor series for $ \log \frac{s}{s-z} $ gives

\begin{theorem} \label{th02.20}\myLabel{th02.20}\relax  The operators $ {\mathbit E}_{t} $ and $ {\mathbit H}_{s} $ can be written as
\begin{equation}
{\mathbit E}_{t}={\mathbit y}_{0}\exp \sum_{n>0}t^{n}{\mathbit y}_{n},\qquad {\mathbit H}_{s}=\exp \sum_{n>0}\frac{s^{-n}}{n}\partial_{n};
\label{equ0.5.10}\end{equation}\myLabel{equ0.5.10,}\relax 
here $ {\mathbit E}_{t} $ is a well-defined linear operator in $ \operatorname{Func}\left({\mathcal L}_{\operatorname{Disk}\left(R\right)}{\mathbb C}^{*}\right) $ for $ |t|<R $;
$ {\mathbit H}_{s} $ is a well-defined linear operator in $ \operatorname{Func}\left({\mathcal L}_{\operatorname{Disk}\left(R\right)}{\mathbb C}^{*}\right) $ for $ |s|\geq R $. \end{theorem}

\begin{remark} For the mapping $ H_{s} $ the point $ 0\in\operatorname{Disk}\left(R\right) $ plays a special role.
However, one can consider the mapping $ h_{s}\colon \varphi \mapsto \widetilde{\varphi} $, $ \widetilde{\varphi}\left(z\right)=\varphi\left(z\right)/\left(z-s\right) $,
instead of $ H_{s} $. The corresponding linear operator $ {\mathbit h}_{s} $ is just $ \left(-s\right)^{y_{0}\partial_{0}}{\mathbit H}_{s} $;
hence it can be written by a formula similar to (0.5.10). \end{remark}

\begin{remark} It is clear that the operators above should be written by the
same formulae no matter which flavor of $ {\mathcal L}X $ (here $ X={\mathbb C} $ or $ X={\mathbb C}^{*} $)
and which flavor of the space $ \operatorname{Func}\left({\mathcal L}X\right) $ of functional on $ {\mathcal L}X $ one considers.
Indeed, it is enough to check the equality of the right-hand side and the
left-hand side on a ``dense subset'' of $ \operatorname{Func}\left({\mathcal L}X\right) $, and this dense subset can
be chosen independently of the flavor. \end{remark}

\subsection{OPEs for classical vertex operators }\label{s0.06}\myLabel{s0.06}\relax  Consider a manifold $ M $;
denote by $ \mu_{F} $ the operator in $ \operatorname{Func}\left(M\right) $ of multiplication by a function
$ F\in\operatorname{Func}\left(M\right) $; consider the operator $ {\mathbit V} $ in $ \operatorname{Func}\left(M\right) $ induced by a vector field $ V $
on $ M $. Then, obviously, $ \left[{\mathbit V},\mu_{F}\right]=\mu_{V\cdot F} $. Hence the operators $ {\mathbit d}_{s} $ and $ {\mathbit E}_{t} $ have
the commutation relation $ \left[{\mathbit d}_{s},{\mathbit E}_{t}\right]=\frac{1}{s-t}{\mathbit E}_{t} $. This exposes the first
OPE-like relationship of this paper. Similarly, $ {\mathbit h}_{s}{\mathbit E}_{t}=\frac{1}{t-s}{\mathbit E}_{t}{\mathbit h}_{s} $. Since
there expressions are not exactly of the OPE form, comment on their
relations to OPEs.

In the former case the right-hand side depends on $ s $ via the
denominator only; thus the left-hand side {\em must\/} have a pole at $ s=t $. In the
latter case $ {\mathbit E}_{t}{\mathbit h}_{s} $ sends the constant functional 1 on $ {\mathcal L}X $ to $ \operatorname{Eval}_{t} $; hence
$ {\mathbit E}_{t}{\mathbit h}_{s} $ has no zero at $ t=s $. Thus, again, $ {\mathbit h}_{s}{\mathbit E}_{t} $ {\em must\/} have a singularity of an
OPE type when $ t=s $.

However, these poles are impossible to see when one considers the
vector space of functionals on analytic flavors of loops, such as
$ \operatorname{Func}\left({\mathcal L}_{\operatorname{Disk}\left(R\right)}{\mathbb C}^{*}\right) $. Indeed, the domains of definition $ \left\{|s|>R\right\} $ and $ \left\{|t|<R\right\} $
of the families in question do not intersect, thus the divisor of poles
(the diagonal) does not intersect the product of domains of two families
(e.g., $ {\mathbit E}_{t} $ and $ {\mathbit d}_{s} $). Moreover, in the discussion above we did not
completely specify the nature of families $ {\mathbit E}_{t} $ etc.: all we defined was
$ {\mathcal L}_{\operatorname{PDisk}\left(R\right)}{\mathbb A}^{1} $ etc., while we also need to specify $ \operatorname{Func}\left({\mathcal L}_{\operatorname{PDisk}\left(R\right)}{\mathbb A}^{1}\right) $ and
which families of operators in $ {\mathcal L}_{\operatorname{PDisk}\left(R\right)}{\mathbb A}^{1} $ are smooth families. Until we
define these notions, we can only discuss the properties of commutators,
$ \left[{\mathbit d}_{s},{\mathbit E}_{t}\right] $ and $ {\mathbit h}_{s}{\mathbit E}_{t}{\mathbit h}_{s}^{-1} $, not the properties of products of elements of the
family, as the OPEs must.

Indeed, the crucial operation in the theory of vertex operators
is the calculation of the {\em fusions\/}; fusions are the Laurent coefficients of
the composition $ A_{s}\circ B_{t} $ near the divisor of poles. For a lot of
calculations it is convenient to write these Laurent coefficients as the
corresponding Cauchy integrals. To be able to do this, one needs to have
$ A_{s}\circ B_{t} $ defined in a neighborhood of the divisor of poles; it is not enough
to have it defined in ``one half-space'' near this divisor; however, this
is the situation with the domain $ \left\{\left(s,t\right) \mid |s| > R, |t| < R\right\} $.

Later, when we switch to the different flavor $ {\mathcal L}X={\mathcal L}^{r}X $ of loop spaces,
we will be able to define the necessary notions (the vector space
$ \operatorname{Func}\left({\mathcal L}X\right) $ and the notion of smooth families of operators on $ \operatorname{Func}\left({\mathcal L}X\right) $) in
such a way that the domains of the families $ {\mathbit E}_{t} $ and $ {\mathbit h}_{s} $ will intersect. We
will see, for example, that for this flavor $ {\mathbit E}_{t}\circ{\mathbit h}_{s} $ is a smooth family; thus
$ {\mathbit h}_{s}\circ{\mathbit E}_{t} $ is going to have a pole of the first order when $ s=t $.

Anyway, a theory of loop spaces which allows the consideration of
OPEs for classical vertex operators must allow the consideration of the
dependence of the family of mappings $ h_{s} $ on the parameter $ s $, similarly for
the vector fields $ {\mathbit d}_{s} $ and mappings $ G_{s}^{\left(\varepsilon\right)} $. On the other hand, application
of the group generated by mappings $ h_{s} $ (or by mappings $ G_{s}^{\left(\varepsilon\right)} $) to an
arbitrary loop allows a creation of a family of loops with arbitrary
{\em moving poles}. This justifies the parameterizations of Example~\ref{ex00.40}
via Remark~\ref{rem00.44}.

\subsection{Locality condition and singularities } As we have seen it above, a
simultaneous consideration of ``basic'' functionals on loop spaces and
``basic'' diffeomorphisms of loop spaces leads to OPEs and OPE-like
singularities. We also saw that the calculus of classical vertex
operators is equivalent to the calculus involving such basic functionals
and basic transformations.

It turns out that there is another reason to consider OPE-like
singularities; it is connected with another crucial feature of loop
spaces: a special role played by ``local functionals''; this phenomenon is
completely lacking in the situation $ \dim <\infty $. In addition to ``locality'' on
$ {\mathcal L}X=\left\{Z \to X\right\} $ (e.g., consideration of open subsets of $ {\mathcal L}X $) one must consider
{\em locality along a loop\/} (i.e., consideration of open subsets of $ Z $). As an
illustration, consider the language of calculus of variations: when one
considers a functional $ F\left(\varphi\right) $ on the space of loops $ Z \xrightarrow[]{\varphi} X $, the partial
derivatives of $ F $ are written as $ \delta F/\delta\varphi\left(z\right) $; they depend on a {\em continuous\/}
index $ z $; it is crucial to analyze the degree of smoothness of the
dependence of these data on $ z $.

A functional of the form $ F\left(\varphi\right)=\int_{Z}L\left(z,\varphi\left(z\right),\varphi'\left(z\right),\dots ,\varphi^{\left(n\right)}\left(z\right)\right)dz $ (here $ \varphi $
is a loop $ \varphi\left(z\right) $) is called a {\em local functional\/} of $ \varphi $ (here $ L $ is an arbitrary
fixed function). When we consider two local variation of the loop $ \varphi $, one
near $ z_{1}\in Z $, another near $ z_{2} $, $ z_{2}\not=z_{1} $, the changes of $ F $ w.r.t.~these two
variations are additive; in other words, the mixed variational derivative
$ \Delta\left(z,z'\right)=\frac{\delta^{2}F}{\delta\varphi\left(z\right)\delta\varphi\left(z'\right)} $ vanishes outside of the diagonal $ \left\{z=z'\right\} $: for an
appropriate $ N $ one has $ \left(z-z'\right)^{N}\Delta\left(z,z'\right)=0 $; here $ \Delta\left(z,z'\right) $ is a generalized
function on $ Z^{2} $.

The appearance of $ \delta $-function singularity of $ \Delta\left(z,z'\right) $ on the diagonal
is very similar to appearances of OPEs in the theory of loop spaces. This
singularity is not a pole, however, as we will see it in the next
section, there is a standard technique to construct related objects with
poles.

\subsection{Skyscraper variations }\label{s00.33}\myLabel{s00.33}\relax  Here we construct a bridge between the
construction of the previous section, which lead to $ \delta $-function
singularities, and the constructions of Section~\ref{s0.05}, which
lead to singularities of the OPE type.

In the language of variational calculus one analyses the dependence
of the functional on the ultimate local variation of a loop; it is $ g_{s}^{\left(\varepsilon\right)}:
\varphi \mapsto \varphi+\varepsilon\delta_{s} $; here $ \delta_{s} $ is the skyscraper $ \delta $-function with support at $ s\in Z $,
and, at least in the case $ X={\mathbb A}^{n} $, $ \varepsilon $ is a tangent vector to $ X $. Assume $ X={\mathbb A}^{1} $,
then $ \varepsilon $ is a number. Consider $ g_{s}^{\left(\varepsilon\right)} $ as a $ 1 $-parameter subgroup of $ \operatorname{Diff}\left({\mathcal L}X\right) $;
taking $ \varepsilon \to 0 $ gives a vector field $ \delta_{s} $ on $ {\mathcal L}X $. Define a linear operator $ {\mathbit D}_{s} $
on $ \operatorname{Func}\left({\mathcal L}X\right) $ as the action of $ \delta_{s} $.

We did not specify which flavor of the loop space $ {\mathcal L}X $ we work with,
but the most useful examples do not allow $ \delta $-like loops; thus one needs to
consider $ g_{s}^{\left(\varepsilon\right)} $ formally, as a generalized function. For any smooth
measure $ \Phi\left(z\right) $ on $ Z $, the integral $ \int_{Z}\Phi\left(z\right)\delta_{z} $ is a well-defined vector field on
$ {\mathcal L}X $, thus $ {\mathbit D}_{\left[\Phi\right]}\buildrel{\text{def}}\over{=}\int_{S^{1}}\Phi\left(z\right){\mathbit D}_{z} $ is a well-defined operator on $ \operatorname{Func}\left({\mathcal L}X\right) $. Using
this description, it is easy to verify that $ \left[{\mathbit D}_{s},{\mathbit E}_{t}\right]=\delta\left(s-t\right){\mathbit E}_{t} $ (here
both sides are considered as generalized functions of $ s $).

Similarly to Section~\ref{s0.05}, the consideration of the products $ {\mathbit D}_{s}{\mathbit E}_{t} $
and $ {\mathbit E}_{t}{\mathbit D}_{s} $ separately is much more delicate.

To bridge this discussion with one of Section~\ref{s0.05}, it is
convenient to consider the kernels $ \Phi\left(z\right) $ of a special form. Assume that
$ Z=S^{1} $; embed $ S^{1} $ in $ {\mathbb C} $ as $ \left\{|z|=1\right\} $; let $ \Phi_{s}=\frac{dz}{z-s} $; then $ {\mathbit d}_{s}={\mathbit D}_{\left[\Phi_{s}\right]} $; here $ {\mathbit d}_{s} $
was introduced in Section~\ref{s0.05}. On the other hand, the importance of
$ \delta $-function variations is related to them spanning the vector space of all
functions; however, the functions $ \Phi_{s} $ also span this vector space. This is
easiest to see via the relation
\begin{equation}
\delta\left(x-s\right)=\frac{1}{2\pi i}\bar{\partial}_{s}\Phi_{s};
\notag\end{equation}
in other words, at least formally, the family $ {\mathbit D}_{s} $ can be reconstructed as
\begin{equation}
{\mathbit D}_{s}=\frac{1}{2\pi i}\bar{\partial}_{s}{\mathbit d}_{s}.
\label{equ0.15}\end{equation}\myLabel{equ0.15,}\relax 
This formula allows one to replace considerations of $ {\mathbit D}_{s} $ by
considerations of $ {\mathbit d}_{s} $. Since $ {\mathbit d}_{s} $ is a bona fide vector field on $ {\mathcal L}X $ (at
least for $ |s|\not=1 $), it is much more convenient to work with. In what
follows we do not work with the family of vector fields $ {\mathbit D}_{s} $; however, make
one more remark about this family.

Our treatments of the classical vertex operators in Section~\ref{s0.05}
skipped one family out of six: $ \sum_{n\in{\mathbb Z}}s^{-n-1}\partial_{n} $ was not related to either one
of operators $ {\mathbit E}_{t} $, $ {\mathbit d}_{s} $ and $ {\mathbit h}_{s} $ on different flavors of loop spaces. However,
the usual ``physical'' interpretation of $ \sum_{n\in{\mathbb Z}}s^{-n-1}z^{n} $ is that it describes
the $ \delta $-function at the point $ s $. This leads to identification of
$ \sum_{n\in{\mathbb Z}}s^{-n-1}\partial_{n} $ with the operator $ {\mathbit D}_{s} $. We do not know how to relate this
heuristic argument with Equation~\eqref{equ0.15}.

\subsection{Birational geometry and partial operators }\label{s0.55}\myLabel{s0.55}\relax  As we saw it in
Section~\ref{s0.04}, the set $ {\mathcal L}^{r}X $ of rational loops is not a variety: it does
not makes sense to consider {\em regular\/} functions on $ {\mathcal L}^{r}X $. This is a crucial
phenomenon; one of the principal aims of this paper is to show that the
set $ {\mathcal L}^{r}X $ has a well-defined ``birational geometry'': it makes sense to
consider {\em meromorphic\/} functionals of a loop. During the introduction, we
call such sets ``{\em birational manifolds\/}''.

The principal feature\footnote{It is this feature which is the principal geometric reason for
appearance of OPE for geometric operators.} of the birational geometry of
infinite-dimensional manifolds is that there is no notion of ``finite
codimension''; show how such a notion would lead to a contradiction in the
case of $ {\mathcal L}{\mathbb P}^{1} $. Fix $ z\in Z $; the $ 0 $-divisor $ D_{0} $ of $ \operatorname{Eval}_{z} $ is given by the equation
$ \varphi\left(z\right)=0 $ on $ \varphi $; if there is a theory of codimension, this set {\em must\/} have
codimension 1. Since the inversion $ \iota $ of $ {\mathcal L}{\mathbb P}^{1} $ (induced by $ x \mapsto x^{-1} $, $ x\in{\mathbb P}^{1} $)
must be a diffeomorphism of $ {\mathcal L}{\mathbb P}^{1} $, the $ \infty $-divisor $ D_{\infty} $ of $ \operatorname{Eval}_{z} $ must be also
of codimension 1. However, the ``diffeomorphism'' $ h_{z} $ of Section~\ref{s0.05}
identifies $ {\mathcal L}{\mathbb P}^{1}\smallsetminus D_{0} $ with $ D_{\infty} $. Thus either we do not consider the invertible
mapping $ h_{z} $ as a diffeomorphism, or we have no viable notion of
codimension. Since the former kills the geometric description of vertex
operators of Section~\ref{s0.05}, we are stuck with the later.

In particular, a diffeomorphism of an infinite-dimensional
birational manifold $ M $ can send it ``into'' a ``divisor of poles'' of a
meromorphic function (or even into the ``indeterminacy subset'' of such a
function). Thus the action of a diffeomorphism on the meromorphic
functions on $ M $ is not always defined. This is one of the key ingredients
of our approach: the action of group $ G $ on $ M $ is encoded not in a family of
linear operators on $ \operatorname{Func}\left(M\right) $, but in a family of {\em partial linear operators\/}
on $ \operatorname{Func}\left(M\right) $ (here Func stands for the vector space of meromorphic
functions).

One of the aims of this paper is to demonstrate that rules of
composition of partial linear operators fully reflect the properties of
products in vertex operator algebras. Moreover, in appropriate vector
spaces any two smooth families of partial linear operators satisfy OPE
relations.

Note that products of two smooth families of partial operators in a
finite-dimensional vector space automatically possess properties {\em similar\/}
to OPE: for the simplest example consider a partial operator $ P $ and a
family $ Q_{t} $; then the dimension of the domain $ \operatorname{Dom}\left(R_{t}\right) $ of $ R_{t}\buildrel{\text{def}}\over{=}P\circ Q_{t} $ is the
same for generic $ t $, but has jumps where $ \operatorname{Image}\left(Q_{t}\right) $ take special position
w.r.t.~$ \operatorname{Dom}\left(P\right) $. In other words, the composition $ P\circ Q_{t} $ has singularities for
special values of $ t $; this property is similar to OPE.

In Section~\ref{s1.65} we formulate conditions on special families of
partial operators in an infinite-dimensional vector space (which we call
{\em codominant conformal fields\/}); one of the principal results of this paper
is that compositions of such families satisfy {\em exactly\/} the OPEs. This
gives the interpretation of the OPEs purely in the terms of the
properties of the {\em vector space\/} in which the families of the operators
act.

\subsection{Flavors of the theories }\label{s0.44}\myLabel{s0.44}\relax  If one needs more convincing examples of
conformal fields, one needs more precise definitions of the building
blocks: $ {\mathcal L}X $, $ \operatorname{Func}\left({\mathcal L}X\right) $, and of smooth families of operators in $ \operatorname{Func}\left({\mathcal L}X\right) $.
The examples of families $ {\mathbit E}_{t} $, $ {\mathbit d}_{s} $ and $ {\mathbit h}_{s} $ in Section~\ref{s0.05} show that it is
very convenient to be able to analytically extend the loops from $ S^{1}\subset{\mathbb C} $.

This is a basis of the {\em analytic approach}. It assumes that $ X $ and $ Z $
are complex-analytic; for an example, let $ Z=A_{\varepsilon}=\left\{1<|z|<1+\varepsilon\right\} $. This gives a
loop space $ {\mathcal L}_{\varepsilon}X $ consisting of loops $ \varphi\left(z\right) $ in $ X $, $ z\in A_{\varepsilon} $. Then the family $ {\mathbit E}_{t} $ is
well defined and analytic\footnote{After an appropriate definition of $ \operatorname{Func}\left({\mathcal L}_{\varepsilon}X\right) $.} for $ t\in A_{\varepsilon} $; the operators $ {\mathbit d}_{s} $ and $ {\mathbit h}_{s} $ are
well-defined and analytic for $ t\notin A_{\varepsilon} $. As we saw it already in Section~%
\ref{s0.06}, in this approach the composition $ {\mathbit h}_{s}{\mathbit E}_{t} $ is defined on too small a
subset of $ {\mathbb C}\times{\mathbb C} $ to define fusions via Cauchy integrals. Even if one allows
$ \varepsilon $ to vary, one can define $ {\mathbit h}_{s}{\mathbit E}_{t} $ only for $ |s|>|t| $. In physics-speak, this
leads to a {\em time-ordered\/} product---time being $ |z| $.

This difficulty, having the operator product well-defined only on
small subset of the desired set, is immanent to all approaches to the
loop spaces based on regular functions on $ {\mathcal L}X $. Another example of such
approach is the {\em formal approach\/}; in this approach $ {\mathcal L}^{f}X $ consists of {\em formal
loops\/} in $ X $, i.e., the mappings $ \varphi\in X\left(\left(z\right)\right) $ of a ``punctured formal disk'' into
$ X $; this set is well-defined if $ X={\mathbb A}^{m} $ (or $ X\subset{\mathbb A}^{m} $). A loop can be written as a
formal Laurent series $ \sum_{n\geq-N}x_{n}z^{n} $, $ x_{n}\in{\mathbb A}^{m} $. The loops with $ N=0 $ are
{\em holomorphic formal loops}, which form a subset $ {\mathcal L}_{+}^{f}X\subset{\mathcal L}^{f}X. $\footnote{The set $ {\mathcal L}_{+}^{f}X $ is well-defined for an arbitrary manifold $ X $; while it is
hard to define $ {\mathcal L}^{f}X $ in such a case, the formal neighborhood $ {\mathcal L}^{ff}X $ of $ {\mathcal L}_{+}^{f}X $
in $ {\mathcal L}_{+}^{f}X $ is well-defined. One can describe this neighborhood by the
condition that $ x_{n} $ with $ n<0 $ are confined to an infinitesimal neighborhood
of 0. See \cite{KapVas01Ver} for details.}

Since the punctured formal disk plays the role of the limit of
annuli $ \left\{r<|z|<r+\varepsilon\right\} $ when $ r,\varepsilon \to $ 0, the considerations for the analytic
case show that the family of operators $ {\mathbit h}_{s} $ is defined on a large set $ s\not=0 $.
However, the family $ {\mathbit E}_{t} $ is defined {\em``nowhere''\/}; it depends on the parameter
$ t $ as a formal Laurent series in $ t $. The coefficients of this formal series
are well-defined operators in $ \operatorname{Func}\left({\mathcal L}^{f}X\right) $. Thus to discuss the ``poles of
the composition'' $ {\mathbit h}_{s}{\mathbit E}_{t} $ one needs to invent an appropriate purely-algebraic
language.

The great advantage of the formal approach is that all the objects
are direct and/or inverse limits of finite-dimensional manifolds; thus
the task of definition of $ \operatorname{Func}\left({\mathcal L}^{f}X\right) $ becomes purely algebraic. E.g., the
vector space $ \operatorname{Func}\left({\mathcal L}^{f}X\right) $ consists of {\em polynomials\/} in variables $ x_{n} $ from
Section~\ref{s0.05} (see \cite{Kac97Ver} for details). Thus the formal approach
avoids the principal difficulty: that the infinite-dimensional
differential geometry is currently very far from being a satisfactory
theory (recall that---with the exception of Banach manifolds---infinite
dimensional analysis is still in its cradle).

Finally, in the {\em rational approach\/} of this paper, $ Z $ is a projective
manifold, $ X $ is a quasiprojective manifold, and $ {\mathcal L}^{r}X $ consists of rational
mappings $ Z \to X $; in other words, of mappings $ Z_{0} \to X $ with $ Z_{0} $ being dense
open in $ Z $; two mappings are the same if their restrictions to some open
dense subset $ Z_{00} $ coincide. As Example~\ref{ex00.40} shows, with reasonable
definitions, $ {\mathcal L}^{r}X $ has no non-constant regular function. However, $ {\mathcal L}^{r}X $ has
many ``rational functions'' (e.g., $ \operatorname{Eval}_{t} $, $ t\in Z $); it is the interplay of the
functions not being defined everywhere, and of diffeomorphisms of $ {\mathcal L}^{r}X $
(e.g., $ g_{t}^{\left(\varepsilon\right)} $) mixing ``divisors'' and ``complements to divisors'' which leads
to appearance of OPE in the algebraic geometry of $ {\mathcal L}^{r}X $.

The advantages of the rational approach are the following. First,
(same as in the analytic approach), the crucial families of operators
(such as $ {\mathbit h}_{s} $, $ {\mathbit E}_{t} $) are well-defined families depending on points $ s\in S $, $ t\in T $
(see Section~\ref{h12}). Second, the compositions (such as $ {\mathbit h}_{s}{\mathbit E}_{t} $) are
well-defined on large enough subsets of $ S\times T $ so that the fusions can be
defined via Cauchy integrals (see Section~\ref{h14}). Third, the crucial
definitions (of the set $ {\mathcal L}^{r}X $, of the vector space $ V=\operatorname{Func}\left({\mathcal L}^{r}X\right) $, and of
which families of operators in $ V $ are smooth) turn out to be purely
algebraic (see Sections~\ref{s11.10} and~\ref{s1.65}).

In particular, the analysis on the set $ {\mathcal L}^{r}X $ considered as an
``infinite dimensional smooth manifold'' does not require introduction of
topology on this set! Out of the data encoded into the definition of a
smooth manifold, we use only the notion of ``smooth families''; in other
words, which mappings $ M \to {\mathcal L}^{r}X $ are smooth (for finite-dimensional
manifolds $ M $). In other words, essentially we consider $ {\mathcal L}^{r}X $ as a functor
Manifolds $ \xrightarrow[]{\text{Families}} $ Sets; the set $ {\mathcal L}^{r}X $ is reconstructed as the
particular value $ \operatorname{Families}\left(M_{0}\right) $ of the functor when $ M_{0} $ is a point.

Note the surprising analogy of this approach with the so-called
``convenient settings'' \cite{KrieMich97Con} in analysis. In the latter
approach one {\em starts\/} with a topology of a would-be-manifold $ L $, and defines
which mappings $ M \to L $ (with $ \dim  M<\infty $) are smooth via an atlas on $ L $;
however, after this initial step, one more or less forgets about the
initial topology, constructing all data basing on this ``structure of the
functor'' on $ Y $.

\subsection{On $ \protect \operatorname{Func}\left({\protect \mathcal L}^{r}X\right) $ }\label{s0.77}\myLabel{s0.77}\relax  In this section we give a heuristic sketch of
results of \cite{Zakh98Mer} which reduces the birational geometry of $ {\mathcal L}^{r}X $ to
purely algebraic terms. We motivate the definition of $ \operatorname{Func}\left({\mathcal L}^{r}X\right) $ from
Section~\ref{s11.10} via more detailed consideration of the family $ \operatorname{Eval}_{s} $, $ s\in Z $,
of functionals on $ {\mathcal L}^{r}X $.

First of all, in which sense is the functional $ \operatorname{Eval}_{s}\colon \varphi \mapsto \varphi\left(s\right) $ a
meromorphic functional of $ \varphi? $ Answer: the restriction of this functional
to ``natural'' finite-dimensional families in $ {\mathcal L}^{r}X $ is meromorphic. List some
of such families.

Fix an ample line bundle $ {\mathfrak L} $ on $ Z $. Recall that one can define {\em the
degree\/} of a loop $ \varphi\in{\mathcal L}^{r}{\mathbb A}^{1} $, $ Z \xrightarrow[]{\varphi} {\mathbb A}^{1} $, as the minimum of $ k $ such that $ \varphi=S/S' $
with $ S,S'\in\Gamma\left(Z,{\mathfrak L}^{k}\right) $. One can proceed similarly for $ {\mathcal L}^{r}X $ for other
quasiprojective manifolds $ X $. Now $ {\mathcal L}^{r}X $ is a union of subsets $ R_{k} $ consisting
of loops of
degree $ k $ or less. For large enough $ k $, $ \operatorname{Eval}_{s} $ is defined on $ R_{k} $ (for small $ k $
it may happen that $ R_{k} $ is in the subset of indeterminacy of $ \operatorname{Eval}_{s} $, if $ s $ is
in the base set of the linear system for $ {\mathfrak L} $). Moreover, one can
parameterize $ R_{k} $ via the mapping $ \widetilde{R}_{k} \buildrel{\text{def}}\over{=}{\mathbb P}\left(\Gamma\left(Z,{\mathfrak L}^{k}\right)^{2}\right) \to R_{k}\colon \left(S,S'\right) \mapsto S/S' $.
This is a bijection on a dense open subset of $ \widetilde{R}_{k} $. Clearly, $ \operatorname{Eval}_{s} $
corresponds to a meromorphic function $ \varepsilon_{s,k} $ on $ \widetilde{R}_{k} $; one could call $ \operatorname{Eval}_{s} $ a
``pro-meromorphic'' function. This leads to the following:

\begin{definition} A {\em meromorphic functional\/} $ F $ on $ {\mathcal L}^{r}X $ is a collection of
meromorphic functions $ F_{k} $ on $ \widetilde{R}_{k} $ ``compatible'' with natural inclusions $ R_{k} \to
R_{k+1} $. \end{definition}

The problem of this definition is that it combines two objects: $ R_{k} $
and $ \widetilde{R}_{k} $; it turns out that the ``correct'' notion of ``compatibility'' is
delicate.\footnote{It is this detail which leads to the fact that diffeomorphism only act
as partial linear operators on meromorphic functions (thus to OPEs);
hence the complication we discuss in this section is actually a central
point of the geometry of $ {\mathcal L}^{r}X $.} Indeed, $ F $ is not everywhere defined on $ R_{k+1} $.

Instead of the mappings $ R_{k} \to R_{k+l} $, consider their liftings $ \widetilde{R}_{k} \to
\widetilde{R}_{k+l} $. The lifting is not uniquely defined; for any $ T\in\Gamma\left(Z,{\mathfrak L}^{l}\right) $, the mapping
$ \left(S,S'\right) \buildrel{\alpha_{T}}\over{\mapsto} \left(TS,TS'\right) $ is a lifting. Thus the set of liftings is
parameterized by $ {\mathbb P}\Gamma\left(Z,{\mathfrak L}^{l}\right) $. A particular mapping $ \alpha_{T} $ may\footnote{As already happens with the functional $ \operatorname{Eval}_{s} $.} send $ \widetilde{R}_{k} $ into the
indeterminacy set of $ F_{k+l} $. Thus the compatibility condition should be
formulated as this: for {\em a generic\/} $ T\in{\mathbb P}\Gamma\left(Z,{\mathfrak L}^{l}\right) $ the set $ \alpha_{T}\widetilde{R}_{k}\subset\widetilde{R}_{k+l} $ is not
contained in the indeterminacy set of the meromorphic function $ F_{k+l} $ on
$ \widetilde{R}_{k+l} $; moreover,\footnote{Note that this condition applied to a locally defined {\em regular\/} functions
on $ \widetilde{R}_{k} $, $ k\geq0 $, can be satisfied by constant functions only. See Example~%
\ref{ex00.40} for details.} $ \alpha_{T}^{*}F_{k+l}=F_{k} $. It turns out that this condition is not
enough for functoriality: one needs also require that $ \alpha_{T}^{-1}{\mathcal I}_{F_{k+l}}={\mathcal I}_{F_{k}} $; here
$ {\mathcal I}_{G} $ denotes the indeterminacy set of a meromorphic function $ G $.

It is obvious that the defined above functions $ \varepsilon_{s,k} $ (induced by
$ \operatorname{Eval}_{s} $) satisfy these compatibility conditions. Similarly to $ \operatorname{Eval}_{s}\colon \varphi \mapsto
\varphi\left(s\right) $, one can send $ \varphi $ to the $ K $-th partial derivative of $ \varphi $ at $ s $; here $ K $ is
a multiindex. This gives another meromorphic function $ \operatorname{Eval}_{s,K} $ on $ {\mathcal L}^{r}X $; as
one can easily check, functions $ \operatorname{Eval}_{t,K} $ are algebraically independent.
Since $ \operatorname{MerFunc}\left({\mathcal L}^{r}X\right) $ is a field, one can consider its subfield freely
generated by $ \operatorname{Eval}_{t,K} $. Call the latter field $ \operatorname{NaiveFunc}\left({\mathcal L}^{r}X\right) $.

One of the principal conclusions of \cite{Zakh98Mer} is that in the case
$ Z={\mathbb P}^{1} $ one has $ \operatorname{NaiveFunc}\left({\mathcal L}^{r}X\right)=\operatorname{MerFunc}\left({\mathcal L}^{r}X\right) $. One can conjecture that a
similar identity holds for an arbitrary $ Z $. No matter what is the fortune
of this conjecture, in this paper we {\em define\/} a ``meromorphic function'' on
$ {\mathcal L}^{r}X $ as an element of $ \operatorname{NaiveFunc}\left({\mathcal L}^{r}X\right) $. The motivation for such a definition
is that it quickly leads to a rich theory with all the manifestations of
OPEs present from the start.

\subsection{Partial operators in $ \protect \operatorname{NaiveFunc}\left({\protect \mathcal L}^{r}X\right) $ } Now, when we finally defined the
vector space $ \operatorname{NaiveFunc}\left({\mathcal L}^{r}X\right) $, we can analyze the operators $ {\mathbit E}_{t} $, $ {\mathbit d}_{s} $ and $ {\mathbit h}_{s} $
in this space. First of all, the operator $ {\mathbit E}_{t} $ are bona fide linear
operators in $ \operatorname{NaiveFunc}\left({\mathcal L}^{r}X\right) $. The operator $ {\mathbit E}_{t} $ multiplies a rational
expression of symbols $ \operatorname{Eval}_{z,K} $ (as any element of $ \operatorname{NaiveFunc}\left({\mathcal L}^{r}X\right) $ is) by
the symbol $ \operatorname{Eval}_{t} $. The result is a rational expression of the same form.

On the other hand, the operators $ {\mathbit d}_{s} $ and $ {\mathbit h}_{s} $ are of a different
nature. E.g., obviously $ \operatorname{Eval}_{t}\circ h_{s}=\frac{1}{t-s}\operatorname{Eval}_{t} $. This can be taken as a
definition of the action of $ {\mathbit h}_{s} $ on $ \operatorname{NaiveFunc}\left({\mathcal L}^{r}X\right)\colon {\mathbit h}_{s}\left(\operatorname{Eval}_{t}\right)=\frac{1}{t-s}\operatorname{Eval}_{t} $,
with obvious modifications (given by Leibniz identity) for the action on
$ \operatorname{Eval}_{t,K} $ (e.g., $ {\mathbit h}_{s}\left(\operatorname{Eval}_{t,1}\right)=\frac{1}{t-s}\operatorname{Eval}_{t,1}-\frac{1}{\left(t-s\right)^{2}}\operatorname{Eval}_{t,0} $), and
extending to rational expressions as an automorphism of an algebra.
(Similarly for the action of $ {\mathbit d}_{s} $, which acts as a derivation in the
algebra.)

However, the formula for $ {\mathbit h}_{s} $ does not make sense when $ t=s $. Instead of
trying to redefine this formula in this case, accept this as a fact: the
operators $ {\mathbit h}_{s} $ and $ {\mathbit d}_{s} $ are only {\em partial operators\/} (i.e., operators defined
on a vector subspace only). We already discussed the geometric reasons
for this in Section~\ref{s0.55}; we are going to consider this in details in a
foregoing paper \cite{Zakh98Mer}.

A closely related statement from geometry of $ {\mathcal L}^{r}X $ is that $ {\mathcal L}^{r}X $ has no
{\em regular functions}, see Example~\ref{ex00.40}. Provide a different
illustration of this phenomenon:

\begin{example} Consider a constant loop $ \varphi\left(z\right)\equiv 0\in R_{0}\subset{\mathcal L}^{r}X $. Then $ \operatorname{Eval}_{0}\left(\varphi\right)=0 $.
Include $ \varphi=\varphi_{0} $ into a one-parameter family $ \varphi_{s}\left(z\right)=\frac{s}{s-z} $. Then
$ \operatorname{Eval}_{0}\left(\varphi_{s}\right)\equiv 1 $; thus the value of $ \operatorname{Eval}_{0} $ on $ \varphi_{0} $ is not compatible with the
value of $ \operatorname{Eval}_{0} $ on the family $ \left\{\varphi_{s}\right\} $. To restore the compatibility, we need
two modifications: first, we need to include $ \varphi $ into a larger family
(e.g., $ 3 $-parameter parameterizing $ R_{1} $; for simplicity, consider a smaller
family
$ \varphi_{s,r}=\frac{s}{r-z} $); second, out of many representatives $ \varphi_{0,r} $ for $ \varphi $ in this
larger family we need to choose those with $ r\not=0 $. After this, the
compatibility is restored: indeed, if $ r_{0}\not=0 $, then $ \lim _{\left(s,r\right)\to\left(0,r_{0}\right)}s/r $
exists and is 0. \end{example}

\begin{remark} Compare the need to consider an open subset of values $ r $ with
the compatibility relation of Section~\ref{s0.77}. \end{remark}

This example shows a geometric reason for a diffeomorphism not
necessarily acting on a meromorphic function if $ \dim =\infty $: the meromorphic
function is defined via its ``restrictions'' to a special exhausting
collection of finite-dimensional families (such as $ R_{d} $) inside our
infinite-dimensional manifold; these restrictions are compatible.
However, one should not expect that a restriction to {\em any\/}
finite-dimensional family is compatible with this exhausting collection
of families; some families are ``too skew'' to be compatible.

Now a diffeomorphism can map any ``good'' family to a skew\footnote{W.r.t.~the given meromorphic function.} one; then
the action of this diffeomorphism on the meromorphic function is not
defined. E.g., $ h_{s} $ maps any ``good'' family to a family on which $ \operatorname{Eval}_{s} $ has a
pole. {\em This\/} is the reason why $ {\mathbit h}_{s}\left(\operatorname{Eval}_{t}\right) $ has a pole at $ s=t $.

\subsection{Basic and generalized functions } The principal reason why OPEs appear
in the formal theory ``in a formal way'' only is the fact (see Section~%
\ref{s0.44}) that the family $ {\mathbit E}_{t} $ of operators is defined only as a formal
series in $ t $. Comparing this with Theorem~\ref{th02.10}, one can see that the
reason for this deficiency is that $ \operatorname{Func}\left({\mathcal L}^{f}X\right) $ is ``too small'';
multiplication by $ \operatorname{Eval}_{t} $ ``moves'' a functional ``outside'' of the set
$ \operatorname{Func}\left({\mathcal L}^{f}X\right) $. Our approach increases the vector space of function on the
loop space so that the space remains invariant under multiplication by $ \operatorname{Eval}_{t} $.

To show how consideration of $ {\mathcal L}^{r}X $ is related to such an increase,
recall the usual duality between basic and generalized functions. While
for finite-dimensional vector spaces increasing $ V $ increases the dual
space $ V^{*} $ as well, for topological vector spaces increasing $ V $ (in such a
way that $ V $ is dense in the result) {\em decreases\/} $ V^{*} $. To make the space of
generalized functions as large as possible, one considers the space of
base functions which is as small as possible. The similar effect holds in
non-linear situation, when one considers $ \operatorname{Func}\left(X\right) $ for an
infinite-dimensional manifold $ X $ instead of the set of linear functionals
on an infinite-dimensional vector space $ V $. If $ X $ is dense in $ X'\supset X $, then
$ \operatorname{Func}\left(X'\right) $ is a subset of $ \operatorname{Func}\left(X\right) $.

Thus the reason why the operators $ {\mathbit E}_{t} $ act in $ \operatorname{Func}\left({\mathcal L}^{r}X\right) $ as honest
operators is that $ \operatorname{Func}\left({\mathcal L}^{r}X\right) $ is very large due to $ {\mathcal L}^{r}X $ being very small
(comparing to $ {\mathcal L}^{r}X $). In fact $ {\mathcal L}^{r}X $ is the orbit in $ {\mathcal L}^{f}X $ of the action of
diffeomorphisms $ G_{s}^{\left(\varepsilon\right)} $ and/or $ H_{s} $ (from Section~\ref{s0.05}). Thus one should
expect that the vector space $ \operatorname{Func}\left({\mathcal L}^{r}X\right) $ is the largest possible space of
functions in which the operators $ {\mathbit E}_{t} $ and $ {\mathbit h}_{s} $, $ {\mathbit d}_{s} $ act as honest (partial)
operators.

\subsection{The plan } In this paper we {\em start\/} at the situation
described in the previous sections. The principal objects we work with
are families of partial linear operators in the vector space
$ V=\operatorname{NaiveFunc}\left({\mathcal L}^{r}X\right) $. The families $ {\mathbit E}_{t} $, $ {\mathbit d}_{s} $ and $ {\mathbit h}_{s} $ become just particular
examples of such families.

For the coarse theory of this paper the target manifold $ X $ of
the loop does not matter much, we consider the case of the multiplicative
group $ X={\mathbb C}^{*} $ only. In this case the set $ {\mathcal L}_{Z}X $ of meromorphic loops becomes
the set $ {\mathcal M}\left(Z\right) $ of meromorphic functions on $ Z $. To shorten the notation,
denote $ \operatorname{NaiveFunc}\left({\mathcal L}_{Z}{\mathbb C}^{*}\right) $ on $ {\mathcal M}\left(Z\right) $ by $ {\mathfrak M}{\mathcal M}\left(Z\right) $.

In Sections~\ref{s11.10} and~\ref{s11.20} we equip the vector space $ V={\mathfrak M}{\mathcal M}\left(Z\right) $
with two structures: a lattice of vector subspaces of $ V $, and notion of {\em a
smooth family\/} of vectors in $ V $. To define the lattice, define the support
of $ \operatorname{Eval}_{t,K}\in V $ as $ \left\{t\right\}\subset Z $; support of a rational expression is the union of
the supports of the symbols $ \operatorname{Eval}_{t,K} $ which appear in this expression. Now
to each Zariski open subset $ U\subset Z $ one can associate the vector subspace
$ V\left(U\right) $ consisting of $ v\in V $ such that $ \operatorname{Supp}\left(v\right)\subset U $.

In Section~\ref{s1.65} we introduce the notion of a {\em codominant conformal
field\/}; this is a family of partial linear operators in $ V $ which behaves
nicely with respect to two structures (lattice and smoothness) in $ V $.
Although $ {\mathbit E}_{t} $, $ {\mathbit d}_{s} $ and $ {\mathbit h}_{s} $ (and other family coming from algebraic geometry
of $ {\mathcal L}^{r}X $) are examples of such families, in what follows we completely
ignore the algebraic geometry of $ {\mathcal L}^{r}X $: it turns out that the two specified
above structures on $ V $ is {\em all\/} that is needed to treat conformal fields. In
particular, {\em any\/} two codominant conformal fields satisfy the OPEs.

In the rest of Section~\ref{h11} we define the rest of the principal
objects of this paper: the field $ {\mathfrak M}{\mathcal M}\left(Z\right) $, conformal operators, conformal
families, residues, fusion and localization. In Section~\ref{h12} we study a
particular case $ Z={\mathbb P}^{1} $ and show that on localization at $ 0\in{\mathbb P}^{1} $ one obtains
the standard vertex operators (this enhances the discussion of Section~%
\ref{s0.05}). In Section~\ref{h13} we show that appropriate projective bundles over
Hilbert schemes of points in $ Z $ provide good finite-dimensional models of
the infinite-dimensional vector space $ {\mathfrak M}{\mathcal M}\left(Z\right) $, and deduce that after
applying a natural restriction, conformal fields automatically satisfy
conformal associativity relations and vertex commutation relations. In
Section~\ref{h14} we show how chiral algebras are related to our definition of
conformal fields, calculate vertex commutation relations in the case $ \dim 
Z=1 $, and show how one of examples of Section~\ref{h12} gives rise to algebras
$ {\mathfrak g}{\mathfrak l}_{\infty} $ and $ {\mathfrak g}\left(\infty\right) $ when localized to Laurent series. The results of this
section show that the standard ``physical degree of strictness'' arguments
deducing Lie-algebraic identities comparing OPEs and cohomologies of a
plane with 3 removed lines can be made honestly for the conformal fields
of the type considered in this paper.

In Section~\ref{h15} we provide an independent description of
boson-fermion correspondence (this correspondence differs a little from
the standard one), and show why it provides heuristics for considering
the space $ {\mathfrak M}{\mathcal M}\left({\mathbb P}^{1}\right) $ as a Fock space (in $ 1+1 $-dimensional theory).

I am deeply grateful to R.~Borcherds, I.~M.~Gelfand, A.~Givental,
V.~Goryunov, I.~Grojnowsky, R.~Hartshorn, Y.-Z.~Huang, M.~Kapranov,
M.~Kontsevich,
A.~Kresch, F.~Malikov, N.~Reshetikhin, M.~Spivakovsky, V.~Serganova and
A.~Vaintrob for fruitful discussions during the work over these concepts.
The revision of the introduction was made possible by hospitality of
IHES, MPIM-Bonn and MSRI. This paper was partially supported by
Alfred~P.~Sloan research fellowship.

{\bf Revision history. }The revisions of this paper are archived as
\url{arXiv:math.AG/9805096} on the math preprint server archive
\url{arXiv.org/abs/math}. In addition to cosmetic changes, Revision II of
this paper contains a complete rework of the introduction (with the
exception of Section~\ref{s0.100}). The numeration of statements outside the
introduction did not change.

\subsection{Alternative approaches } Two properties of the lattice $ V\left(U\right) $ in $ V={\mathfrak M}{\mathcal M}\left(Z\right) $
introduced above are most important for our approach: one is the {\em fatness\/}:
for any element $ v\in V $ and any $ U\subset Z $ there is a family $ v_{t} $, such that $ v_{0}=v $, and
$ v_{t}\in V\left(U\right) $ for $ t\not=0 $. Compare this with Remark~\ref{rem1.110}. Another is a related
property of {\em stability\/}: given the existence of moduli spaces of
subvarieties of $ Z $, say that $ U_{t} $ is a {\em smooth
one-parameter family of Zariski open subsets\/} of $ Z $ if the complement
family $ {\mathbit C}U_{t} $ is such. This leads to a
one-parameter family of vector subspaces $ V\left(U_{t}\right)\subset V $. Due to the description
of the tangent space to the Grassmannian, any smooth variation $ W_{t} $ of a
vector subspace $ W\subset V $ leads to a {\em variation operator\/} $ W \to V/W $. Now it is
easy to see that the variation operator for the family $ V\left(U_{t}\right) $ is actually
0. This is the stability property: $ V\left(U_{t}\right) $ is a non-trivial family with
vanishing infinitesimal variation. Conjecturally, the fatness and
stability properties are what substitutes the {\em factorization property\/} of
the Drinfeld--Beilinson approach to the loop spaces (\cite{BeiDriChiral} and
\cite{KapVas01Ver}).

Since in the case $ Z={\mathbb P}^{1} $ the families of operators $ {\mathbit E}_{t} $, $ {\mathbit d}_{s} $ and $ {\mathbit h}_{s} $
coincide with classical vertex operators, in general case they are
multidimensional analogues of vertex operators. Recall that the general
theory of \cite{Zakh98Mer} shows that at least in the case $ Z={\mathbb P}^{1} $ one can
describe $ {\mathfrak M}{\mathcal M}\left(Z\right) $ as the field of meromorphic functions on $ {\mathcal M}\left(Z\right) $ (see Section~%
\ref{s0.77}). Since for general $ Z $ this is not clear yet, it may be too early
to identify multi-dimensional theory of vertex operators with birational
geometry of sets $ {\mathcal M}\left(Z\right) $; however, this is a very natural conjecture.

A different approach to multidimensional vertex operators is
proposed in \cite{Bor97Ver}. Note that Hilbert schemes in relation to vertex
operators appear also in \cite{Nak96Lec} and \cite{Gro96Ins}, but the setting of
these papers do not look directly related to what we discuss here.

\subsection{Arithmetic analogue }\label{s0.100}\myLabel{s0.100}\relax  Consider here the result of formally
substituting $ Z=\operatorname{Spec}{\mathbb Z} $ into the geometric considerations of the previous
sections. One should consider $ {\mathbb Q} $ as the set $ {\mathcal M}\left(\operatorname{Spec}{\mathbb Z}\right) $ of meromorphic
functions on $ Z $. Results of Section~\ref{s13.3} show that it is natural to
suppose that meromorphic functions on $ {\mathcal M}\left(\operatorname{Spec}{\mathbb Z}\right) $ may be defined via
restrictions from $ {\mathcal M}\left(\operatorname{Spec}{\mathbb Z}\right) $ to $ \Gamma\left(S,{\mathcal O}\right) $ for finite subschemes $ S\subset\operatorname{Spec}{\mathbb Z} $. Thus
a meromorphic function on $ {\mathcal M}\left(\operatorname{Spec}{\mathbb Z}\right) $ may be described by an ideal $ N{\mathbb Z}\subset{\mathbb Z} $ of $ S $
and a function on $ {\mathbb Z}/N{\mathbb Z} $. In other words, the obvious substitution for $ {\mathfrak M}{\mathbb Q} $
is the space of locally constant functions on integer adeles $ {\mathbit A}_{{\mathbb Z}} $.

There is a partially defined action of $ {\mathbb Q}^{\times} $ on $ {\mathfrak M}{\mathbb Q} $ given by the
following rule: an element $ F $ of $ {\mathfrak M}{\mathbb Q} $ is a function on $ {\mathbb Z}/N{\mathbb Z} $ for some $ N>0 $.
Element $ p/q\not=0 $ acts on $ F $ iff $ N $ can be taken prime to $ q $, then it acts as
multiplication by the corresponding invertible residue modulo $ N $. This
defines a partial multiplication in $ {\mathcal C}={\mathbb Q}^{\times}\ltimes{\mathfrak M}{\mathbb Q} $, and a partial action
of $ {\mathcal C} $ on $ {\mathfrak M}{\mathbb Q} $. This partial multiplication is the analogue of conformal
field theory, and the partial action gives analogues of vertex operators.

The previous interpretation has an advantage of $ {\mathfrak M}{\mathbb Q} $ being an algebra.
Sacrificing this property, one may consider locally constant partially
defined functions on adeles $ {\mathbit A}_{{\mathbb Q}} $ with compact {\em domain of definition\/} as a
substitution for $ {\mathfrak M}{\mathbb Q} $, i.e., functions with values in the union of numbers
and a new symbol $ \upsilon $ (for ``undefined''). In such a case one has a
well-defined action of the group $ {\mathbb Q}^{\times} $ on $ {\mathfrak M}{\mathbb Q} $. However, one still has no {\em fair\/}
algebra structure on $ {\mathcal C}={\mathbb Q}^{\times}\ltimes{\mathfrak M}{\mathbb Q} $, only a partial structure.

\section{Field of meromorphic functions on $ {\protect \mathcal M}\left(Z\right) $ }\label{h11}\myLabel{h11}\relax 

In this section we define principal objects (meromorphic
functions on $ {\mathcal M}\left(Z\right) $, conformal operators, conformal fields, residues,
fusions and localizations) of this treatise.

\subsection{Meromorphic functions }\label{s11.10}\myLabel{s11.10}\relax  Consider an irreducible quasiprojective
variety $ Z $ over a field $ {\mathbb K} $ with a line bundle $ {\mathcal L} $. Any two global sections
$ n $, $ d $ of $ {\mathcal L} $ which are not simultaneously identically zero give a
meromorphic function $ n/d $ on $ Z $. This provides a family of meromorphic
functions on $ Z $, this family is parameterized by $ X_{Z,{\mathcal L}}\buildrel{\text{def}}\over{=}{\mathbb P}\left(\Gamma\left(Z,{\mathcal L}\right)\times\Gamma\left(Z,{\mathcal L}\right)\right) $.
(Since most of the time we consider a fixed quasiprojective variety $ Z $,
usually we reduce this notation to $ X_{{\mathcal L}} $.)

We will need to consider $ X_{{\mathcal L}} $ as an algebraic variety, thus $ \Gamma\left(Z,{\mathcal L}\right) $
should be a finite-dimensional vector space. This is automatically true
if $ Z $ is projective. To {\em force\/} this to be true for an arbitrary
quasiprojective $ Z $, {\em define\/} a line bundle on such a $ Z $ to consist of a
projective completion $ \bar{Z} $ of $ Z $, and a usual line bundle $ \bar{{\mathcal L}} $ on $ \bar{Z} $. {\em Define\/}
$ \Gamma\left(Z,{\mathcal L}\right) $ as $ \Gamma\left(\bar{Z},\bar{{\mathcal L}}\right) $.

Fix a point $ z\in Z $. Suppose that $ z $ is not a base point of $ {\mathcal L} $. Define
the {\em evaluation function\/} $ E_{z} $ on $ X_{{\mathcal L}} $ by
\begin{equation}
n/d \mapsto n\left(z\right)/d\left(z\right) \in {\mathbb P}^{1}={\mathbb K}\cup\left\{\infty\right\}.
\notag\end{equation}
It is a meromorphic function on $ X_{{\mathcal L}} $, it is regular on an open subset
$ U_{{\mathcal L},z}\subset X_{{\mathcal L}} $ given by $ n\left(z\right)\not=0 $ or $ d\left(z\right)\not=0 $.

Similarly, if $ z $ is smooth, and $ v_{1},\dots ,v_{n} $ are meromorphic vector
fields on $ Z $ which have no pole at $ z $, then one can also define a {\em partial
derivative\/} $ E_{z}^{\left(v_{1}\dots v_{n}\right)} $ of the evaluation function, given by $ \xi \mapsto
E_{z}\left(v_{1}\dots v_{n}\xi\right) $. When we need to denote an arbitrary derivative, we will use
notation $ E_{z}^{\bullet} $. For a fixed $ z $ meromorphic functions $ E_{z}^{\bullet} $ span a vector space
with the same size as the symmetric power
of the
tangent space $ {\mathcal T}_{z}Z $. We will assume that $ \bullet $ is an index from a basis of this
symmetric power.

One can see that the vector space spanned by $ E_{z}^{\bullet} $ for a fixed smooth
$ z $ is the topological dual space $ \widehat{{\mathcal O}}_{z}^{\vee} $ to the space of formal functions $ \widehat{{\mathcal O}}_{z} $
near $ z $. This description works for a non-smooth $ z $ as well, define indices
$ \bullet $ for $ E_{z}^{\bullet} $ appropriately in the case of non-smooth $ z $ as well.

On the other hand, for a fixed multiindex $ \left(v_{1},\dots ,v_{n}\right) $ one can
consider $ E_{z}^{\left(v_{1}\dots v_{n}\right)} $ as a meromorphic function on $ X_{{\mathcal L}} $ which depends on $ z $.
It is well defined on the complement to the set of base points of $ {\mathcal L} $.
Taking different values of $ z $ and/or $ \left(v_{1},\dots ,v_{n}\right) $, one obtains many
meromorphic functions on $ X_{{\mathcal L}} $. Consider the field generated by these
meromorphic functions.\footnote{Note that this field is generated by $ \delta $-functions on $ Z $ and their
derivatives. Most approaches to the theory of generalized functions
involve consideration of $ {\mathcal D} $-modules on $ Z $. This may be a reason why
$ {\mathcal D} $-modules appear in definitions of chiral algebras of \cite{BeiDriChiral} and
\cite{HuaLep96Mod}.}

Consider how these functions depend on the choice of the line
bundle $ {\mathcal L} $. For any global section $ \beta\not\equiv 0 $ of another line bundle $ {\mathcal M} $ on $ Z $ one
can construct a mapping $ X_{{\mathcal L}} \xrightarrow[]{\iota_{\beta}} X_{{\mathcal L}\otimes{\mathcal M}} $ given by $ \left(n:d\right) \mapsto \left(n\beta:d\beta\right) $. Given a
meromorphic function $ F $ on $ X_{{\mathcal L}\otimes{\mathcal M}} $ such that the set of the base points of
this meromorphic functions does not contain $ \operatorname{Im}\iota_{\beta} $, one can construct a
meromorphic function $ \iota_{\beta}^{*}F\buildrel{\text{def}}\over{=}F\circ\iota_{\beta} $ on $ X_{{\mathcal L}} $.

Consider an arbitrary algebraic expression $ \widetilde{F} $ constructed from
several symbols $ \widetilde{E}_{z}^{\bullet} $, $ z\in Z $, with application of operations of addition,
multiplication and division. Given $ {\mathcal L} $, substitute the symbols $ \widetilde{E}_{z}^{\bullet} $ by the
functions $ E_{z}^{\bullet} $ on $ X_{{\mathcal L}} $, separately in the numerator and denominator of $ \widetilde{F} $.
Unless both the numerator and denominators become 0, one obtains a
function $ F^{\left(X_{{\mathcal L}}\right)} $ on $ X_{{\mathcal L}} $. One can try to perform the same procedure for $ {\mathcal L}\otimes{\mathcal M} $
taken instead of $ {\mathcal L} $. If both $ F^{\left(X_{{\mathcal L}}\right)} $ and $ F^{\left(X_{{\mathcal L}\otimes{\mathcal M}}\right)} $ are defined, and $ F^{\left(X_{{\mathcal L}\otimes{\mathcal M}}\right)} $ is
compatible with $ \iota_{\beta} $ in the above sense, then $ \iota_{\beta}F^{\left(X_{{\mathcal L}\otimes{\mathcal M}}\right)}=F^{\left(X_{{\mathcal L}}\right)} $. Investigate
for which values of $ {\mathcal L} $ and $ {\mathcal M} $ these assumptions hold.

Consider what is involved in the calculation of $ F $ at a given element
$ \left(n:d\right)\in X_{{\mathcal L}} $. Suppose for simplicity that $ F $ includes no partial derivatives
of $ E_{z} $. After all the operations are performed, $ F\left(\left(n:d\right)\right) $ is given by an
expression of the form
\begin{equation}
\frac{G\left(n\left(z_{1}\right),n\left(z_{2}\right),\dots ,n\left(z_{k}\right),d\left(z_{1}\right),d\left(z_{2}\right),\dots ,d\left(z_{k}\right)\right)}{
G'\left(n\left(z_{1}\right),n\left(z_{2}\right),\dots ,n\left(z_{k}\right),d\left(z_{1}\right),d\left(z_{2}\right),\dots ,d\left(z_{k}\right)\right)},
\label{equ6.08}\end{equation}\myLabel{equ6.08,}\relax 
here $ z_{i} $ are different points of $ Z $, $ G $ and $ G' $ are polyhomogeneous
polynomials of homogeneity degrees $ \left(\delta_{1},\dots ,\delta_{2k}\right) $, $ \left(\delta'_{1},\dots ,\delta_{2k}'\right) $. (These
degrees satisfy $ \delta_{l}+\delta_{l+k}=\delta'_{l}+\delta'_{l+k} $, thus the value of the above ratio is
not ``an element of an abstract $ 1 $-dimensional vector space'', but of $ {\mathbb K} $.)

For this expression to be correctly defined, one needs
$ G\left(n\left(z_{1}\right),\dots ,d\left(z_{1}\right),\dots \right) $ and $ G'\left(n\left(z_{1}\right),\dots ,d\left(z_{1}\right),\dots \right) $ to be not
simultaneously identically zero as functions of $ \left(n:d\right)\in X_{{\mathcal L}} $. Taking $ {\mathcal L} $
positive enough, one can ensure that $ n\left(z_{1}\right),n\left(z_{2}\right),\dots ,n\left(z_{k}\right) $ are
algebraically independent when $ n $ runs over $ \Gamma\left(Z,{\mathcal L}\right) $, same for $ d $. This
ensures correctness of $ F $.

If $ F $ contains partial derivatives, then in the above expression one
also needs to consider partial derivatives of $ n\left(z_{i}\right) $ w.r.t.~some (possibly
meromorphic) connection on $ {\mathcal L} $. After this all the arguments still work.

Thus it is clear that for any such formal expression $ \widetilde{F} $ one can pick
up a corresponding line bundle $ {\mathcal L} $ on $ Z $ such that $ F^{\left(X_{{\mathcal L}}\right)} $ is correctly
defined, and in fact this remains true for any line bundle of the form
$ {\mathcal L}\otimes{\mathcal M} $ with an ample $ {\mathcal M} $. The above arguments show that the manifolds $ X_{{\mathcal L}} $ are
objects of an appropriate (ordered) category, and any formal expression $ \widetilde{F} $
defines compatible meromorphic functions on ``tails'' of this category,
i.e., on objects which are ``greater'' than some fixed object. (These
functions are not necessarily compatible w.r.t.~{\em all\/} possible morphisms
between two given objects, but only w.r.t.~an open dense subset of the
set of morphisms, one can look up all the nasty details in \cite{Zakh98Mer}.)

The motivations above lead to the following

\begin{definition} Let $ Z $ be a quasiprojective variety. Let $ {\mathcal M}\left(Z\right) $ be a set of
meromorphic functions on $ Z $. Consider a field $ {\mathfrak M}{\mathcal M}\left(Z\right) $ freely generated by
symbols $ E_{z}^{\bullet} $, $ z\in Z $. Call elements of this field\footnote{Note that in analytic theories, $ E_{z}^{\bullet} $ ``span'' the vector space of generalized
functions
on $ Z $. Thus $ {\mathfrak M}{\mathcal M}\left(Z\right) $ is the field spanned by generalized functions on $ Z $.} {\em meromorphic functions\/} on
$ {\mathcal M}\left(Z\right) $. \end{definition}

Given any fixed meromorphic function $ F $ on $ {\mathcal M}\left(Z\right) $, and any given
meromorphic function $ \xi $ on $ Z $ one cannot guarantee that $ F\left(\xi\right) $ is a
well-defined number. However, one has seen that if $ \xi $ changes in
sufficiently big families, then $ F\left(\xi\right) $ is a well-defined meromorphic
function on the base of such a family.

\begin{remark} Note that studying the interplay between different families
$ X_{Z,{\mathcal L}} $ of meromorphic functions, i.e., studying the algebraic geometry of
$ {\mathcal M}\left(Z\right) $, leads to a {\em very\/} complicated theory (see \cite{Zakh98Mer}). Moreover, the
final result of this more complicated theory coincides with the above
definition, at least in the simplest---but most important---case $ Z={\mathbb P}^{1} $.

Thus we accept the above $ \operatorname{ad} $ hoc point of view, with {\em postulating\/} what
is a meromorphic function on $ {\mathcal M}\left(Z\right) $, consider $ {\mathcal M}\left(Z\right) $ as a set only, and
postulate which families of elements of $ {\mathcal M}\left(Z\right) $ are ``smooth''. \end{remark}

The collection $ \left\{{\mathfrak M}{\mathcal M}\left(U\right)\right\}_{U\subset Z} $ for open $ U\subset Z $ behaves as sections of a
sheaf with compact supports.\footnote{Recall that this means that the canonical mappings go in {\em opposite\/}
direction comparing to usual sheaves. (Thus dual vector spaces $ {\mathfrak M}{\mathcal M}\left(U\right)^{\vee} $ form
a sheaf.)} Denote it by $ {\mathfrak M}{\mathcal M} $.

Collection $ \left\{{\mathfrak M}{\mathcal M}\left(U\right)\right\}_{U\subset Z} $ forms a lattice of subfields of the field
$ {\mathfrak M}{\mathcal M}\left(Z\right) $. This lattice of vector subspaces becomes the playing field for
conformal fields, they will act as operators defined on small enough
spaces of this lattice.

\subsection{Families of meromorphic functions on $ {\protect \mathcal M}\left(Z\right) $ }\label{s11.20}\myLabel{s11.20}\relax  We {\em defined\/}
meromorphic functions on $ {\mathcal M}\left(Z\right) $ to be rational expressions of symbols $ E_{z_{i}}^{\bullet} $.
Similarly, {\em define\/} families of such functions by allowing both the
coefficients of the rational expressions and the points $ z_{i} $ to change in
families. Before giving formal definitions consider some examples of what
we want to define.

For an arbitrary regular mapping $ T \xrightarrow[]{\iota} Z $ consider $ E_{\iota\left(t\right)} $, $ t\in T $.
Obviously, one wants this family to be allowed.

Second, in the case $ T=Z={\mathbb A}^{1} $ one can consider the family $ F_{t}=\left(E_{t}-E_{-t}\right)/t $
if $ t\not=0 $, obviously it ``should be'' extensible to $ t=0 $ with the value
$ F_{0}=2E_{0}^{\left(d/dz\right)} $. Since $ F_{t} $ is even in $ t $, $ G_{t}\buildrel{\text{def}}\over{=}F_{\sqrt{t}} $ should be also allowed.
The following definition is good enough if $ Z $ is complete, as Section~%
\ref{s13.3} will show:

\begin{definition} \label{def11.23}\myLabel{def11.23}\relax  Consider an algebraic variety $ T $. Say that a family $ T
\xrightarrow[]{F_{\bullet}} {\mathfrak M}{\mathcal M}\left(Z\right) $ is {\em admissible\/} if there are mappings $ z_{i}\colon T \to Z $, $ i=1,\dots ,I $,
meromorphic vector fields $ v_{ij} $, $ i=1,\dots ,I $, $ j=1,\dots ,J_{i} $, on $ Z $ and a
meromorphic function $ \widetilde{F}\left(t,f_{1},\dots ,f_{I}\right)\not\equiv 0 $ on $ T\times{\mathbb A}^{I} $ such that $ F_{t} $ can be written
as $ \widetilde{F}\left(t,E_{z_{1}\left(t\right)}^{{\mathcal V}_{1}},\dots ,E_{z_{I}\left(t\right)}^{{\mathcal V}_{I}}\right) $, here $ {\mathcal V}_{i}=v_{i1}\dots v_{iJ_{i}} $, and $ v_{ij} $ has no pole at
$ z_{i}\left(t\right) $ for any $ i,j,t $.

Say that a family $ T \xrightarrow[]{F_{\bullet}} {\mathfrak M}{\mathcal M}\left(Z\right) $ is {\em strictly meromorphic on\/} $ T $ if it is
admissible when lifted to a finite covering of an open dense subset of $ T $,
and $ F_{t}\left(\xi\right) $, $ \xi\in X_{{\mathcal L}} $, is a meromorphic function on $ T\times X_{{\mathcal L}} $ for any ample line
bundle $ {\mathcal L}\geq{\mathcal L}_{0} $ on $ Z $, $ {\mathcal L}_{0} $ being some fixed line bundle on $ Z $.

For a complete $ Z $ say that a family $ T \xrightarrow[]{F_{\bullet}} {\mathfrak M}{\mathcal M}\left(Z\right) $ is {\em strictly smooth\/}
if
\begin{enumerate}
\item
it is strictly meromorphic;
\item
intersections of the divisor of poles of the meromorphic function
$ F_{t}\left(\xi\right) $ on $ T\times X_{{\mathcal L}} $ with fibers $ \left\{t\right\}\times X_{{\mathcal L}} $ have codimension 1\footnote{In algebro-geometric parlance, this divisor is {\em flat\/} over $ T $.} for any $ {\mathcal L}\geq{\mathcal L}_{0} $.
\end{enumerate}

\end{definition}

To generalize this definition to the case of quasiprojective $ Z $, note
that for any meromorphic vector field $ v $ on $ T $ and a meromorphic family $ F_{\bullet} $
the derivative $ vF_{\bullet} $ is a well-defined meromorphic family. This family is
smooth near $ t_{0}\in T $ if $ F_{\bullet} $ is smooth and $ v $ is regular near $ t_{0} $.

\begin{definition} \label{def11.25}\myLabel{def11.25}\relax  Consider a projective variety $ \bar{Z} $ and an open subset
$ Z\subset\bar{Z} $. Say that a family $ T \xrightarrow[]{F_{\bullet}} {\mathfrak M}{\mathcal M}\left(Z\right) $ is {\em strictly smooth\/} if it strictly
smooth when considered as a family in $ {\mathfrak M}{\mathcal M}\left(\bar{Z}\right) $, and for any fixed $ t_{0}\in T $ and
any meromorphic vector fields $ v_{1},\dots ,v_{k} $ on $ T $ which are regular near $ t_{0} $
the smooth families $ v_{1}\dots v_{k}F_{\bullet} $ in $ {\mathfrak M}{\mathcal M}\left(\bar{Z}\right) $ are in fact families in $ {\mathfrak M}{\mathcal M}\left(Z\right) $. \end{definition}

\begin{remark} The additional constraints of the latter definition may be
easily seen in the case $ T=\bar{Z}={\mathbb P}^{1} $, $ Z={\mathbb A}^{1} $. In this case the family
$ F_{t}=\frac{E_{t}}{1+tE_{t}} $ sends any $ t\in{\mathbb P}^{1} $ to $ {\mathfrak M}{\mathcal M}\left({\mathbb A}^{1}\right) $. However, the derivative $ E_{\infty} $ at $ t=\infty $
is not in $ {\mathfrak M}{\mathcal M}\left({\mathbb A}^{1}\right) $. \end{remark}

As with any $ \operatorname{ad} $ hoc definition, a lot of problems disappear, but some
other problems jump to existence. Say, with the above definition it is
hard to show that a subfamily of a strictly smooth family is strictly
smooth. On the other hand, it is obvious that strictly meromorphic
families form a commutative algebra with division, strictly smooth
mappings form a commutative algebra.

To avoid the problem with subfamilies, use the following trick:

\begin{definition} Say that a family $ T \xrightarrow[]{F_{\bullet}} {\mathfrak M}{\mathcal M}\left(Z\right) $, $ F_{\bullet}\not\equiv\infty $, parameterized by a
quasiprojective variety $ T $ is {\em smooth\/} if there is a quasiprojective variety
$ T' $, a strictly smooth mapping $ T' \xrightarrow[]{G_{\bullet}} {\mathfrak M}{\mathcal M}\left(Z\right) $, and a regular mapping $ T \to
T' $ such that the composition is $ F_{\bullet} $.

Say that a family $ T \xrightarrow[]{F_{\bullet}} {\mathfrak M}{\mathcal M}\left(Z\right) $ is {\em meromorphic\/} if there is a blow-up
$ T' \xrightarrow[]{\pi} T $ such for any point $ t\in T' $ there is a neighborhood $ U\ni t $ and a
regular scalar function $ f\not\equiv 0 $ on $ U $ such that $ f\cdot\left(F\circ\pi\right) $ gives a smooth mapping
$ U \to {\mathfrak M}{\mathcal M}\left(Z\right) $. Say that $ F_{\bullet} $ {\em is defined\/} at $ t\in T $ if it is smooth in a
neighborhood of $ t $. \end{definition}

Now by definition the restriction of a smooth mapping to a
subvariety is smooth or $ \equiv\infty $, and meromorphic families form a commutative
algebra with division, strictly smooth mappings form a commutative
algebra. Of course, the above trick is just for simplicity, and

\begin{nwthrmi} Any meromorphic family is strictly meromorphic, any smooth
family is strictly smooth. \end{nwthrmi}

We will not use this statement in what follows. It follows from the
results of Section~\ref{s13.3}.

\begin{example} Suppose that $ Z $ is a point $ O $. In such a case $ {\mathcal M}\left(Z\right)={\mathbb P}^{1} $, and $ {\mathfrak M}{\mathcal M}\left(Z\right) $
has one generator $ E_{O} $, thus is $ {\mathbb K}\left(E_{O}\right) $. Let $ T={\mathbb A}^{2} $ with coordinates $ \left(t,t'\right) $.
Then any element of $ {\mathbb K}\left(E_{O},t,t'\right) $ defines a meromorphic family of
meromorphic functions on $ {\mathcal M}\left(Z\right) $. Write this element as $ P/Q $, $ P,Q\in{\mathbb K}\left[E_{O},t,t'\right] $.
Suppose that $ P $ and $ Q $ are mutually prime, write $ Q=\sum_{k=0}^{K}a_{k}\left(t,t'\right)E_{O}^{k} $. The
family $ P/Q $ is smooth iff polynomials $ a_{k}\left(t,t'\right) $, $ k=0,\dots ,K $, have no common
zeros.

In particular, both $ F_{t,t'}=t/t' $ and $ G_{t,t'}=1/\left(t+t'E_{0}\right) $ are meromorphic,
but not smooth. $ F_{t,t'} $ becomes smooth either after a multiplication by $ t' $,
or after a blow up of $ \left(0,0\right)\in{\mathbb A}^{2} $. $ G_{t,t'} $ becomes smooth after multiplication
by any scalar polynomial $ P\left(t,t'\right) $ such that $ P\left(0,0\right)=0 $ {\em followed by\/} a blow-up
of $ \left(0,0\right)\in{\mathbb A}^{2} $. \end{example}

A fundamental property of smooth and meromorphic families is their
relation with the lattice $ \left\{{\mathfrak M}{\mathcal M}\left(U\right)\right\} $, $ U $ running through open subsets of $ Z $.
If $ F_{t} $, $ t\in T $, is a meromorphic family defined at $ t_{0} $, and $ F_{t_{0}}\in{\mathfrak M}{\mathcal M}\left(U\right) $, then
$ F_{t}\in{\mathfrak M}{\mathcal M}\left(U\right) $ for $ T $ in some opens neighborhood of $ t_{0} $. Moreover, say that $ U_{s} $,
$ s\in S $, is a {\em family of open subsets\/} of $ Z $ if there is an open subset $ {\mathbit U}\subset Z\times S $
such that $ {\mathbit U}\cap\left(Z\times\left\{s\right\}\right)=U_{s}\times\left\{s\right\} $. Then for any meromorphic family $ F_{t} $ and any
family of open subsets $ U_{s}\subset Z $ the subset
\begin{equation}
\left\{\left(t,s\right)\in T\times S \mid F_{t}\in{\mathfrak M}{\mathcal M}\left(U_{s}\right)\right\}
\notag\end{equation}
is open.

Call this property {\em a locality property}. Heuristically, it states
that a function $ F\in{\mathfrak M}{\mathcal M}\left(Z\right) $ has $ 0 $-dimensional support in $ Z $. Strictly
speaking, it would be better to define $ {\mathfrak M}{\mathcal M} $ as a sheaf on the Hilbert
scheme of subschemes of $ Z $ of dimension 0, but we want to avoid all the
algebro-geometric complications until Section~\ref{h13}.

\begin{remark} Here is another trivial but important observation: for any
quasiprojective variety $ S $ with an open dense subset $ U $ any smooth family
on $ U $ extends to a meromorphic family on $ S $. \end{remark}

\subsection{Action of automorphisms of $ {\protect \mathcal M}\left(Z\right) $ }\label{s11.31}\myLabel{s11.31}\relax  We are going to define some
partially defined mappings of vector spaces $ {\mathfrak M}{\mathcal M}\left(Z\right) $.

Consider a fixed meromorphic function $ \eta\not\equiv 0 $ on $ Z $. Multiplication by $ \eta $
gives an invertible mapping $ {\mathbit m}_{\eta}\colon {\mathcal M}\left(Z\right) \to {\mathcal M}\left(Z\right) $. Consider how this mapping
acts on meromorphic functions on $ {\mathcal M}\left(Z\right) $. Given such a meromorphic function
$ F\left(\xi\right)\in{\mathfrak M}{\mathcal M}\left(Z\right) $, try to find a meromorphic function $ {\mathbit M}_{\eta}F $ which satisfies
$ \left({\mathbit M}_{\eta}F\right)\left(\xi\right) \buildrel{\text{def}}\over{=}F\left(\eta\xi\right) $. For $ {\mathbit M}_{\eta}F $ to exist, the mapping $ {\mathbit m}_{\eta} $ and families $ X_{{\mathcal L}} $
should satisfy some compatibility conditions.

First, one can assume that $ \eta=n_{0}/d_{0} $, $ n_{0} $, $ d_{0} $ being global sections of
a line bundle $ {\mathcal M} $ on $ Z $. Multiplication by $ \eta $ sends $ X_{{\mathcal L}} $ into $ X_{{\mathcal L}\otimes{\mathcal M}} $. If $ {\mathcal L} $ is
positive enough, then $ F $ is correctly defined on $ X_{{\mathcal L}\otimes{\mathcal M}} $. However, the image
$ {\mathbit m}_{\eta}X_{{\mathcal L}}\subset X_{{\mathcal L}\otimes{\mathcal M}} $ may lie inside the base set of the meromorphic function $ F|_{X_{{\mathcal L}\otimes{\mathcal M}}} $,
where it is not defined, thus it is possible that $ F\left(\eta\xi\right) $ is not defined
for $ \xi\in X_{{\mathcal L}} $ for any $ {\mathcal L} $.

Suppose that there is $ {\mathcal L}_{0} $ such that for $ {\mathcal L}\geq{\mathcal L}_{0} $ the above situation will
not happen, i.e., $ {\mathbit m}_{\eta}X_{{\mathcal L}}\subset X_{{\mathcal L}\otimes{\mathcal M}} $ intersects with the base set of $ F|_{X_{{\mathcal L}\otimes{\mathcal M}}} $ along
a subvariety of codimension two\footnote{The motivation of codimension 2 is functoriality w.r.t.~composition of
mappings $ {\mathbit m}_{\eta} $. The details are outside of the scope of this paper, see
\cite{Zakh98Mer} instead. For the problems which appear otherwise see Example~%
\ref{ex11.8}.} in $ {\mathbit m}_{\eta}X_{{\mathcal L}} $. (Recall that the base set is the
intersection of the divisor of zeros and the divisor of infinity, and it
has codimension two.) Then $ {\mathbit m}_{\eta}^{*}F $ is a well-defined meromorphic function on
$ X_{{\mathcal L}} $.

It is clear that obtained in such a way functions on $ X_{{\mathcal L}} $ (for all
$ {\mathcal L}\geq{\mathcal L}_{0} $) are compatible in the same sense as algebraic combinations of $ E_{z}^{\bullet} $
are. Thus $ {\mathbit m}_{\eta}^{*} $ gives a (partially defined) mapping from the set of
meromorphic functions on $ {\mathcal M}\left(Z\right) $ into itself. Denote this mapping by $ {\mathbit M}_{\eta} $.

Moreover, if $ {\mathbit M}_{\eta}F $ and $ {\mathbit M}_{\eta'}{\mathbit M}_{\eta}F $ are both well-defined, then so is $ {\mathbit M}_{\eta'\eta}F $,
and $ {\mathbit M}_{\eta'}{\mathbit M}_{\eta}F = {\mathbit M}_{\eta'\eta}F $. In particular, one can see that any mapping $ {\mathbit M}_{\eta} $ is
automatically invertible (on the domain of definition).

\begin{example} \label{ex11.8}\myLabel{ex11.8}\relax  Suppose that $ \eta $ has no pole or zero at $ z\in Z $. Then $ {\mathbit M}_{\eta}E_{z} $ is
well-defined, and equals $ C\cdot E_{z} $, $ C=\eta\left(z\right) $. It is clear that if $ \eta $ has zero at
$ z $, then one should require $ C=0 $, which contradicts invertibility of $ {\mathbit M}_{\eta} $.

Yet stronger breach of contract comes in the case when $ z $ is in the
base set of $ \eta $, when $ C $ is not defined at all. \end{example}

From the above description one can easily conclude that $ {\mathbit M}_{\eta} $ is well
defined on $ {\mathfrak M}{\mathcal M}\left(U\right) $ if $ \eta $ has no zero and no pole on $ U\subset Z $, and that it sends
smooth families in $ {\mathfrak M}{\mathcal M}\left(U\right) $ to smooth families. Note also that the
operators $ {\mathbit h}_{s} $ and $ {\mathbit H}_{s} $ of Section~\ref{s0.05} are of this form.

\subsection{Multiplication in the field, smooth families, and conformal operators }
Consider a
different construction of mappings which send meromorphic functions on
$ {\mathcal M}\left(Z\right) $ into themselves. The study of these mappings will justify the
definition of general conformal fields.

The generators $ E_{z} $ of $ {\mathfrak M}{\mathcal M}\left(Z\right) $ form a smooth family parameterized by
$ z\in Z $. For any element $ F\in {\mathfrak M}{\mathcal M}\left(Z\right) $ one can consider a smooth family $ {\mathbit E}_{z}F
\buildrel{\text{def}}\over{=}E_{z}\cdot F $, here $ \cdot $ is the multiplication in $ {\mathfrak M}{\mathcal M}\left(Z\right) $. This defines a family of
operators $ {\mathbit E}_{z} $. The difference with the family $ {\mathbit M}_{\xi} $ is the fact that $ {\mathbit E}_{z} $ is
defined on the whole $ {\mathfrak M}{\mathcal M}\left(Z\right) $, but does not send $ {\mathfrak M}{\mathcal M}\left(U\right) $ into $ {\mathfrak M}{\mathcal M}\left(U\right) $ if $ z\notin U $.

Note that to make $ {\mathbit M}_{\xi} $ defined on $ {\mathfrak M}{\mathcal M}\left(U\right) $ the open subset $ U $ should be
small enough, but to make $ {\mathbit E}_{z} $ act in $ {\mathfrak M}{\mathcal M}\left(U\right) $ the open subset $ U $ should be big
enough. To make the situation more uniform, introduce the notion of
$ V $-conformal operator.

\begin{definition} \label{def11.37}\myLabel{def11.37}\relax  Consider partially defined mappings $ F_{k}\colon {\mathfrak M}{\mathcal M}\left(U_{k}\right) \to
{\mathfrak M}{\mathcal M}\left(Z\right) $, $ k=1,2 $, $ U_{1,2}\subset Z $. Say that two such mappings are {\em equivalent\/} if they
become equal when restricted to some $ {\mathfrak M}{\mathcal M}\left(U\right) $, $ U\subset Z $.

Say that such a mapping is {\em smooth\/} if it sends smooth families in its
domain of definition to smooth families. A {\em conformal operator\/} on $ {\mathfrak M}{\mathcal M}\left(Z\right) $ is
a class of equivalence of smooth mappings.

Let $ V $ be an open subset of $ Z $. Say that a conformal operator is
$ V $-{\em conformal\/} if for any $ U\subset U_{0} $ the corresponding mapping sends $ {\mathfrak M}{\mathcal M}\left(U\right) $ into
$ {\mathfrak M}{\mathcal M}\left(U\cup V\right) $, and smooth families in $ {\mathfrak M}{\mathcal M}\left(U\right) $ to smooth families in $ {\mathfrak M}{\mathcal M}\left(U\cup V\right) $. \end{definition}

Obviously, operators $ {\mathbit M}_{\eta} $ are $ \varnothing $-conformal, and $ \varnothing $-conformal operators
form a monoid with unity. The operators $ {\mathbit M}_{\eta} $ form a commutative subgroup of
this monoid.

\begin{remark} \label{rem1.110}\myLabel{rem1.110}\relax  Note that two equivalent smooth mappings $ {\mathbit c} $, $ {\mathbit c}' $ coincide
on the intersection of their domains. Indeed, they are both defined on
some $ {\mathfrak M}{\mathcal M}\left(U\right) $, and coincide on this subset. However, any element $ F $ of $ {\mathfrak M}{\mathcal M}\left(Z\right) $
can be included in a smooth family $ F_{s} $, $ s\in S $, which intersects $ {\mathfrak M}{\mathcal M}\left(U\right) $. Now
families $ {\mathbit c}F_{s} $ and $ {\mathbit c}'F_{s} $ are defined on a common open subset of $ S $, and
coincide there. Since smooth, they coincide where both are defined. Thus
two images of $ F $ coincide if both mappings are defined on $ F $.

In particular, any conformal operator can be maximally extended: it
is represented by a smooth mapping with maximal possible domain of
definition. In what follows we always consider such families only. \end{remark}

The particular definition of smoothness we used here is tuned to
make the proofs of principal theorems (associativity relation and
commutation relations) easier. One should expect that a definition
written completely in terms of the geometry of the vector space $ {\mathfrak M}{\mathcal M}\left(Y\right) $ and
the lattice of vector subspaces may work as well.

\subsection{Families of operators $ {\mathbit M}_{\eta} $ } In this section we provide an abstract
linear algebra framework formalizing properties of operators $ {\mathbit h}_{s} $ from
Section~\ref{s0.05}.

Instead of one meromorphic function $ \eta $ on $ Z $ consider a family $ {\text H}_{\bullet} $ of
meromorphic functions $ {\text H}_{s} $ parameterized by $ s\in S $ (in other words, there is a
mapping $ S \xrightarrow[]{\iota} X_{{\mathcal M}} $, and $ {\text H}_{s} $ is $ n/d $ if $ i\left(s\right)=\left(n:d\right)\in X_{{\mathcal M}} $, here $ {\mathcal M} $ is a line bundle
on $ Z $). For every fixed value of parameter $ s\in S $ one obtains a meromorphic
function $ {\text H}_{s} $ on $ Z $, thus a linear operator $ {\mathbit M}_{{\text H}_{s}} $ defined on $ {\mathfrak M}{\mathcal M}\left(U_{s}\right) $ for a
corresponding subset $ U_{s} $.

\begin{definition} For $ z\in Z $ let
\begin{equation}
S^{\left(z\right)} = \left\{s\in S \mid {\text H}_{s}\text{ has a pole or a zero at }z\right\}.
\notag\end{equation}
Say that $ {\text H}_{s} $, $ s\in S $, {\em dominates\/} $ U\subset Z $ if $ S $ cannot be covered by a finite number of
sets $ S^{\left(z\right)} $, $ z\in U $. \end{definition}

\begin{example} An arbitrary meromorphic function $ \eta $ on $ Z $ can be considered as
a family parameterized by a $ 1 $-point set. It dominates the complement to
its divisor of zeros and poles. \end{example}

Consider a family $ {\text H}_{s} $, $ s\in S $, which dominates $ U $. Then for every
meromorphic function $ F\in{\mathfrak M}{\mathcal M}\left(U\right) $ there is a non-empty open subset $ S^{\left(F\right)} $ of $ S $
such that $ {\mathbit M}_{{\text H}_{s}}F $ is well-defined for $ s\in S^{\left(F\right)} $. One can see an immediate
advantage of families over individual functions: families dominate bigger
open subsets than individual functions, and can even dominate a whole
projective variety:

\begin{example} Consider a projective variety $ Z $, and a very ample line bundle
$ {\mathcal M} $ on it. There is a tautological meromorphic function $ {\text H}_{s} $ parameterized by
$ s\in S=X_{{\mathcal M}} $. Since $ {\mathcal M} $ has no base points, subsets $ S^{\left(z\right)} $ have smaller dimension
than $ S $, thus $ {\text H}_{s} $ dominates $ Z $. \end{example}

\subsection{Conformal fields }\label{s1.65}\myLabel{s1.65}\relax  The following property singles out family of
conformal operators which behave similar to $ \varnothing $-conformal operators when
considered together.

\begin{definition} Consider a family $ {\mathbit C}_{s} $, $ s\in S $, of partially defined mappings
each acting from a subset of $ {\mathfrak M}{\mathcal M}\left(Z\right) $ into $ {\mathfrak M}{\mathcal M}\left(Z\right) $. Suppose that for any $ s\in S $
the mappings $ {\mathbit C}_{s} $ is a conformal operator. Say that this family is
{\em codominant\/} if
\begin{enumerate}
\item
for any non-empty open subset $ V $ of $ Z $ there is a non-empty open dense
subset $ V^{\left(S\right)} $ of $ S $ such that the mapping $ {\mathbit C}_{s} $ is a $ V $-conformal operator in
$ {\mathfrak M}{\mathcal M}\left(Z\right) $ if $ s\in V^{\left(S\right)} $.
\item
If $ V $ changes in a family, then $ V^{\left(S\right)} $ can be chosen to change in a
family.
\end{enumerate}
\end{definition}

In other words, a codominant family may increase support, but it
will increase it by an arbitrary small amount, thus codominant families
are closest analogues of $ \varnothing $-conformal operators. Obviously, the family $ {\mathbit E}_{z} $,
$ z\in S\subset Z $ is codominant iff $ S $ is dense.

Consider a family $ {\text H}_{s} $ which dominates $ U $. If $ F\in{\mathfrak M}{\mathcal M}\left(U\right) $, then $ {\mathbit M}_{{\text H}_{s}}F $ is
well-defined for $ s $ in a non-empty open subset $ S^{\left(F\right)} $ of $ S $, moreover, it can
be written as a rational expression of $ E_{z}^{\bullet} $ {\em with regular functions on\/} $ S^{\left(F\right)} $
{\em as coefficients}. Consider this as a hint for a definition of a smooth
family of conformal operators.

Another hint to understand what one wants from a family of
$ V $-conformal operators to be smooth is given by the locality property of a
smooth family $ E_{z} $ of meromorphic functions on $ {\mathcal M}\left(Z\right) $. These result in

\begin{definition} \label{def6.64}\myLabel{def6.64}\relax  Consider a family $ {\mathbit C}\left(s\right) $, $ s\in S $, each member of which is a
$ V_{s} $-conformal operator for a correspsonding subset $ V_{s}\subset Z $. Call this family
a {\em conformal family\/} or a {\em conformal field on\/} $ S $ {\em acting on\/} $ {\mathfrak M}{\mathcal M}\left(Z\right) $ if
\begin{enumerate}
\item
There is a non-empty open dense subset $ \widetilde{S}\subset S\times Z $ such that $ {\mathbit C}\left(s\right) $ is
defined on $ {\mathfrak M}{\mathcal M}\left(U\right) $ if $ \left\{s\right\}\times U\subset\widetilde{S} $;
\item
For any smooth family $ F_{t} $ in $ {\mathfrak M}{\mathcal M}\left(Z\right) $, $ t\in T $, let $ {\mathbit U} $ the subset of $ S\times T $ where
$ {\mathbit C}\left(s\right)F_{t} $
is well defined. Then $ {\mathbit U} $ is open, and the family $ {\mathbit C}\left(s\right)F_{t} $ is smooth on $ {\mathbit U} $.
\end{enumerate}

Say that $ {\mathbit C}\left(\bullet\right) $ {\em dominates\/} an open subset $ U\subset Z $ if $ \widetilde{S} $ can be chosen to
project {\em onto\/} $ U $, $ {\mathbit C}\left(\bullet\right) $ is {\em dominant\/} if it dominates $ Z $. (With appropriate
modification in the case of reducible $ S $.)

Say that $ {\mathbit C}\left(\bullet\right) $ is {\em defined\/} on $ \widetilde{S} $. Say that $ {\mathbit C}\left(\bullet\right) $ is {\em smooth\/} if $ \widetilde{S} $ can be
chosen to project onto $ S $.

\end{definition}

\begin{example} If a family $ {\text H}_{\bullet} $ of meromorphic functions dominates $ U $, then $ {\mathbit M}_{{\text H}_{\bullet}} $
is a $ \varnothing $-conformal field which dominates $ U $. \end{example}

Obviously, there is a maximal open $ U\subset Z $ which is dominated by a
conformal field $ {\mathbit C} $, and it is dense. For any $ F\in{\mathfrak M}{\mathcal M}\left(U\right) $ there is a non-empty
open subset $ S^{\left(F\right)} $ of $ S $ where $ {\mathbit C}\left(s\right)F $ is defined. Suppose that the family $ F_{t} $
from the last condition of the definition is in $ {\mathfrak M}{\mathcal M}\left(U\right) $. Then the subset $ {\mathbit U} $
is not empty, thus the family $ {\mathbit C}\left(s\right)F_{t} $ is meromorphic on $ S\times T $.

Thus dominant conformal fields send smooth families to meromorphic
ones! This is a first harbinger of the relation with the OPE expansions
of the field theory.

The heuristic for the definition of smoothness is the following one:
example of $ {\mathbit M}_{\eta} $ shows that one cannot expect all the interesting operators
to be defined on the whole vector space $ {\mathfrak M}{\mathcal M}\left(Z\right) $. The condition of
smoothness ensures that the domain of definition does not jump down to
$ \left\{0\right\} $ when the parameter of the family changes. (It is not clear whether
for such families the domain of definition {\em actually\/} always changes in a
smooth family.)

\subsection{Conformal field theory }\label{s11.71}\myLabel{s11.71}\relax  The principal property of conformal
fields is that they may be composed, but the composition of smooth fields
is not necessarily smooth:

\begin{theorem} Suppose that $ {\mathbit C}\left(s\right) $, $ s\in S $, $ {\mathbit C}'\left(s'\right) $, $ s'\in S' $, are conformal fields
acting on $ {\mathfrak M}{\mathcal M}\left(Z\right) $, and $ {\mathbit C}' $ is codominant. Then $ {\mathbit C}\left(s\right)\circ{\mathbit C}'\left(s'\right) $ is a conformal
field on an open dense subset of $ S\times S' $, and it is dominant if $ {\mathbit C} $ and $ {\mathbit C}' $
are dominant, and codominant if $ {\mathbit C} $ is codominant. \end{theorem}

\begin{proof} Suppose that $ {\mathbit C} $ dominates $ U $, $ {\mathbit C}' $ dominates $ U' $. Consider an
arbitrary $ F\in{\mathfrak M}{\mathcal M}\left(U'\right) $. Then $ {\mathbit C}'\left(s'\right)F $ is defined for $ s' $ in an open dense
subset $ T' $ of $ S' $, and depends smoothly on $ s' $.

Since $ {\mathbit C}' $ is codominant, one can decrease $ T' $ and $ U' $ in such a way
that $ {\mathbit C}'\left(s'\right)F $ is in $ {\mathfrak M}{\mathcal M}\left(U\right) $ if $ s'\in T' $ and $ F\in{\mathfrak M}{\mathcal M}\left(U'\right) $. Hence $ {\mathbit C}\left(s\right){\mathbit C}'\left(s'\right)F $ is
defined for $ s $ in an open dense subset of $ S $ for any given $ s'\in T' $.

Moreover, since the family $ {\mathbit C}'\left(s'\right)F $ is smooth, the family $ {\mathbit C}\left(s\right){\mathbit C}'\left(s'\right)F $
is smooth on an open dense subset of $ S\times S' $. If one lets $ F $ vary in some
family $ F_{t} $, $ t\in T $, inside $ {\mathfrak M}{\mathcal M}\left(U'\right) $, then the same is true for the family
$ {\mathbit C}\left(s\right){\mathbit C}'\left(s'\right)F_{t} $, which verifies the second condition of Definition~\ref{def6.64}.

Let $ {\mathbit C}\left(s\right) $ be defined on $ {\mathfrak M}{\mathcal M}\left(V_{s}\right) $, here $ V_{s} $ give a family of open subsets
of $ Z $. Since $ {\mathbit C}' $ is codominant, one can pick up $ V_{s}^{\left(S'\right)} $ to be the subset of
$ S' $ from the definition of codominantness, let $ {\mathbit V} $ be the corresponding
total subset of $ S\times S' $. One can conclude that $ {\mathbit C}\left(s\right){\mathbit C}'\left(s'\right)F $ is defined if
$ F\in{\mathfrak M}{\mathcal M}\left(V_{s}\cap U'\right) $ and $ \left(s,s'\right)\in{\mathbit V} $. This verifies the first condition from
Definition~\ref{def6.64}. If $ {\mathbit C} $ is dominant, $ \left\{V_{s}\right\} $ covers $ Z $, which proves that
$ {\mathbit C}\left(s\right){\mathbit C}'\left(s'\right) $ dominates $ U' $.

To show that the composition is codominant if $ {\mathbit C} $ is is
straightforward. \end{proof}

Note that one cannot expect the composition of smooth families to be
smooth:

\begin{example} \label{ex11.180}\myLabel{ex11.180}\relax  Consider the family $ {\mathbit M}_{\xi} $ parameterized by $ \xi\in X_{{\mathcal M}} $ and the
family $ {\mathbit E}_{z} $ parameterized by $ z\in Z $. If $ {\mathcal M} $ is ample, both families are smooth,
dominant and codominant. Moreover, there is a commutation relation
\begin{equation}
{\mathbit M}_{\xi}\circ{\mathbit E}_{z}=\left< \xi,z \right>{\mathbit E}_{z}\circ{\mathbit M}_{\xi},
\label{equ11.18}\end{equation}\myLabel{equ11.18,}\relax 
here $ \left< \xi,z \right> $ is a scalar meromorphic function on $ X_{{\mathcal M}}\times Z $ of taking the
value of $ \xi $ at $ z $. Since $ \left< \xi,z \right> $ has poles on $ S=X_{{\mathcal M}}\times Z $, the set $ \widetilde{S} $ for $ {\mathbit E}_{z}\circ{\mathbit M}_{\xi} $
cannot project on the points of $ S $ which are on the divisor of poles of $ \left<
\xi,z \right> $.

Since the calculation of the above relationship is confusion-prone,
let us outline it here. Given $ F\left(\xi_{1}\right) $, $ {\mathbit E}_{z}F\left(\xi_{1}\right)=\xi_{1}\left(z\right)F\left(\xi_{1}\right) $, thus
\begin{equation}
{\mathbit M}_{\xi}{\mathbit E}_{z}\left(F\left(\xi_{1}\right)\right)={\mathbit M}_{\xi}\left(\xi_{1}\left(z\right)F\left(\xi_{1}\right)\right)=\left(\xi\xi_{1}\right)\left(z\right)F\left(\xi\xi_{1}\right)=\xi\left(z\right)\xi_{1}\left(z\right)F\left(\xi\xi_{1}\right).
\notag\end{equation}
Similarly, $ {\mathbit E}_{z}{\mathbit M}_{\xi}\left(F\left(\xi_{1}\right)\right)=\xi_{1}\left(z\right)F\left(\xi\xi_{1}\right) $. \end{example}

\subsection{Homogeneous variant } The vector spaces in the lattice on which
conformal operators act are {\em huge}, and the theory can be simplified a lot
by reducing them. One of the variants is to consider only the semisimple
part with respect to the natural action of the multiplicative group
$ {\mathbb G}_{m}=\operatorname{GL}\left(1\right) $ on $ {\mathcal M}\left(Z\right) $.

Let $ \lambda\in{\mathbb G}_{m} $ acts on $ X_{{\mathcal M}} $ by $ \left(n:d\right) \to \left(\lambda n:d\right) $. Consider only those
meromorphic functions on $ X_{{\mathcal M}} $ which are homogeneous w.r.t.~this action. In
other words, for some $ \delta\in{\mathbb Z} $
\begin{equation}
F\left(\lambda n:d\right) = \lambda^{\delta}F\left(n:d\right)
\notag\end{equation}
(this gives some additional restriction on degrees of polyhomogeneity of
$ G $ and $ G' $ in~\eqref{equ6.08}). Similarly, consider only conformal fields of
degree $ \delta'\in{\mathbb Z} $, i.e., operators which send homogeneous functions of degree $ \delta $
into functions of degree $ \delta+\delta' $.

\begin{example} The degree of $ {\mathbit M}_{\bullet} $ is 0, the degree of $ {\mathbit E}_{\bullet} $ is 1. Composition adds
degrees. In particular, the family $ {\mathbit E}_{z}\circ{\mathbit E}_{z'}^{-1} $, $ \left(z,z'\right)\in Z\times Z $, is of degree 0. \end{example}

When one considers conformal families of homogeneity degree 0, the
representation of $ \operatorname{GL}\left(1\right) $ in the vector space (or lattice of vector spaces)
is split into a direct sum of homogeneous components. Thus it is
reasonable to restrict attention to one homogeneous component.

Note that if $ Z $ has a distinguished point $ z_{0} $, then multiplication by
$ E_{z_{0}} $ identifies all the different homogeneous components. Thus it is
enough to consider the action in the subspace of homogeneity degree 0.

\subsection{Residues }\label{s11.8}\myLabel{s11.8}\relax  Recall that given a meromorphic function $ G $ on an
algebraic variety $ S $, and an irreducible hypersurface $ D\subset S $, one may
consider {\em Leray residues\/} $ \operatorname{Res}_{D}^{\left(n\right)}G $ of $ G $ along $ D $. If $ S $ is a direct product
$ D_{0}\times C $ with $ C $ being a curve, and $ D=D_{0}\times\left\{c_{0}\right\} $, $ c_{0}\in C $, then $ \operatorname{Res}_{D}^{\left(n\right)}G|_{d\in D} $ is
$ -n $-th Laurent coefficients of $ G|_{\left\{d\right\}\times C} $ along $ C $ at $ c_{0} $. In general, some
additional data make Laurent residues into meromorphic functions on $ D $,
and one can easily describe how to transform these functions if these
additional data are changed.

Say, let $ \alpha $ be a meromorphic function on $ S $ which has a zero of the
first order at $ D $, $ \omega=d\alpha $. $ \omega $ is a meromorphic $ 1 $-form on $ S $ which has no zero
or pole at $ D $, and vanishes when pulled back to $ D $. Let $ v $ be a meromorphic
vector field on $ S $ such that $ \left< \omega,v \right>=1 $, and $ v $ has no pole or zero on $ D $.
Using these data, one can describe the Leray residue of $ G $ as the residue
along $ v $ of the differential form $ \omega\left(s\right)\alpha\left(s\right)^{n-1}G $ at generic points of the
divisor $ D $. Fix such data $ \left(\alpha,v\right) $, in the future consider residues w.r.t.~
these data only. Note that one can also use these data to define the
residue on $ D\times T $ of a function on $ S\times T $ for any irreducible variety $ T $.

Residues of positive order vanish if $ G $ is regular at generic points
of $ D $. In this case residues of negative order are Taylor coefficients of
$ G $ in the direction of $ v $. Thus residues with $ n\geq1 $ describe the singular
part of $ G $ at $ D $. Abusing notations, call residues with positive and
negative $ n $ taken together {\em Laurent coefficients\/} of $ G $ at $ D $ (numbered with
opposite numbers).

Note that $ \operatorname{Res}_{D}^{\left(0\right)}G $ coincides with restriction to $ D $ if $ G $ is smooth at
generic points of $ D $. However, $ \operatorname{Res}_{D}^{\left(0\right)}G $ is defined even if $ G $ has a pole at
$ D $. The biggest contrast with restriction is that if $ D' $ is a different
hypersurface, then for any irreducible component $ D'' $ of $ D\cap D' $ one can
consider $ \operatorname{Res}_{D''}^{\left(k\right)}\operatorname{Res}_{D}^{\left(l\right)}G $ and $ \operatorname{Res}_{D''}^{\left(l\right)}\operatorname{Res}_{D'}^{\left(k\right)}G $. These two functions
coincide\footnote{Assume that the vector fields $ v $, $ v' $ associated to $ D $, $ D' $ are tangent
to the other hypersurface of a pair $ \left(D,D'\right) $.} if no pole of $ G $ contains $ D'' $, however, this is true no more if $ G $
is not regular on $ D'' $. This observation is the core of {\em conformal
commutation relations}.

In particular, if $ T\subset S $ is not contained in $ D $, this does not imply
that
\begin{equation}
\left(\operatorname{Res}_{D}^{\left(n\right)}G\right)|_{T}=\operatorname{Res}_{D\cap T}^{\left(n\right)}\left(G|_{T}\right).
\notag\end{equation}
This identity holds only if $ D\cap T $ is not contained in the divisor of poles
of $ G $.

\begin{example} Consider a function $ G=\frac{1}{x-z} $ on the plane, let $ T=\left\{x=0\right\} $,
$ D=\left\{z=0\right\} $. Then $ \operatorname{Res}_{D}^{\left(1\right)}G=0 $, but $ G|_{T} $ has a pole at the origin, thus
$ \operatorname{Res}_{0}^{\left(1\right)}\left(G|_{T}\right)=-1 $. \end{example}

Multiplying $ G $ by a power of $ \alpha $, one can see that

\begin{amplification} \label{amp11.81}\myLabel{amp11.81}\relax 
\begin{equation}
\left(\operatorname{Res}_{D}^{\left(n\right)}G\right)|_{T}=\operatorname{Res}_{D\cap T}^{\left(n\right)}\left(G|_{T}\right)
\notag\end{equation}
if $ G $ is defined and finite on $ D\cap T $, or if the only component of the
divisor of poles of $ G $ which contain $ D\cap T $ is $ D $. \end{amplification}

Consider a meromorphic family $ F_{s} $ in $ {\mathfrak M}{\mathcal M}\left(Z\right) $ parameterized by $ S $. Let $ D $
be an irreducible hypersurface in $ S $. Then $ F_{s} $ induces a meromorphic
function on $ S\times X_{{\mathcal L}} $ for $ {\mathcal L}\geq{\mathcal L}_{0} $. Thus for any $ n\in{\mathbb Z} $ one can consider the Leray
residue $ \operatorname{Res}_{D\times X_{{\mathcal L}}}^{\left(n\right)}\left(F_{s}\right) $. This residue is a meromorphic function on $ D\times X_{{\mathcal L}} $.

Consider a generic inclusion $ X_{{\mathcal L}} \to X_{{\mathcal L}'} $ for $ {\mathcal L}<{\mathcal L}' $. This gives an
inclusion $ S\times X_{{\mathcal L}} \to S\times X_{{\mathcal L}'} $, which is compatible with families $ F_{s} $ on these
spaces. To show that the residues are compatible too, note that this is
true if $ F_{s} $ is defined on $ D\times X_{{\mathcal L}} $ for $ {\mathcal L}\geq{\mathcal L}_{0} $. Otherwise $ F_{s} $ has a pole of the
{\em same\/} order $ d $ on $ D\times X_{{\mathcal L}}\subset S\times X_{{\mathcal L}} $ for $ {\mathcal L}\geq{\mathcal L}_{0} $, and multiplication of $ F $ by $ \alpha^{d} $
finishes the argument.

Thus $ \operatorname{Res}_{D\times X_{{\mathcal L}}}^{\left(n\right)}\left(F_{s}\right) $ provides a meromorphic family $ \operatorname{Res}_{D}^{\left(n\right)}F_{s} $ on $ {\mathfrak M}{\mathcal M}\left(Z\right) $
parameterized by $ D $.

\begin{definition} Call this family $ \operatorname{Res}_{D}^{\left(n\right)}F_{s} $ the {\em residue\/} of $ F_{s} $ at $ D\subset S $. \end{definition}

Consider a conformal field $ {\mathbit C}\left(s\right) $ on $ S $ acting on $ {\mathfrak M}{\mathcal M}\left(Z\right) $. Let $ D $ be an
irreducible hypersurface in $ S $. Suppose the $ {\mathbit C} $ dominates $ U\subset Z $. For any
$ F\in{\mathfrak M}{\mathcal M}\left(U\right) $ one can consider $ \operatorname{Res}_{D}^{\left(n\right)}\left({\mathbit C}\left(s\right)F\right) $.

\begin{definition} \label{def11.251}\myLabel{def11.251}\relax  The $ n ${\em -th fine residue\/} $ \operatorname{RES}_{D}^{\left(n\right)}{\mathbit C} $ of $ {\mathbit C} $ on $ D $ is the
family of mappings $ {\mathfrak M}{\mathcal M}\left(U\right) \to {\mathfrak M}{\mathcal M}\left(Z\right) $ parameterized by $ D $, and defined by $ F \mapsto
\operatorname{Res}_{D}^{\left(n\right)}\left({\mathbit C}\left(s\right)F\right) $. \end{definition}

The big problem with this notion is that $ \operatorname{Res}_{D}^{\left(n\right)}\left({\mathbit C}\left(s\right)F\right) $ does not
depend smoothly on $ F $, compare Example~\ref{ex11.270} below. One needs a
coarser notion that that: below we will find a restriction of this
operator to a vector subspace which is smooth.

For any smooth family $ F_{t} $, $ t\in T $, one can consider the residue of
$ {\mathbit C}\left(s\right)F_{t} $ on $ D\times T $. Note that if $ T_{1}\subset T $, then the residue on $ D\times T $ restricted to
$ D\times T_{1} $ does not necessarily coincide with the residue of the restriction of
$ {\mathbit C}\left(s\right)F_{t} $ on $ S\times T_{1} $. However, they coincide if $ D\times T_{1} $ is not contained in the
divisor of poles of $ {\mathbit C}\left(s\right)F_{t} $.

Next step is to study for which $ T $ and $ T_{1} $ the operations of taking a
residue and of restriction commute. Suppose that $ {\mathbit C} $ is defined on $ \widetilde{S}\subset S\times Z $.
Let $ \widetilde{R}=\left(S\times Z\right)\smallsetminus\widetilde{S} $, let $ \widetilde{R}_{0} $ be the union of irreducible components of $ \widetilde{R} $ which
differ from $ D\times Z $. Consider an open subset $ \widetilde{D}_{0}=\left(D\times Z\right)\smallsetminus\widetilde{R}_{0} $ of $ D\times Z $. Let $ D_{0} $ be
the projection of $ \widetilde{D}_{0} $ to $ D $, $ Z_{0} $ be the projection of $ \widetilde{D}_{0} $ to $ Z $.

If $ F_{t_{0}}\in Z_{0} $ for $ t_{0}\in T_{1} $, then for a generic $ s_{0}\in S $ the family $ {\mathbit C}\left(s\right)F_{t} $, $ s\in S $,
$ t\in T_{1} $, is meromorphic, and the only poles of this family near $ \left(s_{0},t_{0}\right) $ are
poles at $ D\times T_{1} $ (if any). Since $ t_{0}\in T $ as well, same is true for the family
$ {\mathbit C}\left(s\right)F_{t} $ on $ S\times T $. Amplification~\ref{amp11.81} implies that in this case taking
the residue along $ D $ and restriction to $ T_{1} $ commute indeed.

\begin{definition} Define a parameterized by $ D $ family of partially defined
linear operator $ \operatorname{Res}_{D}^{\left(n\right)}{\mathbit C} $ in the lattice $ \left\{{\mathfrak M}{\mathcal M}\left(U\right)\right\}_{U\subset Z} $ by
\begin{equation}
\left(\operatorname{Res}_{D}^{\left(n\right)}{\mathbit C}\right)F=\operatorname{Res}_{D}^{\left(n\right)}\left({\mathbit C}F\right),\qquad F\in{\mathfrak M}{\mathcal M}\left(Z_{0}\right).
\notag\end{equation}
\end{definition}

The following statement is obvious:

\begin{proposition} Suppose that $ {\mathbit C} $ is defined on $ \widetilde{S}\subset S\times Z $. Let $ \widetilde{R}=\left(S\times Z\right)\smallsetminus\widetilde{S} $, let $ \widetilde{R}_{1} $
be the union of irreducible components of $ \widetilde{R} $ which differ from $ D\times Z $.
Consider an open subset $ \widetilde{D}_{0}=\left(D\times Z\right)\smallsetminus\widetilde{R}_{0} $ of $ D\times Z $. Let $ D_{0} $ be the projection of
$ \widetilde{D}_{0} $ to $ D $, $ Z_{0} $ be the projection of $ \widetilde{D}_{0} $ to $ Z $.

Then $ \operatorname{Res}_{D}^{\left(n\right)}{\mathbit C} $ is a conformal field on $ D $ acting in $ {\mathfrak M}{\mathcal M}\left(Z\right) $ which is
defined on $ \widetilde{D}_{0} $, dominates $ Z_{0} $, and is smooth on $ D_{0} $.
\end{proposition}

\begin{remark} Heuristically, the above description of the conformal field
$ \operatorname{Res}_{D}^{\left(n\right)}{\mathbit C}\left(s\right) $ can be boiled down to: define it by
$ \left(\operatorname{Res}_{D}^{\left(n\right)}{\mathbit C}\left(s\right)\right)F=\operatorname{Res}_{D}^{\left(n\right)}\left({\mathbit C}\left(s\right)F\right) $ {\em for generic\/} $ F $, and extend smoothly to
non-generic $ F $. The necessity of such double limit comes from the fact
that $ \operatorname{Res}_{D}^{\left(n\right)}\left({\mathbit C}\left(s\right)F\right) $ does not depend smoothly on $ F $. \end{remark}

\begin{example} \label{ex11.270}\myLabel{ex11.270}\relax  Consider $ {\mathbit C}\left(s\right)={\mathbit M}_{\frac{1}{z-s}} $, $ z,s\in{\mathbb A}^{1} $. Then
$ \operatorname{Res}_{s=0}^{\left(1\right)}\left({\mathbit M}_{\frac{1}{z-s}}E_{z_{0}}\right) $ is 0 for generic $ z_{0} $, $ -1 $ if $ z_{0}=0 $. Yet worse,
$ \operatorname{Res}_{s=0}^{\left(0\right)}\left({\mathbit M}_{\frac{z+s}{z-s}}E_{z_{0}}\right) $ is 1 for generic $ z_{0} $, $ -1 $ if $ z_{0}=0 $, and
$ \operatorname{Res}_{s=s_{0}}^{\left(n\right)}\left({\mathbit M}_{\frac{z+s}{z-s}}E_{z_{0}}\right)=0 $ for any $ s_{0} $ and $ n>0 $. \end{example}

\begin{remark} Let $ {\mathbit C} $ is smooth on $ S_{0}\subset S $. Then $ D_{0}\supset S_{0}\cap D $. Note that one cannot
deduce anything else about $ D_{0} $ and $ Z_{0} $ by looking at the set which $ {\mathbit C} $
dominates, and at $ S_{0} $. Indeed, if $ S $ is $ {\mathbb A}^{2} $ with coordinates $ s_{1} $, $ s_{2} $, $ Z $ is $ {\mathbb A}^{1} $
with coordinate $ z $, and $ \widetilde{S} $ is given by $ s_{1}\not=0 $ and $ s_{2}=zs_{1} $, then the field is
smooth at $ s_{1}\not=0 $, but the residue at $ s_{1}=0 $ is smooth on $ s_{2}\not=0 $.

Similarly, one cannot assume that the residue is codominant even if
the initial field is. \end{remark}

\begin{remark} Note that the system of axioms we use here ensures that for any
given element of $ {\mathfrak M}{\mathcal M}\left(Z\right) $ (or even any family of elements)
\begin{equation}
\left(\operatorname{Res}_{D}^{\left(n\right)}{\mathbit C}\right)F=0\text{ for }n\gg 0.
\notag\end{equation}
However, these axioms do not ensure the fact that the above equality is
valid {\em uniformly\/} in $ F $.

To improve the system of axioms to exclude the above pathologies one
may need to consider as support not open subsets $ U\subset Z $, but closed
subschemes of $ Z\smallsetminus U $, similarly for the definition set of a conformal field.
Yet better alternative is to consider the Hilbert scheme of $ 0 $-dimensional
subschemes of $ Z $ (compare with Section~\ref{h13}).

Since we do not want to delve into such details, note that for all
examples of conformal fields considered here, the above assumption is
true. \end{remark}

\begin{remark} \label{rem11.315}\myLabel{rem11.315}\relax  Note that the above equality for fine residue
\begin{equation}
\left(\operatorname{RES}_{D}^{\left(n\right)}{\mathbit C}\right)F=0\text{ for }n\gg 0
\notag\end{equation}
is manifestly non-uniform in $ F $. Indeed, consider again $ {\mathbit C}\left(s\right)={\mathbit M}_{\frac{1}{z-s}} $,
$ Z=S={\mathbb A}^{1} $, $ D=\left\{0\right\} $. Let $ F $ be the $ k $-th derivative of evaluation function
$ E_{0}^{\left(d/dz\right)^{k}} $. Then $ \left(\operatorname{RES}_{D}^{\left(n\right)}{\mathbit C}\right)F=0 $ iff $ n>k+1 $. \end{remark}

\begin{remark} Note that though conformal operators are linear, and residues
may be written as integrals in the case $ {\mathbb K}={\mathbb C} $, conformal fields not
necessarily commute with taking a residue. Indeed, fix $ \xi_{0}\in{\mathcal M}\left(Z\right) $, let $ \xi_{0} $
has a simple pole at hypersurface $ H\subset Z $. Restrict what we did in Example~%
\ref{ex11.180} to $ \left\{\xi_{0}\right\} $: let $ {\mathbit C}={\mathbit M}_{\xi_{0}} $, consider a family $ E_{z} $, $ z\in Z $, of elements of
$ {\mathfrak M}{\mathcal M}\left(Z\right) $. Then $ \operatorname{Res}_{H}^{\left(0\right)}E_{z} =E_{z}|_{z\in H} $, so $ {\mathbit C}\left(\operatorname{Res}_{H}^{\left(0\right)}E_{z}\right) $ is not defined. On the
other hand,
\begin{equation}
\operatorname{Res}_{H}^{\left(0\right)}\left({\mathbit C}E_{z}\right) = \operatorname{Res}_{H}^{\left(0\right)}\xi_{0}\left(z\right)E_{z} = \operatorname{Res}_{H}^{\left(-1\right)}E_{z}=E_{z}^{v}|_{z\in H},
\notag\end{equation}
here $ v $ is an appropriate (meromorphic) vector field on $ Z $. (Proposition~%
\ref{prop13.330} will provide a way to fix this non-commutativity.) \end{remark}

\subsection{Fusions }\label{s11.85}\myLabel{s11.85}\relax  Consider now two smooth conformal fields $ {\mathbit C}\left(s\right) $, $ s\in S $, and
$ {\mathbit C}'\left(s'\right) $, $ s'\in S' $, both acting on $ {\mathfrak M}{\mathcal M}\left(Z\right) $. The composition $ {\mathbit C}\left(s\right)\circ{\mathbit C}\left(s'\right) $ is smooth
on an open subset of $ S\times S' $, after a blow up of the complement to this
subset one may assume that this composition is smooth on a complement to
a collection of divisors $ D_{i} $, $ i\in I $. Fixing a local equation $ \alpha_{i} $ and a
transversal vector field $ v_{i} $ for each $ D_{i} $, one obtains families
$ \operatorname{Res}_{D_{i}}^{\left(n\right)}\left({\mathbit C}\left(s\right)\circ{\mathbit C}\left(s'\right)\right) $, $ n\geq1 $. If $ i $-th family is not smooth on $ D_{i} $, one can
consider further blow-ups and residues, and so on, until all the residues
of high iteration will be smooth.

\begin{definition} The {\em fusions\/} of conformal fields $ {\mathbit C} $ and $ {\mathbit C}' $ are the elements of
the defined above collection of conformal fields. \end{definition}

\begin{remark} Since composition of conformal fields is associative, one
should be able to describe some ``pseudo-associativity'' restrictions on
the fusions of fusions, which are what remains of associativity relation.
The actual description of these relations (i.e., an axiomatic conformal
field theory, compare \cite{Bor97Ver}) is outside of the scope of this paper.

However, note that Section~\ref{s13.50} shows that a wide class of
conformal fields has pairwise fusions, but no ``new'' triple-product
fusions, which is the base of conformal field theory. \end{remark}

Note also that we use the word {\em fusion\/} is a slightly more general way
than it is used in the conformal field theory: their the families $ {\mathbit C} $ and
$ {\mathbit C}' $ are defined on the same parameter space $ S=S' $; moreover, $ \dim  S=1 $, and
the divisor $ D $ coincides with the diagonal in $ S $. Under this restriction
our definition of fusions coincides with the classical one; moreover, the
mentioned above absense of ``new'' triple-product fusions means that the
singularities of a tripple product, such as $ {\mathbit C}\left(s\right){\mathbit C}'\left(s'\right){\mathbit C}''\left(s''\right) $, happen
only on divisors $ \left\{s=s'\right\} $, $ \left\{s'=s''\right\} $, and $ \left\{s=s''\right\} $. Moreover, one can easily
check that this condition is equivalent to the conformal associativity
relation.

\subsection{Localization }\label{s11.91}\myLabel{s11.91}\relax  The main difference between conformal fields and
families of operators in a vector space is that the conformal fields are
only {\em partially defined}, i.e., they act only on small subspaces of the
lattice $ {\mathfrak M}{\mathcal M}\left(U\right) $, $ U\subset Z $, of subspaces of a vector space $ {\mathfrak M}{\mathcal M}\left(Z\right) $. To bridge this
gap, one may consider a direct limit over the family $ U_{i} $ elements of which
become smaller and smaller. To do this, take a point $ z_{0}\in Z $, and consider a
limit w.r.t.~the lattice of open subsets of $ Z $ which contain $ z_{0} $. (To
simplify notations, in what follows $ z_{0} $ is supposed to be smooth.)

Note that this process is absolutely analogues to the definition of
localization of $ Z $ at $ z_{0} $, however, we will see that our axioms lead to
consideration of formal completion $ \widehat{Z} $ of $ Z $ at $ z_{0} $ instead. Indeed, one can
localize any conformal field from $ Z $ to any open subset $ U $ of $ Z $. However,
this does not change the theory a lot, since $ {\mathcal M}\left(U\right)={\mathcal M}\left(Z\right) $, thus the only
things which change are domains of definition of families of meromorphic
functions on $ {\mathcal M}\left(Z\right) $, and dominantness and codominantness.

On the other hand, $ \cap_{U\ni z_{0}}{\mathfrak M}{\mathcal M}\left(U\right) $ contains only generators $ E_{z_{0}}^{\bullet} $, and
coincides with the field generated by these elements. While one can
calculate the value of an element of $ {\mathfrak M}{\mathcal M}\left(Z\right) $ at ``almost any'' meromorphic
function on $ Z $, it also makes sense to calculate values of elements of
$ \cap_{U\ni z_{0}}{\mathfrak M}{\mathcal M}\left(U\right) $ on ``almost any'' formal power series at $ z_{0} $. One can see that,
first, we achieved the settings of the formal approach (see Section~%
\ref{s0.44}), thus the in the conformal field theory {\em localization\/} of $ Z $ at $ z_{0} $
leads to study of functions on {\em completion\/} of $ Z $ at $ z_{0} $.

Second, the generators $ E_{z_{0}}^{\bullet} $ of $ {\mathfrak M}{\mathcal M}\left(\widehat{Z}\right)\buildrel{\text{def}}\over{=}\cap_{U\ni z_{0}}{\mathfrak M}{\mathcal M}\left(U\right) $ are enumerated by
multiindices $ K=\left(k_{1},\dots ,k_{d}\right) $, $ d=\dim  Z $, $ k_{i}\geq0 $. Thus $ {\mathfrak M}{\mathcal M}\left(\widehat{Z}\right) $ consists of
rational functions of symbols $ E^{\left(K\right)} $, $ K\in{\mathbb Z}_{\geq0}^{d} $.

At last, the lattice of subsets of $ \widehat{Z} $ is trivial, so is what remains
of the lattice of subspaces in $ {\mathfrak M}{\mathcal M}\left(Z\right) $ when one restricts attention to
$ {\mathfrak M}{\mathcal M}\left(\widehat{Z}\right)\subset{\mathfrak M}{\mathcal M}\left(Z\right) $. In particular, those conformal fields which may be
restricted to $ {\mathfrak M}{\mathcal M}\left(\widehat{Z}\right) $ act there as usual families of linear operators. Note
however, that not all conformal operators can be restricted to $ {\mathfrak M}{\mathcal M}\left(\widehat{Z}\right) $, the
operator should be $ V $-conformal for any $ V\ni z_{0} $, and should dominate some
neighborhood of $ z_{0} $. In particular, any dominant $ \varnothing $-conformal field acts in
$ {\mathfrak M}{\mathcal M}\left(\widehat{Z}\right) $ by families of linear operators meromorphically depending on the
parameter. For example, $ {\mathbit M}_{\xi} $ acts in $ {\mathfrak M}{\mathcal M}\left(\widehat{Z}\right) $ if $ \xi $ has no zero or pole at $ z_{0} $.

Consider now an arbitrary codominant conformal field $ {\mathbit C}\left(s\right) $, $ s\in S $. For
any $ V\ni z_{0} $ there is a non-empty dense subset $ S^{\left(V\right)} $ such that the operators
$ {\mathbit C}\left(s\right) $ are $ V $-conformal operators if $ s\in S^{\left(V\right)} $. A priori the intersection
$ \cap_{V\ni z_{0}}S^{\left(V\right)} $ may be empty, but it is possible to show that it is a non-empty
closed subset $ S^{\left\{z_{0}\right\}} $ of $ S $ if $ S $ is projective. One may conclude that
operators $ {\mathbit C}\left(s\right) $ and all their derivatives in $ s $ act on $ {\mathfrak M}{\mathcal M}\left(\widehat{Z}\right) $ if $ s\in S^{\left\{z_{0}\right\}} $,
and $ {\mathbit C} $ is defined at $ \left(s,z_{0}\right) $. The last condition gives an open subset of $ S $,
so in general one cannot conclude that the set of such $ s $ is non-empty,
even if $ {\mathbit C} $ is dominant. However, for the family $ {\mathbit E}_{z} $ this holds.

In general, one can see that the restriction to $ {\mathfrak M}{\mathcal M}\left(\widehat{Z}\right) $ of a
codominant conformal field $ {\mathbit C} $ on a projective variety $ S $ is defined as on
an open subset $ \widehat{U} $ of a formal completion $ \widehat{S}^{\left\{z_{0}\right\}} $ in $ S $ of a closed subvariety
$ S^{\left\{z_{0}\right\}} $ of $ S $. If $ {\mathbit C} $ is dominant, the open subset is an intersection of
non-empty open subset $ U $ of $ S $ with $ \widehat{S}^{\left\{z_{0}\right\}} $. However, this intersection does
not necessarily contain closed points of $ Z $, but it is non-empty as a
scheme. Indeed, functions on a formal completion are formal power series,
when some closed subsets are removed, one should considers the subalgebra
of fractions of such formal power series with denominators having zero at
the thrown away divisors only. (Note that this complication does not
appear if one considers the field $ {\mathbit E}_{z} $.)

If $ S $ is not projective, consider any projective completion $ \bar{S} $ of $ S $.
The family $ {\mathbit C} $ can be extended to $ \bar{S} $ (possibly with a loss of smoothness).
However, smoothness was not used in the above arguments, thus, as above,
a dominant codominant family $ {\mathbit C} $ can be restricted to an appropriate
localization of formal neighborhood of a closed subvariety $ S^{\left\{z_{0}\right\}} $ of $ \bar{S} $. If
one needs to localize $ \widehat{S}^{\left\{z_{0}\right\}} $ outside of $ S^{\left\{z_{0}\right\}} $, one needs to consider
formal Laurent series instead of formal Taylor series.

Consider operators $ {\mathbit M}_{\xi} $ in more details. Assume that $ \xi $ has no zero or
pole at $ z_{0} $, so can be written as a formal power series near $ z_{0} $. Thus the
subgroup of $ \left\{{\mathbit M}_{\xi}\right\} $ which acts in $ {\mathfrak M}{\mathcal M}\left(\widehat{Z}\right) $ is generated by multiplications by
formal power series $ \xi=\sum_{K}\xi_{K}{\text H}^{K} $, here $ {\text H}=\left(\eta_{1},\dots ,\eta_{d}\right) $, $ \eta_{i} $ are local
coordinates near $ z_{0} $. (One should consider $ \xi $ which give meromorphic
functions on $ Z $ only.)

Restrict ourselves to the homogeneous theory. Suppose that $ \xi $ has no
zero at $ z_{0} $. Since $ {\mathbit M}_{\xi} $ is of homogeneity degree 0, $ \xi=\operatorname{const} $ gives a trivial
operator in grading 0, thus one may normalize $ \xi $ at any point, say, by
$ \xi\left(z_{0}\right)=1 $. This means that only the {\em divisor\/} of $ \xi $ is relevant, thus one
obtains an action of the group of principal divisors on $ Z $ on $ {\mathfrak M}{\mathcal M}\left(\widehat{Z}\right) $. This
action is given by the product formula, which expresses Tailor
coefficients of $ \xi\xi' $ in terms of Taylor coefficients of $ \xi $ and $ \xi' $. In this
action each $ E^{\left(K\right)} $ is mapped to a linear combination of $ E^{\left(L\right)} $, $ l_{i}\leq k_{i} $ for any
$ i $.

As far as the divisor does not pass through $ z_{0} $, this action
preserves polynomials in $ E^{\left(K\right)} $. However, the action of divisor which pass
through $ z_{0} $ is not defined.

Consider $ {\mathbit E}_{z} $ now. Obviously, one cannot multiply an element of $ {\mathfrak M}{\mathcal M}\left(\widehat{Z}\right) $
by $ E_{z} $ unless $ z=z_{0} $. However, the action of $ {\mathbit E}_{z} $ can be written in terms of
the (formal) Taylor series of $ {\mathbit E}_{z} $ {\em near\/} $ z=z_{0} $:
\begin{equation}
E_{z}=\sum_{K}{\text H}^{K}\left(z\right)E^{\left(K\right)},
\notag\end{equation}
here a monomial $ {\text H}^{K} $ in $ \left(\eta_{1},\dots ,\eta_{d}\right) $ is considered as a function of $ z\approx z_{0} $.
Multiplication by such an expression gives not a rational function of
$ E^{\left(K\right)} $, but a power series.

Particular Taylor coefficients of $ {\mathbit E}_{z} $ in $ z $ near $ z=z_{0} $ are just
operators of multiplication by $ E^{\left(K\right)} $, thus preserve polynomials in $ E^{\left(K\right)} $.

\begin{proposition} The operators $ {\mathbit E}_{z} $ considered as power series near $ z=z_{0} $ act
in the field of fractions of formal power series in letters $ E^{\left(K\right)} $. Same if
one considers homogeneous fractions of power series, with $ \deg  E^{\left(K\right)}=1 $. \end{proposition}

Similarly

\begin{proposition} The operators $ \xi_{D} $ are defined for any principal divisor of
$ Z $ which does not pass through $ z_{0} $. They act in any algebra: polynomials in
$ E^{\left(K\right)} $, formal power series in $ E^{\left(K\right)} $ and fractions thereof. \end{proposition}

\section{Examples of conformal families }\label{h12}\myLabel{h12}\relax 

Consider the simplest possible example, when $ Z $ is a point $ z_{0} $. In such a
case the only generator in $ {\mathfrak M}{\mathcal M}\left(Z\right) $ is $ E\buildrel{\text{def}}\over{=}E_{z_{0}} $. The vector subspace of
homogeneous functions of degree $ d $ is generated by $ E^{d} $. Thus a homogeneous
conformal operator $ C $ of degree $ \delta $ is described by sequences of numbers $ C_{d} $,
$ d\in{\mathbb Z} $, such that $ C\left(E^{d}\right)=C_{d}E^{d+\delta} $. Restricted to functions of homogeneity
degree 0, this becomes a number. Similarly, conformal fields are as dull as
families of numbers, i.e., scalar functions.

Thus any real theory starts in dimension 1, with the simplest case
being of genus $ g=0 $. The biggest difference of this case with the
higher-dimensional situation is the fact that the group of all divisors
is generated by the subset of effective divisors of degree 1. Moreover,
if genus $ g=0 $, then the same is true for principal divisors. This allows
one to consider $ {\mathbit M}_{\xi} $ for $ \xi\in X_{{\mathcal O}\left(1\right)} $ only, since the rest can be described by
compositions of such operators.

In this section we investigate conformal fields in the particular
case $ Z={\mathbb P}^{1} $.

\subsection{Fields based on $ {\mathbb P}^{1} $ } Consider the simplest possible case of the family
$ {\mathbit M}_{\xi} $, with $ \xi\in X_{Z,{\mathcal L}} $, and $ Z={\mathbb P}^{1} $, $ {\mathcal L}={\mathcal O}\left(1\right) $. In this case $ X_{Z,{\mathcal L}} $ is $ 3 $-dimensional,
and described by coordinates $ \left(a:b:c:d\right) $ giving a meromorphic function $ z \mapsto
\frac{az+b}{cz+d} $. If one considers the homogeneous case only, then
multiplication by a non-0 constant does not matter, so one can consider
$ \widehat{\xi}=\left(\xi_{n},\xi_{d}\right)=\left(\left(a:b\right),\left(c:d\right)\right)\in{\mathbb P}^{1}\times{\mathbb P}^{1} $ instead, here $ \xi_{n}\ominus\xi_{d} $ is the divisor of $ \xi $.

For a $ 0 $-cycle $ {\mathfrak y}=\sum k_{i}z_{i} $ on $ Z $ with coefficients in $ {\mathbb Z} $ one can define
$ {\mathbit E}_{{\mathfrak y}}=\prod{\mathbit E}_{z_{i}}^{k_{i}} $ (recall that $ {\mathbit E}_{z} $ is the operator of multiplication by $ E_{z} $, and
different $ {\mathbit E}_{z} $ commute). Consider the homogeneous case again, then one
needs to consider $ 0 $-cycles of degree 0, and the simplest such a cycle is
given by $ z_{+}\ominus z_{-} $, it is parameterized by $ Z\times Z $. In the case $ Z={\mathbb P}^{1} $ one can see
that the set of parameters $ {\mathbb P}^{1}\times{\mathbb P}^{1} $ is the same as in the case of the
simplest family $ {\mathbit M}_{\xi} $.

Obviously, $ {\mathbit M}_{\xi} $, $ \xi\in X_{{\mathbb P}^{1},{\mathcal O}\left(1\right)} $, form a dominant and codominant conformal
field, similarly for $ {\mathbit E}_{{\mathfrak y}} $ if $ {\mathfrak y} $ runs through $ Z\times Z $, as above. If one considers
compositions of $ {\mathbit M}_{\xi} $ with itself, then it is a smooth family, and
compositions commute, the same for $ {\mathbit E}_{{\mathfrak y}} $. Repeat again that though we
restrict ourselves to the smallest-dimension families $ X_{{\mathcal O}\left(1\right)} $ and $ Z\times Z $, {\em any\/}
operator $ M_{\xi} $ can be described as a composition of such simplest $ {\mathbit M}_{\xi} $,
similarly for any $ 0 $-cycle of degree 0.

Consider now compositions $ {\mathbit M}_{\xi}\circ{\mathbit E}_{{\mathfrak y}} $ and $ {\mathbit E}_{{\mathfrak y}}\circ{\mathbit M}_{\xi} $. As we have seen it
already,
\begin{equation}
{\mathbit M}_{\xi}\circ{\mathbit E}_{{\mathfrak y}}=\left< \xi,{\mathfrak y} \right>{\mathbit E}_{{\mathfrak y}}\circ{\mathbit M}_{\xi},\qquad \left< \xi,{\mathfrak y} \right>\buildrel{\text{def}}\over{=}\xi\left(z_{+}\right)/\xi\left(z_{-}\right)\text{ if }{\mathfrak y}=z_{+}\ominus z_{-}.
\notag\end{equation}
Indeed, $ {\mathbit M}_{\xi}E_{z}\buildrel{\text{def}}\over{=}\xi\left(z\right)E_{z} $, $ {\mathbit E}_{z_{+}\ominus z_{-}}E_{z}\buildrel{\text{def}}\over{=}E_{z_{+}}E_{z_{-}}^{-1}E_{z} $. Thus
\begin{align} {\mathbit M}_{\xi}\circ{\mathbit E}_{z_{+}\ominus z_{-}}\left(E_{z}\right) & ={\mathbit M}_{\xi}\left(E_{z_{+}}E_{z_{-}}^{-1}E_{z}\right)=\xi\left(z\right)\xi\left(z_{+}\right)\xi^{-1}\left(z_{-}\right)E_{z_{+}}E_{z_{-}}^{-1}E_{z},
\notag\\
{\mathbit E}_{z_{+}\ominus z_{-}}\circ{\mathbit M}_{\xi}\left(E_{z}\right) & = \xi\left(z\right){\mathbit E}_{z_{+}\ominus z_{-}}\left(E_{z}\right) =\xi\left(z\right)E_{z_{+}}E_{z_{-}}^{-1}E_{z}.
\notag\end{align}
Moreover, one can see that the family $ {\mathbit E}_{{\mathfrak y}}\circ{\mathbit M}_{\xi} $ is smooth. Since $ \left< \xi,{\mathfrak y} \right>=
\frac{z_{+}-\xi_{n}}{z_{+}-\xi_{d}} \cdot\frac{z_{-}-\xi_{d}}{z_{-}-\xi_{n}} $, $ {\mathbit M}_{\xi}\circ{\mathbit E}_{{\mathfrak y}} $ has fusions on the hypersurfaces
$ z_{-}=\xi_{n} $ and $ z_{+}=\xi_{d} $. Only the $ n=1 $ fusions are not zero. For example, consider
$ z_{-}=\xi_{n} $. On this hypersurface the fusion is
\begin{equation}
\left(z_{-}-\xi_{d}\right) \frac{z_{+}-z_{-}}{z_{+}-\xi_{d}}{\mathbit E}_{z_{+}\ominus z_{-}}\circ{\mathbit M}_{\xi}
\notag\end{equation}
However, this conformal field has a pole at the intersection $ z_{+}=\xi_{d} $
with the
other divisor $ z_{+}=\xi_{d} $ with the residue
\begin{equation}
\left(z_{-}-z_{+}\right) \left(z_{+}-z_{-}\right){\mathbit E}_{z_{+}\ominus z_{-}}\circ{\mathbit M}_{\xi},
\notag\end{equation}
which is smooth, here $ \xi\left(z\right)=\frac{z-z_{-}}{z-z_{+}} $. One gets the same result when one
considers the different order of taking residues along divisors. Note
that the factor $ \left(z_{-}-z_{+}\right) \left(z_{+}-z_{-}\right) $ is just the normalization by the
determinant of conormal bundle to the intersection of divisors.

Note that any $ {\mathfrak y}=z_{+}\ominus z_{-} $ corresponds to a meromorphic function $ \xi\left({\mathfrak y}\right) $ (up
to proportionality) with a pole at $ z_{+} $ and a zero at $ z_{-} $. Let $ {\mathbit M}_{{\mathfrak y}}\buildrel{\text{def}}\over{=}{\mathbit M}_{\xi\left({\mathfrak y}\right)} $,
$ {\mathbit F}_{{\mathfrak y}}={\mathbit E}_{{\mathfrak y}}\circ{\mathbit M}_{{\mathfrak y}} $. Then $ {\mathbit F}_{{\mathfrak y}} $ coincides with the above smooth double-fusion without
the normalizing coefficient:
\begin{equation}
{\mathbit F}_{\left(z_{+},z_{-}\right)}={\mathbit E}_{z_{+}\ominus z_{-}}\circ{\mathbit M}_{\frac{z-z_{-}}{z-z_{+}}}.
\notag\end{equation}

One can describe the action of $ {\mathbit F}_{\left(z_{+},z_{-}\right)} $ on a rational function
$ F\left(E_{z_{1}},\dots ,E_{z_{k}}\right) $ via
\begin{equation}
E_{z} \mapsto \frac{z-z_{-}}{z-z_{+}}E_{z_{+}}E_{z_{-}}^{-1}E_{z}.
\notag\end{equation}
Note that the action on $ E_{z}^{-1} $ is
\begin{equation}
E_{z}^{-1} \mapsto \frac{z-z_{+}}{z-z_{-}}E_{z_{+}}E_{z_{-}}^{-1}E_{z}^{-1}.
\notag\end{equation}
Thus this family is defined on $ z\not=z_{+} $, $ z\not=z_{-} $ (this is another way to see
that the family is smooth).

\subsection{Fusions of $ {\mathbit F}_{{\mathfrak y}} $ }\label{s12.20}\myLabel{s12.20}\relax  Since any $ {\mathbit M}_{\xi} $ may be represented as a composition
of $ {\mathbit M}_{\xi_{i}} $ with $ \xi_{i}\in X_{{\mathcal O}\left(1\right)} $, and any $ {\mathbit E}_{{\mathfrak y}} $ can be represented as composition of
$ {\mathbit E}_{{\mathfrak y}_{i}} $, $ {\mathfrak y}_{i}=z_{i+}\ominus z_{i-} $, similarly for families, the fusions of the fields which
are compositions of $ {\mathbit M}_{\xi} $ and $ {\mathbit E}_{{\mathfrak h}} $ may be expressed via fusions of elementary
fields. Thus the most interesting fusions are the fusions between $ {\mathbit F}_{{\mathfrak y}} $ and
$ {\mathbit F}_{{\mathfrak y}'} $.

One can write the action of $ {\mathbit F}_{{\mathfrak y}}{\mathbit F}_{{\mathfrak y}'} $ on $ E_{z} $ as
\begin{equation}
E_{z} \mapsto \frac{z'_{+}-z_{-}}{z'_{+}-z_{+}}\cdot \frac{z'_{-}-z_{+}}{z'_{-}-z_{-}} \cdot\frac{z-z_{-}}{z-z_{+}} \cdot\frac{z-z'_{-}}{z-z'_{+}}
E_{z_{+}}E_{z_{-}}^{-1}E_{z'_{+}}E_{z'_{-}}^{-1}E_{z},
\notag\end{equation}
here $ {\mathfrak y}=z_{+}\ominus z_{-} $, $ {\mathfrak y}'=z'_{+}\ominus z'_{-} $. Similarly, the action on $ E_{z}^{-1} $ is
\begin{equation}
E_{z}^{-1} \mapsto \frac{z'_{+}-z_{-}}{z'_{+}-z_{+}} \cdot\frac{z'_{-}-z_{+}}{z'_{-}-z_{-}} \cdot\frac{z-z_{+}}{z-z_{-}} \cdot\frac{z-z'_{+}}{z-z'_{-}}
E_{z_{+}}E_{z_{-}}^{-1}E_{z'_{+}}E_{z'_{-}}^{-1}E_{z}^{-1},
\notag\end{equation}
One can see that this conformal field is defined if
\begin{equation}
\frac{z'_{+}-z_{-}}{z'_{+}-z_{+}} \cdot\frac{z'_{-}-z_{+}}{z'_{-}-z_{-}},\quad \frac{z-z_{+}}{z-z_{-}}\cdot\frac{z-z'_{+}}{z-z'_{-}}\text{,
}\frac{z-z_{-}}{z-z_{+}}\cdot\frac{z-z'_{-}}{z-z'_{+}}
\notag\end{equation}
are not infinity, thus if
\begin{equation}
z'_{+}\not=z_{+}\text{, }z'_{-}\not=z_{-}\text{, }z\not=z_{-}\text{, }z\not=z_{+}\text{, }z\not=z'_{-}\text{, }z\not=z'_{+}.
\notag\end{equation}
Correspondingly, $ {\mathbit F}_{{\mathfrak y}}{\mathbit F}_{{\mathfrak y}'} $ is smooth on
\begin{equation}
z'_{+}\not=z_{+}\text{, }z'_{-}\not=z_{-}.
\notag\end{equation}
Thus $ {\mathbit F}_{{\mathfrak y}}{\mathbit F}_{{\mathfrak y}'} $ has fusions at two divisor, both fusions are non-zero for $ n=1 $
only. Again, these fusions have residues on the intersection with the
other divisor, and the fusion of the second order is defined on the
surface
\begin{equation}
z'_{+}=z_{+}\text{, }z'_{-}=z_{-}.
\notag\end{equation}
of codimension 2. Plugging into the formula for $ {\mathbit F}_{{\mathfrak y}}{\mathbit F}_{{\mathfrak y}'} $, one can describe
this fusion as
\begin{equation}
E_{z} \mapsto \left(\frac{z-z_{-}}{z-z_{+}}\right)^{2}E_{z_{+}}^{2}E_{z_{-}}^{-2}E_{z},
\notag\end{equation}
one can continue taking the fusions obtaining newer $ 2 $-parametric conformal
fields as fusions of the older ones.

\subsection{A different reduction }\label{s12.30}\myLabel{s12.30}\relax  One has seen that fusions of conformal
fields $ {\mathbit M}_{\xi} $ and $ {\mathbit E}_{{\mathfrak y}} $ can be described as a reduction of the smooth field
$ {\mathbit E}_{{\mathfrak y}}\circ{\mathbit M}_{\xi} $ on subvariety $ \xi_{n}=z_{-} $, $ \xi_{d}=z_{+} $. Consider now a different reduction
$ \xi_{n}=z_{+} $, $ \xi_{d}=z_{-} $:

\begin{definition} Let $ {\mathbit G}_{{\mathfrak y}}={\mathbit E}_{{\mathfrak y}}\circ{\mathbit M}_{\xi} $ with $ \xi\left(z\right)=\frac{z-z_{+}}{z-z_{-}} $, $ {\mathfrak y}=z_{+}\ominus z_{-} $. \end{definition}

Obviously, $ {\mathbit G}_{{\mathfrak y}} $ is a smooth dominant and codominant family sending
\begin{equation}
E_{z} \mapsto \frac{z-z_{+}}{z-z_{-}}E_{z_{+}}E_{z_{-}}^{-1}E_{z}.
\notag\end{equation}
The restriction of $ {\mathbit G}_{{\mathfrak y}} $ to the diagonal $ z_{+}=z_{-} $ is the identity operator
$ \operatorname{id} $. One can write the action of $ {\mathbit G}_{{\mathfrak y}}{\mathbit G}_{{\mathfrak y}'} $ on $ E_{z} $ as
\begin{equation}
E_{z} \mapsto \frac{z'_{+}-z_{+}}{z'_{+}-z_{-}} \cdot\frac{z'_{-}-z_{-}}{z'_{-}-z_{+}} \cdot\frac{z-z_{+}}{z-z_{-}} \cdot\frac{z-z'_{+}}{z-z'_{-}}
E_{z_{+}}E_{z_{-}}^{-1}E_{z'_{+}}E_{z'_{-}}^{-1}E_{z},
\notag\end{equation}
here $ {\mathfrak y}=z_{+}\ominus z_{-} $, $ {\mathfrak y}'=z'_{+}\ominus z'_{-} $. Similarly, the action on $ E_{z}^{-1} $ is
\begin{equation}
E_{z}^{-1} \mapsto \frac{z'_{+}-z_{+}}{z'_{+}-z_{-}} \cdot\frac{z'_{-}-z_{-}}{z'_{-}-z_{+}} \cdot\frac{z-z_{-}}{z-z_{+}} \cdot\frac{z-z'_{-}}{z-z'_{+}}
E_{z_{+}}E_{z_{-}}^{-1}E_{z'_{+}}E_{z'_{-}}^{-1}E_{z}^{-1},
\notag\end{equation}
One can see that this conformal field is defined if
\begin{equation}
\frac{z'_{+}-z_{+}}{z'_{+}-z_{-}}\cdot\frac{z'_{-}-z_{-}}{z'_{-}-z_{+}},\qquad \frac{z-z_{-}}{z-z_{+}}\cdot\frac{z-z'_{-}}{z-z'_{+}},\qquad
\frac{z-z_{+}}{z-z_{-}}\cdot\frac{z-z'_{+}}{z-z'_{-}}
\notag\end{equation}
are not infinity, thus if
\begin{equation}
z'_{+}\not=z_{-},\quad z'_{-}\not=z_{+},\quad z\not=z_{-},\quad z\not=z_{+},\quad z\not=z'_{-},\quad z\not=z'_{+}.
\notag\end{equation}
Correspondingly, $ {\mathbit G}_{{\mathfrak y}}{\mathbit G}_{{\mathfrak y}'} $ is smooth on
\begin{equation}
z'_{+}\not=z_{-},\qquad z'_{-}\not=z_{+}.
\notag\end{equation}
Thus $ {\mathbit G}_{{\mathfrak y}}{\mathbit G}_{{\mathfrak y}'} $ has fusions at two divisor, both fusions are non-zero for $ n=1 $
only. Calculate the residue at $ z'_{+}=z_{-} $:
\begin{equation}
E_{z} \mapsto \left(z_{-}-z_{+}\right) \frac{z'_{-}-z_{-}}{z'_{-}-z_{+}} \cdot \frac{z-z_{+}}{z-z'_{-}} E_{z_{+}}E_{z'_{-}}^{-1}E_{z}.
\notag\end{equation}

This formula can be significantly simplified by a renormalization
\begin{equation}
{\mathcal G}\left(z_{+},z_{-}\right)=-\frac{1}{z_{+}-z_{-}}{\mathbit G}_{z_{+}\ominus z_{-}}.
\notag\end{equation}
Then the fusion of $ {\mathcal G}\left(z_{+},z_{-}\right) $ with $ {\mathcal G}\left(z'_{+},z'_{-}\right) $ at $ z'_{+}=z_{-} $ is given by
\begin{equation}
E_{z} \mapsto \frac{1}{z_{-}'-z_{+}}\cdot\frac{z-z_{+}}{z-z'_{-}} E_{z_{+}}E_{z'_{-}}^{-1}E_{z} = {\mathcal G}\left(z_{+},z_{-}'\right)E_{z}.
\notag\end{equation}
Note, however, that the conformal field $ {\mathcal G}\left(z_{+},z_{-}\right) $ is no longer smooth, it
has a pole at $ z_{+}=z_{-} $ with the residue $ {\mathbit G}_{z_{+}\ominus z_{+}}=-\operatorname{id} $.

One can conclude that the non-smoothness of the algebra generated by
conformal fields $ {\mathcal G}\left(z_{+},z_{-}\right) $ is completely describable by relations
\begin{align} \operatorname{Res}_{z_{+}=z_{-}}{\mathcal G}\left(z_{+},z_{-}\right) & = -\operatorname{id},
\notag\\
\operatorname{Res}_{z'_{+}=z_{-}} {\mathcal G}\left(z_{+},z_{-}\right){\mathcal G}\left(z'_{+},z'_{-}\right) & = {\mathcal G}\left(z_{+},z'_{-}\right).
\notag\\
\operatorname{Res}_{z_{+}=z'_{-}} {\mathcal G}\left(z_{+},z_{-}\right){\mathcal G}\left(z'_{+},z'_{-}\right) & = -{\mathcal G}\left(z'_{+},z_{-}\right).
\notag\end{align}
The last two relationship look similar to commutation relations of matrix
elements. The results of Section~\ref{s6.30} will show that this is
indeed true after taking Laurent coefficients in $ z_{+} $ and $ z_{-} $ at 0.

Note also that the conformal fields $ {\mathcal G}\left(z_{+},z_{-}\right) $ and $ {\mathcal G}\left(z'_{+},z'_{-}\right) $ commute:
\begin{equation}
{\mathcal G}\left(z_{+},z_{-}\right){\mathcal G}\left(z'_{+},z'_{-}\right) = {\mathcal G}\left(z'_{+},z'_{-}\right){\mathcal G}\left(z_{+},z_{-}\right)
\notag\end{equation}
on the subset of the set of parameters where both sides are defined. We
leave verification of this statement to the reader.

\begin{remark} Note that renormalization from $ {\mathbit G} $ to $ {\mathcal G} $ is not invariant w.r.t.~
automorphisms of $ {\mathbb P}^{1} $. To see the invariant form of fusion relation, note
that the residue describes the singular part of a conformal field near a
divisor.

Rewrite the above calculations as
\begin{equation}
{\mathbit G}_{z_{+}\ominus z_{-}}{\mathbit G}_{z'_{+}\ominus z'_{-}}= \frac{\frac{z'_{+}-z_{+}}{z_{+}-z'_{-}} }{ \frac{z_{-}-z'_{+}}{z'_{-}-z_{-}} } {\mathbit E}_{z_{-}}^{-1}{\mathbit E}_{z'_{+}}{\mathbit G}_{z_{+}\ominus z'_{-}}\text{,
}= \lambda\left(z'_{+},z'_{-},z_{+},z_{-}\right) {\mathbit E}_{z_{-}}^{-1}{\mathbit E}_{z'_{+}}{\mathbit G}_{z_{+}\ominus z'_{-}},
\notag\end{equation}
if $ z_{-}\approx z'_{+} $ this is
\begin{equation}
{\mathbit G}_{z_{+}\ominus z_{-}}{\mathbit G}_{z'_{+}\ominus z'_{-}}\approx \lambda\left(z'_{+},z'_{-},z_{+},z_{-}\right) {\mathbit G}_{z_{+}\ominus z'_{-}},
\notag\end{equation}
here $ \lambda\left(a,b,c,d\right) $ is the double ratio of 4 points on $ {\mathbb P}^{1} $.

One should understand the above equality in the following sense: for
any $ F\in{\mathfrak M}{\mathcal M}\left({\mathbb P}^{1}\right) $ the left-hand side and the right-hand side have same
residues at (an open subset of) $ z_{-}=z'_{+} $ when applied to $ F $, same for
$ z_{+}=z'_{-} $. Note that the double ratio has poles exactly at these two
divisors. \end{remark}

\begin{remark} Results Section~\ref{s6.30} will show that the Laurent coefficients
of the family $ {\mathcal G}\left(s,r\right) $ (and its smooth part $ {\mathcal G}\left(s,r\right)-\frac{1}{r-s} $) give rise to
standard descriptions of the Lie algebra $ {\mathfrak g}{\mathfrak l}\left(\infty\right) $ and its central extension
$ {\mathfrak g}\left(\infty\right) $. \end{remark}

\subsection{Further reductions }\label{s2.40}\myLabel{s2.40}\relax  Consider a further reduction of the family
$ {\mathbit E}_{{\mathfrak y}}M_{\xi} $. Restrict $ {\mathfrak h}=\left(z_{+},z_{-}\right) $ to have $ z_{-}=0 $, and consider only $ \xi $ which have a
pole at $ \infty $. The second restriction makes this family dominant on $ {\mathbb A}^{1}\subset{\mathbb P}^{1} $,
the first one codominant outside of 0 only.

In fact it is better to abandon the restriction that the homogeneity
degree is 0, and consider $ {\mathbit E}_{0}^{-1} $ as a normalization of degree only. Thus
consider the family $ \Phi\left(z_{+},z_{-}\right)= {\mathbit E}_{z_{+}}{\mathbit M}_{1/\left(z-z_{-}\right)} $. It is now dominant and
codominant smooth family of homogeneity degree 1 on $ {\mathbb A}^{1} $. Note that
\begin{equation}
{\mathbit M}_{1/\left(z-z_{-}\right)}{\mathbit E}_{z_{+}} = \frac{1}{z_{+}-z_{-}}{\mathbit E}_{z_{+}}{\mathbit M}_{1/\left(z-z_{-}\right)},
\notag\end{equation}
and one can easily calculate fusions of $ \Phi\left(z_{+},z_{-}\right) $ and $ \Phi\left(z'_{+},z'_{-}\right) $ using
this relationship.

Make a reduction $ z_{+}=z_{-} $ in $ \Phi\left(z_{+},z_{-}\right) $ (or consider the residue of
$ {\mathbit M}_{1/\left(z-z_{-}\right)}{\mathbit E}_{z_{+}} $):
\begin{equation}
\varphi\left(z_{0}\right)\buildrel{\text{def}}\over{=}\Phi\left(z_{0},z_{0}\right)={\mathbit E}_{z_{0}}{\mathbit M}_{1/\left(z-z_{0}\right)}.
\notag\end{equation}
It is still a dominant codominant smooth family on $ {\mathbb A}^{1} $. One can easily see
that $ \varphi\left(z\right)\varphi\left(z'\right) $ is not smooth, thus $ \varphi $ has a fusion with itself of order
$ n=1 $, which is proportional to $ {\mathbit E}_{z_{0}}^{2}{\mathbit M}_{1/\left(z-z_{0}\right)}^{2} $, and that
\begin{equation}
\varphi\left(z\right)\varphi\left(z'\right)=-\varphi\left(z'\right)\varphi\left(z\right).
\notag\end{equation}
One may also calculate fusions of $ \varphi $ with itself of higher order, and see
that one gets higher and higher degrees of $ {\mathbit E} $ and $ {\mathbit M} $. Thus the field $ \varphi $ is
a $ 1 $-parametric analogue of the field from Section~\ref{s12.20}.

Two other fields are defined by reductions of $ {\mathbit E}_{{\mathfrak y}}M_{\xi} $, $ {\mathfrak y}=z_{+}\ominus z_{-} $, with
either $ z_{+}=0 $, or $ \xi $ having a pole at $ \infty $. They are again dominant on $ {\mathbb A}^{1} $ only:
\begin{equation}
\psi\left(z_{0}\right)={\mathbit E}_{z_{0}}{\mathbit M}_{z-z_{0}},\qquad \psi^{+}\left(z_{0}\right)={\mathbit E}_{z_{0}}^{-1}{\mathbit M}_{1/\left(z-z_{0}\right)}.
\notag\end{equation}
These fields have the same commutation relations as $ \varphi $:
\begin{align} \psi\left(z\right)\psi\left(z'\right) & =-\psi\left(z'\right)\psi\left(z\right),
\notag\\
\psi\left(z\right)\psi^{+}\left(z'\right) & =-\psi^{+}\left(z'\right)\psi\left(z\right),
\notag\\
\psi^{+}\left(z\right)\psi^{+}\left(z'\right) & =-\psi^{+}\left(z'\right)\psi^{+}\left(z\right),
\notag\end{align}
moreover, the families $ \psi\left(z\right)\psi\left(z'\right) $ and $ \psi^{+}\left(z\right)\psi^{+}\left(z'\right) $ are smooth, and
$ \psi\left(z\right)\psi^{+}\left(z'\right) $ has the only fusion
\begin{equation}
\psi\left(z\right)\psi^{+}\left(z'\right) \approx \frac{1}{z'-z}\operatorname{id}
\notag\end{equation}
Thus the fields $ \psi $ and $ \psi^{+} $ are $ 1 $-parametric analogues of the field from
Section~\ref{s12.30}.

Let us spell out the calculation of the products $ \psi\left(r\right)\psi^{+}\left(s\right) $ and
$ \psi^{+}\left(s\right)\psi\left(r\right) $.

Suppose that $ F\in{\mathfrak M}{\mathcal M}\left({\mathbb A}^{1}\right) $, $ \xi\in{\mathcal M}\left({\mathbb A}^{1}\right) $, let $ z $ runs through $ {\mathbb A}^{1} $. Then
\begin{equation}
\left(\psi^{+}\left(s\right)F\right)\left(\xi\right) = \xi\left(r\right)\frac{1}{\left(s-r\right)\xi\left(s\right)}F\left(\frac{1}{z-s}\left(z-r\right)\xi\left(z\right)\right),
\label{equ12.31}\end{equation}\myLabel{equ12.31,}\relax 
and
\begin{align} \left(\psi\left(r\right)\psi^{+}\left(s\right)F\right)\left(\xi\right) & = \xi\left(r\right)\frac{1}{\left(s-r\right)\xi\left(s\right)}F\left(\frac{1}{z-s}\left(z-r\right)\xi\left(z\right)\right),
\notag\\
\left(\psi^{+}\left(s\right)\psi\left(r\right)F\right)\left(\xi\right) & = \frac{1}{\left(r-s\right)\xi\left(s\right)}\xi\left(r\right)F\left(\left(z-r\right)\frac{1}{z-s}\xi\left(z\right)\right),
\notag\end{align}
thus $ \psi\left(r\right)\psi^{+}\left(s\right) = -\psi^{+}\left(s\right)\psi\left(r\right) $, if $ r\not=s $. It is also easy to calculate the
pole of $ \psi\left(r\right)\psi^{+}\left(s\right) $ at $ r\approx s $ for a fixed $ F $:
\begin{equation}
\left(\psi\left(r\right)\psi^{+}\left(s\right)F\right)\left(\xi\right) = \frac{1}{\left(s-r\right)}F\left(\xi\right) + \left[\text{ a smooth function }\right].
\notag\end{equation}
Indeed, fix $ \xi\left(z\right)\in{\mathcal M}\left({\mathbb A}^{1}\right) $, consider $ F\left(\frac{z-r}{z-s}\xi\left(z\right)\right) $ as a function of $ r $ and
$ s $. Show that for a generic $ \xi $ this function is a well-defined meromorphic
function of $ r $ and $ s $. Moreover, it has no pole at $ r\approx s $.

It is easy to deduce this given an explicit formula for $ F $ as a
rational expression of $ E_{z}^{\bullet} $, $ z\in{\mathbb A}^{1} $, but one can also obtain it from the
heuristic for $ {\mathfrak M}{\mathcal M}\left(Z\right) $ as being compatible collections of meromorphic
functions on $ X_{{\mathcal L}} $, $ {\mathcal L} $ running through ample sheaves on $ Z $ (see Section~%
\ref{s11.10} and \cite{Zakh98Mer}).

As in Section~\ref{s11.10}, in the case $ Z={\mathbb A}^{1} $ one should consider
manifolds $ X_{{\mathcal O}\left(k\right)} $, $ k>0 $, and mappings between them given by multiplication
of numerator and denominator by $ P_{l}\left(z\right) $, $ P_{l}\left(z\right) $ being a polynomial in $ z $ of
degree $ l $. The ``function'' $ F $ is modeled by meromorphic functions $ F^{\left(k\right)} $ on
$ X_{{\mathcal O}\left(k\right)} $ compatible with mappings which correspond to a generic $ P $. In
particular, the divisor of poles of $ F^{\left(k+1\right)} $ does not contain the image of
$ X_{{\mathcal O}\left(k\right)} $ in $ X_{{\mathcal O}\left(k+1\right)} $ for $ P_{1}=x-s $ and a generic $ s $. This means that for generic
$ s $ the expression $ F\left(\frac{z-r}{z-s}\xi\left(z\right)\right) $ has no pole at $ r=s $ for a generic
$ \xi\in X_{{\mathcal O}\left(k\right)} $.

Hence
\begin{equation}
\psi\left(r\right)\psi^{+}\left(s\right) = \frac{1}{\left(s-r\right)} \boldsymbol1 + O\left(1\right),
\notag\end{equation}
here {\bf1 }is the identity operator. Thus the fields $ \psi $ and $ \psi^{+} $ are
$ 1 $-parametric analogues of the field from Section~\ref{s12.30}.

Similarly,
\begin{align} \varphi\left(r\right)\varphi\left(s\right) & = \frac{1}{\left(s-r\right)} \varphi_{2} + O\left(1\right),
\notag\\
\left(\varphi_{2}\left(r\right)F\right)\left(\xi\right) & = \xi^{2}\left(r\right)F\left(\frac{\xi\left(z\right)}{\left(z-r\right)^{2}}\right),\qquad \xi\text{ is rational.}
\notag\end{align}

One can interpret these formulae as saying that {\bf1 }is a {\em fusion\/} of $ \psi $
and $ \psi^{+} $, and $ \varphi_{2} $ is a {\em fusion\/} of $ \varphi $ with itself.

In Section~\ref{s16.60} we will show that $ \psi^{+} $ and $ \psi $ give rise to standard
fermion creation and destruction operators.

\subsection{Vertex operators }\label{s2.44}\myLabel{s2.44}\relax  Consider how the transition to the formal
completion $ \widehat{Z} $ of $ Z={\mathbb P}^{1} $ at 0 can be written in an explicit form. Recall that
in the case of a completed manifold the lattice disappears, thus
conformal fields become well-defined families of operators on the vector
space $ {\mathfrak M}{\mathcal M}\left(\widehat{Z}\right) $. As in the case of general $ Z $, $ {\mathfrak M}{\mathcal M}\left(\widehat{Z}\right) $ is generated by $ E_{0}^{\bullet} $, but
in $ 1 $-dimensional case the index $ \bullet $ takes values in $ {\mathbb Z}_{\geq0} $. (Similarly, in the
homogeneous situation one may normalize an element $ \xi\in{\mathcal M}\left(\widehat{Z}\right) $ to satisfy
$ \xi\left(0\right)=1 $, thus $ E_{0}\equiv 1 $, so the set of indices of $ E_{0}^{\bullet} $ can be restricted to be
$ {\mathbb N}={\mathbb Z}_{>0} $.) Denote these generators of $ {\mathfrak M}{\mathcal M}\left(\widehat{Z}\right) $ by $ e^{\left(k\right)} $, $ k\in{\mathbb Z}_{\geq0} $.

The action of $ {\mathbit M}_{\xi} $ is defined on $ {\mathfrak M}{\mathcal M}\left(\widehat{Z}\right) $ only if $ \xi $ has no pole or zero
at 0. Let $ \xi $ has a Taylor coefficients $ \xi_{k} $, $ \xi=\sum\xi_{k}z^{k} $, then $ {\mathbit M}_{\xi} $
acts as a field automorphism of $ {\mathfrak M}{\mathcal M}\left(\widehat{Z}\right) $ generated by
\begin{equation}
{\mathbit M}_{\xi}\colon e^{\left(k\right)} \mapsto \sum_{l=0}^{k} e^{\left(l\right)}\xi_{k-l}.
\notag\end{equation}
In the particular case $ \xi\left(z\right) = -s/\left(z-s\right) $, one gets
\begin{equation}
{\mathbit M}_{-s/\left(z-s\right)}\colon t_{k} \mapsto \sum_{l=0}^{k} s^{l-k}e^{\left(l\right)},
\notag\end{equation}

However, to define the action of operators $ {\mathbit E}_{z} $ ``outside'' of $ z=0 $ one
needs to consider formal power series, since only the operators $ {\mathbit E}_{0}^{\bullet} $ act
in $ {\mathfrak M}{\mathcal M}\left(\widehat{Z}\right) $. So consider
\begin{equation}
E_{z}=\sum_{k}\frac{1}{k!}z^{k}E_{z}^{\left(k\right)},\qquad {\mathbit E}_{z}=\sum_{k}\frac{1}{k!}z^{k}{\mathbit E}_{z}^{\left(k\right)}
\notag\end{equation}
as formal power series in $ z $, and see that the coefficients of the latter
series act in the field $ {\mathfrak M}{\mathcal M}\left(\widehat{Z}\right) $, $ {\mathbit E}_{z}^{\left(k\right)} $ acts as multiplication by $ e^{\left(k\right)} $.

Both these actions preserve the subalgebra $ {\mathcal V}_{\geq0}\subset{\mathfrak M}{\mathcal M}\left(\widehat{Z}\right) $ of polynomials
in $ e^{\left(k\right)} $. Assigning to each $ e^{\left(k\right)} $ grading 1, one obtains a grading
$ {\mathcal V}_{\geq0}=\bigoplus_{l\geq0}{\mathcal V}_{l} $. Allowing negative powers of $ e^{\left(0\right)} $, one may define $ {\mathcal V}_{l} $ with $ l<0 $.
Let $ {\mathcal V}=\bigoplus_{l\in{\mathbb Z}}{\mathcal V}_{l} $.

The usual way to obfuscate the simple formulae given above is to
choose a different set of generators in the algebra $ {\mathcal V} $. The above set of
generators $ e^{\left(k\right)} $ corresponds to taking Taylor coefficients $ e^{\left(k\right)} $ at 0 of an
element $ \xi\left(z\right) $ of $ {\mathcal M}\left({\mathbb P}^{1}\right) $ without a pole at 0. These Taylor coefficient
generate some field of functions on $ {\mathcal M}\left({\mathbb P}^{1}\right) $.

Up to this moment one could work in any characteristic, since though
the above formulae contain denominators, in fact $ \frac{e^{\left(k\right)}}{k!} $ is well
defined as a Taylor coefficient of given meromorphic function $ \xi $ at 0.

Now assume $ \operatorname{char}{\mathbb K}=0 $, and take different coordinates: for a fixed
$ \xi\in{\mathcal M}\left({\mathbb P}^{1}\right) $, $ \xi\left(0\right)\in{\mathbb K}^{\times} $, write
\begin{equation}
\xi\left(z\right)=T \exp  \sum_{k>0}t_{k}z^{k},\qquad T,t_{k}\in{\mathbb K}.
\notag\end{equation}
Considered as functions of $ \xi $, the symbols $ T $, $ t_{1},\dots $ are meromorphic
functions on $ {\mathcal M}\left({\mathbb P}^{1}\right) $ (since they can be expressed as polynomials in $ e^{\left(k\right)} $),
and $ e^{\left(k\right)} $ can be expressed as polynomials in $ T,t_{1},\dots . $ Thus $ T,t_{1},\dots $ form
a different collection of independent generators of $ {\mathcal V}_{\geq0} $. The degree of $ T $
is 1, of $ t_{k} $ is 0, $ {\mathcal V}_{l}=T^{l}{\mathbb C}\left[t_{1},\dots \right] $.

One can easily write the considered above operators in terms of
these new coordinates:

\begin{proposition} \label{prop12.30}\myLabel{prop12.30}\relax 
\begin{enumerate}
\item
The vector subspace $ {\mathcal V} $ is invariant w.r.t.~operators $ {\mathbit M}_{z-s} $, $ s\not=0 $. The
action of $ {\mathbit M}_{z-s} $ can be written as
\begin{equation}
{\mathbit M}_{z-s}T=-sT,\qquad {\mathbit M}_{z-s} t_{k}=t_{k}-\frac{1}{ks^{k}},\qquad {\mathbit M}_{z-s}\left(ab\right)=\left({\mathbit M}_{z-s}a\right)\left({\mathbit M}_{z-s}b\right).
\notag\end{equation}
\item
The vector subspace $ {\mathcal V} $ is invariant w.r.t.~operators $ {\mathbit E}_{s} $, $ s\approx0 $. In
other words, if one writes the Taylor series for $ {\mathbit E}_{s} $:
\begin{equation}
{\mathbit E}_{s}v \approx \sum_{l\geq0}\varepsilon_{k}\left(v\right)s^{l},
\notag\end{equation}
then the operators $ \varepsilon_{k} $ preserve the subspace $ {\mathcal V} $. The action of $ {\mathbit E}_{s} $ on $ {\mathcal V} $
can be described as
\begin{equation}
{\mathbit E}_{s}T=e^{\sum_{l}t_{l}s^{l}}T^{2},\qquad {\mathbit E}_{s}t_{k}=T e^{\sum_{l}t_{l}s^{l}}t_{k},\qquad {\mathbit E}_{s}\left(ab\right)=\left({\mathbit E}_{s}a\right)b=a\left({\mathbit E}_{s}b\right).
\notag\end{equation}

\end{enumerate}
\end{proposition}

Writing $ {\mathbit M}_{z-s}|_{{\mathcal V}_{l}} $, $ {\mathbit E}_{s}|_{{\mathcal V}_{l}} $ as Laurent series in $ s $, one gets
\begin{align} {\mathbit M}_{z-s} & = \left(-s\right)^{l}\exp \left(-\sum_{k>0}\frac{1}{ks^{k}}\frac{\partial}{\partial t_{k}}\right),
\notag\\
E_{s} & =T \exp  \sum_{k>0} s^{k}t_{k} 
\notag\end{align}
(recall that in the Section~\ref{h0.2} these operators were called $ a^{-1}\left(s\right) $,
$ b\left(s\right) $). Taking into account definitions of $ \varphi $, $ \psi $, and $ \psi^{+} $, one can conclude
that

\begin{corollary} \label{cor12.40}\myLabel{cor12.40}\relax  When considered as Laurent series in $ s $ and restricted
on the subspace $ {\mathcal V}_{l} $, the operators $ {\mathbit M}_{z-s} $, $ {\mathbit E}_{s} $, $ \varphi\left(s\right) $, $ \psi\left(s\right) $, and $ \psi^{+}\left(s\right) $ become
the standard vertex operators:
\begin{align} \varphi\left(s\right) & =T \left(-s\right)^{-l}\exp \left(\sum_{k>0} s^{k}t_{k}\right) \exp \left(\sum_{k>0}\frac{1}{ks^{k}}\frac{\partial}{\partial t_{k}}\right),
\notag\\
\psi\left(s\right) & =T \left(-s\right)^{l}\exp \left(\sum_{k>0} s^{k}t_{k}\right) \exp \left(-\sum_{k>0}\frac{1}{ks^{k}}\frac{\partial}{\partial t_{k}}\right),
\notag\\
\psi^{+}\left(s\right) & =T^{-1} \left(-s\right)^{-l}\exp \left(-\sum_{k>0} s^{k}t_{k}\right) \exp \left(\sum_{k>0}\frac{1}{ks^{k}}\frac{\partial}{\partial t_{k}}\right).
\notag\end{align}
\end{corollary}

\subsection{Analogy with flag space }\label{s12.6}\myLabel{s12.6}\relax  Instead of interpreting standard vertex
operators as restrictions of conformal fields to the space of local
functionals $ {\mathcal V} $, one can interpret smooth families $ {\mathbit M}_{\xi} $, $ {\mathbit E}_{z} $, $ \varphi $, $ \psi $, and $ \psi^{+} $ as
{\em geometric realizations\/} of the standard vertex operators. The relation of
the geometric realization to the realization by differential operators
(of infinite order) is parallel to the relation of geometric realization
of representations of semisimple groups to representations in the
category $ {\mathcal O} $.

Indeed, take a representation of $ G $ in sections of a bundle $ {\mathcal L} $ over
the flag space $ {\mathcal F} $. The Lie algebra $ {\mathfrak g} $ acts by vector fields $ v_{X} $ on $ {\mathcal F} $ (lifted to
vector fields on $ {\mathcal L} $), $ X\in{\mathfrak g} $. Take an open Schubert cell $ U $ in $ {\mathcal F} $. It is isomorphic
to an affine space with coordinates $ T_{k} $. One can trivialize the bundle
$ {\mathcal L}|_{U} $, so the vector fields $ v_{X} $ become differential operators of order 1 in
variables $ T_{k} $. These differential operators preserve the space of
polynomials $ k\left[T_{k}\right] $, and the resulting representation of $ {\mathfrak g} $ is in category
$ {\mathcal O} $.

One can easily see that the set $ {\mathcal M}_{0}\left({\mathbb P}^{1}\right) $ of rational functions without
a pole at 0 is an analogue of $ U $, $ t_{k} $ are analogues of $ T_{k} $. This makes the
conformal fields into analogues of exponents of elements of Lie algebra
(thus elements of the Lie group), and makes $ {\mathcal M}\left({\mathbb P}^{1}\right) $ into an analogue of the
flag space.

\section{Hilbert schemes and codominant conformal fields }\label{h13}\myLabel{h13}\relax 

In this section we show that appropriate projective bundles over the
Hilbert schemes of points in $ Z $ are providing important {\em models\/} of the
space $ {\mathfrak M}{\mathcal M}\left(Z\right) $, and prove that codominant families automatically satisfy
conformal associativity relations, and vertex commutation relations.

Blow-ups we consider are going to be blow-ups with a center in an
arbitrary closed subscheme, not necessarily blow-ups with centers in
smooth submanifolds.

\subsection{Big families }\label{s13.09}\myLabel{s13.09}\relax  Describe some families in $ {\mathfrak M}{\mathcal M}\left(Z\right) $ such that any
other family can be {\em almost\/} pushed through one of these families. First,
note that for any finite collection of points $ A=\left\{z_{1},\dots ,z_{m}\right\}\subset Z $ any
meromorphic function $ F $ on $ {\mathcal M}\left(Z\right) $ can be determined by the restriction to
the subset $ U_{A} $ of $ {\mathcal M}\left(Z\right) $ consisting of functions without a pole at $ z_{i}\in A $. Let
$ Z_{0} $ be any subscheme of $ Z $ of finite length with support at the points $ z_{i}\in A $.
Say that $ F $ can be pushed through $ Z_{0} $, if $ F\left(\xi\right) $, $ \xi\in U_{A} $, is determined by
$ \xi|_{Z_{0}} $.
Since any $ E_{z_{0}}^{v_{1}\dots v_{d}} $ can be pushed through a suitable finite subscheme of
$ Z $ with support at $ \left\{z_{0}\right\} $, any element $ F\in{\mathfrak M}{\mathcal M}\left(Z\right) $ can be pushed through a
finite subscheme. Obviously, there is a minimal such subscheme which will
be denoted $ \operatorname{Supp} F $.

Second, for any vector space $ V $ over $ {\mathbb K} $ the symmetric algebra $ S^{\bullet}V $ is
filtered by $ {\mathfrak P}_{k}V\buildrel{\text{def}}\over{=}\oplus_{i=0}^{k}S^{i}V $. The field of fractions of $ S^{\bullet}V $ is filtered by
$ \bar{{\mathcal Q}}_{k}V\buildrel{\text{def}}\over{=}{\mathfrak P}_{k}V/\left({\mathfrak P}_{k}V\smallsetminus\left\{0\right\}\right) $. Suppose $ \dim  V<\infty $, consider $ {\mathfrak P}_{k}V $ as an affine
algebraic variety over $ {\mathbb K} $, denote the projectivization $ {\mathbb P}\left({\mathfrak P}_{k}V\oplus{\mathfrak P}_{k}V\right) $ by $ {\mathcal Q}_{k}V $.
It is a finite dimensional projective space, it contains projective
subspaces $ {\mathcal Q}_{k}^{\infty}V $ and $ {\mathcal Q}_{k}^{0}V $ which corresponds to $ {\mathbb P}\left({\mathfrak P}_{k}V\oplus0\right) $ and $ {\mathbb P}\left(0\oplus{\mathfrak P}_{k}V\right) $, let
$ \mathring{{\mathcal Q}}_{k}V={\mathcal Q}_{k}V\smallsetminus{\mathcal Q}_{k}^{\infty}V $. Points of $ \mathring{{\mathcal Q}}_{k}V $ over $ {\mathbb K} $ map to elements of $ \bar{{\mathcal Q}}_{k}V $.

For a vector space $ W $ denote the dual space by $ W^{\vee} $, then $ \bar{{\mathcal Q}}_{k}W^{\vee} $ may be
identified with the set of rational functions on $ W $ of degree $ \leq k $. Let $ Z_{0} $ be
a finite subscheme of length $ m $ in a quasiprojective variety $ Z $. Then
$ {\mathfrak O}_{Z_{0}}\buildrel{\text{def}}\over{=}\Gamma\left(Z_{0},{\mathcal O}_{Z_{0}}\right) $ is an $ m $-dimensional vector space, and given an element
$ \varphi\in{\mathcal Q}_{k}{\mathfrak O}_{Z_{0}}^{\vee} $, one can associate to it an element of $ {\mathfrak M}{\mathcal M}\left(Z\right) $ given by $ {\mathcal M}\left(Z\right)\ni\xi \mapsto
\varphi\left(\xi|_{Z_{0}}\right) $ (defined if $ Z_{0} $ does not intersect divisor of poles of $ \xi $). This
defines a family of meromorphic functions on $ {\mathcal M}\left(Z\right) $ parameterized by $ {\mathcal Q}_{k}{\mathfrak O}_{Z_{0}}^{\vee} $.

\begin{example} \label{ex13.07}\myLabel{ex13.07}\relax  Given distinct points $ z_{i}\in Z $, $ 1\leq i\leq m $, one can associate to
them a subscheme $ Z_{0}=\cup_{i}\left\{z_{i}\right\} $. The vector space $ {\mathfrak O}_{Z_{0}}^{\vee} $ has a natural basis $ E_{z_{i}} $,
and $ {\mathcal Q}_{k}{\mathfrak O}_{Z_{0}}^{\vee} $ consists of rational functions of $ E_{z_{i}} $, $ 1\leq i\leq m $ of degree at most
$ k $. Let $ U\subset Z^{m} $ consists of $ m $-tuples of distinct points on $ Z $, then one can
associate a projective space $ {\mathcal Q}_{k}{\mathfrak O}_{Z_{0}}^{\vee} $ to any $ A\in U $. Let $ Q_{k,m} $ be the total
space of the corresponding projective bundle over $ U $. There is a natural
admissible family of meromorphic functions on $ {\mathcal M}\left(Z\right) $ parameterized by $ Q_{k,m} $.
\end{example}

To investigate what happens when some of the points $ z_{i} $ collide, note
that there is an action of the symmetric group $ {\mathfrak S}_{m} $ on $ Q_{k,m} $ which preserves
the family, so one should expect that the family may be pushed down to an
appropriate ``desingularization'' of the quotient by the group action. It
is indeed so if one takes the principal irreducible component $ Z^{\left[m\right]} $ of the
{\em Hilbert scheme\/} $ \operatorname{Hilb}^{m}\left(Z\right) $ of $ m $ points on $ Z $ as such a ``desingularization''.
Recall that points of $ \operatorname{Hilb}^{m}\left(Z\right) $ are finite subschemes of length $ m $ in $ Z $,
the principal irreducible component is one which contains $ m $-tuples of
distinct points\footnote{Note that if $ \dim  Z\leq2 $ and $ Z $ is smooth, the Hilbert scheme is irreducible
(and smooth).}. It is known that $ Z^{\left[m\right]} $ and $ \operatorname{Hilb}^{m}\left(Z\right) $ are projective
varieties if $ Z $ is, and are quasiprojective in general case of
quasiprojective $ Z $.

For a given $ Z_{0}\in Z^{\left[m\right]} $ the vector space $ {\mathfrak O}_{Z_{0}} $ is a fiber of a locally free
{\em tautological structure sheaf\/} $ {\mathfrak O} $ on $ Z^{\left[m\right]} $. Applying $ {\mathcal Q}_{k} $ to the dual sheaf,
one obtains a projective bundle $ {\mathfrak Q}_{k,m} $ over $ Z^{\left[m\right]} $. Denote by the same symbol
$ {\mathfrak Q}_{k,m} $ the total space of this bundle. It is a quasiprojective variety, and
is projective if $ Z $ is. Define $ \mathring{{\mathfrak Q}}_{k,m} $ similarly. Similarly, one can define
corresponding bundles over $ \operatorname{Hilb}^{m}\left(Z\right) $, denote them $ {\mathfrak Q}_{k,m}^{\operatorname{Hilb}} $ and $ \mathring{{\mathfrak Q}}_{k,m}^{\operatorname{Hilb}} $.

We have seen that there is a natural family of meromorphic functions
on $ {\mathcal M}\left(Z\right) $ parameterized by $ {\mathfrak Q}_{k,m} $. Denote this family by $ {\mathfrak F}_{k,m} $.

\begin{proposition} \label{prop13.15}\myLabel{prop13.15}\relax  The family $ {\mathfrak F}_{k,m} $ is strictly meromorphic on $ {\mathfrak Q}_{k,m} $ and
strictly smooth on $ \mathring{{\mathfrak Q}}_{k,m} $. \end{proposition}

\begin{proof} First, the variety $ Q_{k,m} $ from Example~\ref{ex13.07} is a finite
covering of an open subset of $ {\mathfrak Q}_{k,m} $. We have seen that the family becomes
admissible when pulled back to $ Q_{k,m} $.

To show that $ {\mathfrak F}_{k,m} $ induces a meromorphic function on $ {\mathfrak Q}_{k,m}\times X_{{\mathcal L}} $ for $ {\mathcal L}\geq{\mathcal L}_{0} $
note that for any finite-dimensional vector space $ V $ the evaluation
function $ \left(f,v\right) \mapsto f\left(x\right) $ is well-defined and meromorphic on $ {\mathcal Q}_{k}V^{\vee}\times V $. Thus
for any vector bundle $ B $ over $ S $ the evaluation function is well-defined on
$ {\mathcal Q}_{k}B^{\vee}\times_{S}B $ and is meromorphic.

The next step of the proof requires the following fact. Represent
the locally free sheaf $ {\mathfrak O} $ on $ Z^{\left[m\right]} $ as the sheaf of sections of a
$ m $-dimensional vector bundle $ B \xrightarrow[]{\pi} Z^{\left[m\right]} $. The fibers of this bundle have a
structure of algebra, and the multiplication mapping $ *\colon B\times_{Z^{\left[m\right]}}B \to B $ is
regular.

\begin{nwthrmii} One can find a regular scalar function $ \varphi $ on $ B $ and a regular
mapping $ i\colon B \to B $ over $ Z^{\left[m\right]} $ such that $ b\in B $ is invertible w.r.t.~* iff
$ \varphi\left(b\right)\not=0 $, and then $ b^{-1}=\varphi\left(b\right)^{-1}i\left(b\right) $. Here $ b \mapsto b^{-1} $ is taking the ring inverse
in fibers of $ B $. \end{nwthrmii}

\begin{proof} One may suppose that $ Z $ is affine. Fix a function $ f $ on $ Z $. Define
a function
\begin{equation}
\Sigma_{l}f\colon \left(z_{1},\dots ,z_{m}\right) \mapsto \sigma_{l}\left(f\left(z_{1}\right),\dots ,f\left(z_{m}\right)\right)
\notag\end{equation}
on $ Z^{m} $, here $ \sigma_{l} $ is the $ l $-th elementary symmetric function. This function
on $ Z^{m} $ is invariant w.r.t.~the action of the symmetric group $ {\mathfrak S}_{m} $ on $ Z^{m} $,
thus it induces a function on {\em Chow fraction\/} $ Z^{m}/{\mathfrak S}_{m} $ (compare \cite{Kap93Chow}).
Composing with the natural mapping $ Z^{\left[m\right]} \to Z^{m}/{\mathfrak S}_{m} $, one obtains a regular
function on $ Z^{\left[m\right]} $, which will be denoted by the same symbol $ \Sigma_{l}f\colon Z_{0} \mapsto
\Sigma_{l}f\left(Z_{0}\right) $.

If we fix a finite-dimensional vector space $ V $ of functions on $ Z $,
then $ \Sigma_{l}\colon \left(f,Z_{0}\right) \mapsto \Sigma_{l}f\left(Z_{0}\right) $ is a regular function on $ V\times Z^{\left[m\right]} $. On the other
hand, for any fixed $ Z_{0}\in Z^{\left[m\right]} $ there is a natural mapping of restriction
$ r_{Z_{0}}:
V \to {\mathfrak O}_{Z_{0}} $. Moreover, for a fixed $ Z_{00}\in Z^{\left[m\right]} $ one can find an $ m $-dimensional
space $ V $ such that the above mapping is an isomorphism of vector spaces
for $ Z_{0}\in U_{Z_{00}} $, $ U_{Z_{00}} $ being a neighborhood of $ Z_{00} $. Thus $ \Sigma_{l} $ can be lifted to a
regular function on $ \pi^{-1}U_{Z_{00}} $. It is clear that different choices of $ V $
result in the liftings which coincide on $ \pi^{-1}\left(U_{Z_{00}}\cap U\right) $, thus everywhere
($ U $ is from Example~\ref{ex13.07}).

This shows that $ \Sigma_{l} $ can be considered as a function on $ B $. Note that
if $ Z_{0} $ contains points $ z_{1},\dots ,z_{m} $ (counted with appropriate multiplicities),
and $ f\in{\mathfrak O}_{Z_{0}} $, then $ \Sigma_{m}f=f\left(z_{1}\right)\dots f\left(z_{m}\right) $, thus $ f\in{\mathfrak O}_{Z_{0}} $ is invertible in $ {\mathfrak O}_{Z_{0}} $ iff $ \Sigma_{m} $
does not vanish on $ \left(Z_{0},f\right) $. Moreover,
\begin{equation}
f^{*m}-\Sigma_{1}\left(Z_{0},f\right)f^{*m-1}+\dots \pm\Sigma_{m-1}\left(Z_{0},f\right)f\mp\Sigma_{m}\left(Z_{0},f\right)\boldsymbol1=0
\notag\end{equation}
if $ Z_{0}\in U $, $ f\in{\mathfrak O}_{Z_{0}} $, * denotes multiplication in $ B $, {\bf1 }the unit element of
$ {\mathfrak O}_{Z_{0}} $,
and $ \bullet^{*d} $ denotes taking $ d $-th power w.r.t.~multiplication in $ B $. Since the
left-hand side is a regular mapping from $ B $ to $ B $, this is true for any
$ Z_{0}\in Z^{\left[m\right]} $.

In other words,
\begin{equation}
\Sigma_{m}\left(Z_{0},f\right)\boldsymbol1= f*\left(\Sigma_{m-1}\left(Z_{0},f\right)\boldsymbol1-\Sigma_{m-e}\left(Z_{0},f\right)f+\dots \pm\Sigma_{1}\left(Z_{0},f\right)f^{*m-2}\mp f^{*m-1}\right).
\notag\end{equation}
Recall that for invertible $ f\in{\mathfrak O}_{Z_{0}} $ one has $ \Sigma_{m}\left(Z_{0},f\right)\not=0 $, thus
\begin{equation}
f^{-1} = \Sigma_{m}\left(Z_{0},f\right)^{-1}\left(\Sigma_{m-1}\left(Z_{0},f\right)\boldsymbol1-\Sigma_{m-e}\left(Z_{0},f\right)f+\dots \pm\Sigma_{1}\left(Z_{0},f\right)f^{*m-2}\mp f^{*m-1}\right).
\notag\end{equation}
\end{proof}

\begin{lemma} \label{lm13.21}\myLabel{lm13.21}\relax  The expression
\begin{equation}
{\mathcal E}\colon \left(F,n,d\right) \to F\left(n/d\right)
\notag\end{equation}
defines a meromorphic function $ {\mathcal E}_{B} $ on $ {\mathcal Q}_{k}B^{\vee}\times_{Z^{\left[m\right]}}B\times_{Z^{\left[m\right]}}B $. \end{lemma}

\begin{proof} Indeed, $ \left(n,d\right) \to n/d $ is a regular dominant mapping $ B\times B^{\operatorname{inv}} \to B $,
here $ B^{\operatorname{inv}} $ is an open subset of invertible elements of $ B $. Thus $ F\left(n/d\right) $
gives a meromorphic function on an open subset of $ {\mathcal Q}_{k}B^{\vee}\times_{Z^{\left[m\right]}}B\times_{Z^{\left[m\right]}}B $, thus
on $ {\mathcal Q}_{k}B^{\vee}\times_{Z^{\left[m\right]}}B\times_{Z^{\left[m\right]}}B $. \end{proof}

For a line bundle $ {\mathcal L} $ on $ Z $ and a finite subscheme $ Z_{0} $ of $ Z $ there is a
natural {\em trivialization\/} mapping $ \Gamma\left(Z,{\mathcal L}\right) \to {\mathfrak O}_{Z_{0}} $ defined up to multiplication
by an invertible element of $ {\mathfrak O}_{Z_{0}} $. Say that $ {\mathcal L} $ {\em exhausts\/} $ Z_{0} $ if this mapping is
surjective, say that $ {\mathcal L} $ is $ m $-ample if it exhausts all the finite
subschemes of length $ m $. Now Lemma~\ref{lm13.21} shows that

\begin{corollary} If $ {\mathcal L} $ is $ m $-ample, then $ {\mathfrak F}_{k,m} $ induces a meromorphic function
on $ {\mathfrak Q}_{k,m}\times X_{{\mathcal L}} $. \end{corollary}

\begin{amplification} If $ {\mathcal L} $ is $ m $-ample, then for any $ q\in\mathring{{\mathfrak Q}}_{k,m} $ the divisor of poles
of the above meromorphic function on $ {\mathfrak Q}_{k,m}\times X_{{\mathcal L}} $ does not contain $ \left\{q\right\}\times X_{{\mathcal L}} $. \end{amplification}

\begin{proof} It is enough to prove a similar statement about the meromorphic
function $ {\mathcal E}_{B} $ from Lemma~\ref{lm13.21}. Suppose that $ q=\left(Z_{0},q_{Z_{0}}\right)\in{\mathcal Q}_{k}B^{\vee} $ is such that
$ \left\{q\right\}\times{\mathfrak O}_{Z_{0}}\times{\mathfrak O}_{Z_{0}} $ is in the divisor of poles of $ {\mathcal E}_{B} $. In the notations of Lemma~%
\ref{lm13.21} put $ d=\boldsymbol1_{Z_{0}} $. One can see that $ q $ considered as a meromorphic function
on $ {\mathfrak O}_{Z_{0}} $ is $ \infty $, thus $ q $ is on the infinity $ {\mathcal Q}_{k}^{\infty}B^{\vee}={\mathcal Q}_{k}B^{\vee}\smallsetminus\mathring{{\mathcal Q}}_{k}B^{\vee} $ of $ {\mathcal Q}_{k}B^{\vee} $. \end{proof}

To finish the proof of Proposition~\ref{prop13.15}, note that

\begin{lemma} Let $ {\mathcal L}_{1} $ be a very ample bundle on $ Z $. Then any $ {\mathcal L}\geq{\mathcal L}_{1}^{m} $ is $ m $-ample. \end{lemma}

\begin{proof} It is enough to prove that $ {\mathcal O}\left(m\right) $ is $ m $-ample on $ {\mathbb P}^{l} $ for any $ l $.
This in turn follows from the (well-known and easy-to-prove, see
\cite{Iar77Pun})
fact that any ideal of $ {\mathbb K}\left[\left[t_{1},\dots ,t_{l}\right]\right] $ of codimension $ m $ contains $ {\mathfrak m}^{m} $, $ {\mathfrak m} $
being the maximal ideal of $ {\mathbb K}\left[\left[t_{1},\dots ,t_{l}\right]\right] $. \end{proof}

Thus the divisor of poles of $ {\mathfrak F}_{k,m} $ is flat over $ \mathring{{\mathfrak Q}}_{k,m} $ indeed. \end{proof}

\begin{remark} \label{rem13.72}\myLabel{rem13.72}\relax  Consider a blow-up $ \widetilde{{\mathcal Q}}_{k}B^{\vee} $ of $ {\mathcal Q}_{k}B^{\vee} $ in $ {\mathcal Q}_{k}^{\infty}B^{\vee} $. Let $ \widetilde{{\mathcal Q}}_{k}^{\infty}B^{\vee} $ be
the preimage of $ {\mathcal Q}_{k}^{\infty}B^{\vee} $ in $ \widetilde{{\mathcal Q}}_{k}B^{\vee} $. Then the family $ {\mathfrak F}_{k,m} $ has a {\em pole\/} at $ \widetilde{{\mathcal Q}}_{k}^{\infty}B^{\vee} $
in the following sense: for a fixed $ q_{0}\in\widetilde{{\mathcal Q}}_{k}^{\infty}B^{\vee} $ and any fixed scalar
meromorphic function $ \varphi $ on $ \widetilde{{\mathcal Q}}_{k}B^{\vee} $ which vanishes on $ \widetilde{{\mathcal Q}}_{k}^{\infty}B^{\vee} $, and has no pole
at $ q_{0} $, the expression $ \varphi{\mathfrak F}_{k,m} $ gives a strictly meromorphic family which is
strictly smooth in an appropriate neighborhood of $ q_{0} $. On the normal
bundle $ {\mathcal N} $ to $ \widetilde{{\mathcal Q}}_{k}^{\infty}B^{\vee} $ the family $ {\mathfrak F}_{k,m} $ induces a homogeneous family of
meromorphic functions of homogeneity degree $ -1 $.

In particular, when studying homogeneous theory of $ {\mathfrak M}{\mathcal M}\left(Z\right) $, the line
bundle $ {\mathcal N} $ on $ \widetilde{{\mathcal Q}}_{k}^{\infty}B^{\vee} $ is an appropriate substitution of $ {\mathfrak Q}_{k,m} $. It has an
advantage of being smooth on a complete variety $ \widetilde{{\mathcal Q}}_{k}^{\infty}B^{\vee} $ (after a twist by
$ {\mathcal N}^{\vee} $). \end{remark}

In fact Proposition~\ref{prop13.15} may be partially extended to the
other irreducible components of $ {\mathfrak Q}_{k,m}^{\operatorname{Hilb}} $, if one removes strictness
conditions, and allows blowups of these components. This follows from the
following

\begin{lemma} \label{lm13.28}\myLabel{lm13.28}\relax  Consider a mapping $ S \xrightarrow[]{\alpha} \operatorname{Hilb}^{m}\left(Z\right) $ with an irreducible
quasiprojective variety $ S $. There is $ M>0 $, a blow-up $ \widetilde{S} \xrightarrow[]{\pi} S $ and a mapping
$ \widetilde{S} \xrightarrow[]{\alpha'} Z^{\left[M\right]} $ such that $ \alpha\left(s\right)\subset\alpha'\left(\pi\left(s\right)\right) $ if $ s\in\widetilde{S} $, here $ \subset $ denotes inclusion of
finite subschemes of $ Z $. \end{lemma}

\begin{proof} The technique of Section~\ref{s13.3} will show that it is enough to
prove the statement for an open subset $ \widetilde{S}\subset S $, instead of a blow-up.
Associate to a finite subscheme $ Z_{0}\in\operatorname{Hilb}^{m}\left(Z\right) $ the number $ \operatorname{rk}\left(Z_{0}\right) $ of distinct
points of $ Z $ which support $ Z_{0} $ (i.e., the length of $ Z_{0,\text{red}} $). Let $ \widetilde{S} $ be the
subset of $ S $ where $ \operatorname{rk}\left(\alpha\left(s\right)\right) $ is maximal. Obviously, $ \widetilde{S} $ is open in $ S $. Let
$ {\mathfrak m}_{z} $ be the maximal ideal of a point $ z\in Z $, $ \left\{z\right\}^{k} $ be the closed subscheme
which correspond to the ideal $ {\mathfrak m}_{z}^{k} $, mult$ _{z}Z_{0} $ be the length of $ Z_{0}\cap\left\{z\right\}^{m} $. Define
$ \alpha'\left(s\right)=\cap_{z\in Z}{\mathfrak m}_{z}^{\text{mult}_{z}\left(Z_{0}\right)} $. Obviously, $ \alpha' $ gives a smooth mapping of $ \widetilde{S} $ to $ Z^{\left[M\right]} $
for an appropriate $ M $. \end{proof}

\subsection{Versality }\label{s13.3}\myLabel{s13.3}\relax  As usual, assume that $ Z $ is a quasiprojective variety.
Consider an admissible family of elements of $ {\mathfrak M}{\mathcal M}\left(Z\right) $ parameterized by $ S $. By
definition, there is a mapping $ S \to {\mathfrak Q}_{k,m} $ (for appropriate $ k,m $) such that
the family is induced from the family $ {\mathfrak F}_{k,m} $. For non-admissible families
one gets the following statement:

\begin{theorem} \label{th13.32}\myLabel{th13.32}\relax  For any meromorphic family $ F_{\bullet} $ parameterized by an
irreducible quasiprojective variety $ S $ there is a blow-up $ \widetilde{S} $ of $ S $ such that
the pull-back of the family to $ \widetilde{S} $ is induced from a family $ {\mathfrak F}_{k,m} $ on $ {\mathfrak Q}_{k,m} $.

If the family $ F_{\bullet} $ is smooth, the same is true if we substitute $ {\mathfrak Q}_{k,m} $
by $ \mathring{{\mathfrak Q}}_{k,m} $. \end{theorem}

\begin{proof} Assume for the start that $ F_{\bullet} $ is strictly meromorphic, and $ Z $ and
$ S $ are projective. There is an open subset $ S_{0} $ of $ S $ such that the family is
admissible on $ S_{0} $, thus there is a mapping $ \alpha\colon S_{0} \to {\mathfrak Q}_{k,m} $ for appropriate $ k $
and $ m $ such that $ F_{\bullet}|_{S_{0}} $ is induced from $ {\mathfrak F}_{k,m} $.

Let $ \alpha'=\alpha\times\operatorname{id}\colon S_{0} \to {\mathfrak Q}_{k,m}\times S_{0} $, $ S_{1} $ be the closure of $ \alpha'S_{0} $ in $ {\mathfrak Q}_{k,m}\times S_{0} $.
Then $ S_{1} $ is a projective variety, thus is a blow-up $ \widetilde{S}=S_{1} \to S $, and one can
pull back the mapping $ \alpha\colon S_{0} \to {\mathfrak Q}_{k,m} $ from $ S_{0} $ to $ \widetilde{S} $.

\begin{lemma} For any subvariety $ T $ of $ {\mathfrak Q}_{k,m} $ the family $ {\mathfrak F}_{k,m} $ induces a
canonically defined strictly meromorphic family $ {\mathfrak F}_{k,m}|_{T} $ on $ T $. It is
strictly smooth on $ T\smallsetminus{\mathfrak Q}_{k,m}^{\infty} $ and only there. \end{lemma}

\begin{proof} The only thing one needs to prove is that when restricting $ {\mathfrak F}_{k,m} $
to $ T $, one will never need to resolve 0/0. But we have seen that if $ {\mathcal L} $ is
$ m $-ample, the restriction of $ {\mathfrak F}_{k,m} $ to $ X_{{\mathcal L}}\subset{\mathcal M}\left(Z\right) $ is well-defined, and has a
pole only at $ {\mathfrak Q}_{k,m}^{\infty} $, and zero only at $ {\mathfrak Q}_{k,m}^{0} $. Since for any vector space $ V $
the projective subspaces $ {\mathcal Q}_{k}^{\infty}V $ and $ {\mathcal Q}_{k}^{0}V $ do not intersect, this proves the
lemma. \end{proof}

This gives two families of meromorphic functions on $ {\mathcal M}\left(Z\right) $
parameterized by $ \widetilde{S} $, one induced from $ S $, another one from $ {\mathfrak Q}_{k,m} $. They
coincide on the preimage of $ S_{0} $, thus everywhere. Minor modifications of
this proof work in the cases of quasiprojective $ Z $ and $ S $ too.

Suppose that the family $ F_{\bullet} $ is strictly smooth. Then the pull back to
$ \widetilde{S} $ is also strictly smooth, thus one can see that $ \alpha\widetilde{S} $ does not intersect with
$ {\mathfrak Q}_{k,m}^{\infty} $. This proves the theorem in the case of strictly smooth $ F_{\bullet} $.

Consider now the case of a smooth $ F_{\bullet} $. There is an inclusion $ S \to T $,
and a blow-up $ \pi\colon \widetilde{T} \to T $ with $ \widetilde{T} $ which maps to $ {\mathfrak Q}_{k,m} $. Again, taking
completions, one may assume that $ S $ and $ T $ are projective. Let $ S_{1} $ be the
reduction of the scheme $ \pi^{-1}S $. It is a variety that maps onto $ S $. Consider
any irreducible component $ S_{2} $ of $ S_{1} $ which maps onto $ S $. Since $ S_{2} $ is
projective, cutting with hyperplanes, one can find a subvariety $ S_{3} $ of $ S_{2} $
such that the projection $ S_{3} \to S $ is finite of degree $ d $ over the generic
point. Taking a finite covering of an open subset of $ S_{3} $, and taking a
closure, one obtains a projective mapping $ S_{4} \to S $ which is finite Galois
covering over an open subset of $ S $, a mapping from $ S_{4} $ to $ {\mathfrak Q}_{k,m} $, and an
action of a Galois group $ G=\left\{g_{1},\dots ,g_{d}\right\} $ on $ S_{4} $.

We know that the pull-back of $ F_{\bullet} $ to $ S_{4} $ maps through $ {\mathfrak F}_{k,m} $ via $ S_{4} \xrightarrow[]{\beta}
{\mathfrak Q}_{k,m} $. For any $ g\in G $ one can also consider a twisted map $ \beta^{g}=\beta\circ g $. Consider
\begin{equation}
\bar{\beta}\buildrel{\text{def}}\over{=}\beta^{g_{1}}+\dots +\beta^{g_{d}}.
\notag\end{equation}

Let us explain what is the meaning of $ + $ in the above expression.
First of all, let $ \gamma $ be the composition of $ \beta $ with the projection $ {\mathfrak Q}_{k,m} \to
Z^{\left[m\right]} $. For $ s\in S_{4} $ consider the finite subscheme $ \bar{\gamma}s =\gamma^{g_{1}}s\cup..\cup\gamma^{g_{d}}s $ of length
$ m\left(s\right) $. Existence of a flattening subdivision \cite{Mum66Lec} shows that for an
open subset $ S_{5} $ of $ S_{4} $ the length of $ \bar{\gamma}s $ does not depend on $ s\in S_{5} $ (in fact by
completeness of the Hilbert scheme it is $ M=\max _{s\in S_{4}}m\left(s\right) $), and $ \bar{\gamma}s $ depends
smoothly on $ s\in S_{5} $, i.e., $ \bar{\gamma} $ extends to a regular mapping $ S_{5} \to \operatorname{Hilb}^{M}Z $.
Lemma~\ref{lm13.28} shows that decreasing $ S_{5} $ and increasing $ M $ one can
construct a mapping $ S_{5} \to Z^{\left[M\right]} $ instead.

Note that for any surjection of vector spaces $ W \to V $ there is a
natural mapping $ {\mathcal Q}_{k}V \to {\mathcal Q}_{k}W $, thus each $ \beta^{g_{l}}s $, $ s\in S_{5} $, can be considered as an
element of $ {\mathcal Q}_{k}{\mathfrak O}_{\bar{\gamma}s}^{\vee} $. This defines mappings $ \widetilde{\beta}^{g_{l}}\colon S_{5} \to {\mathfrak Q}_{k,M} $. Since the
images of these mappings do not intersect $ {\mathcal Q}_{k}^{\infty}{\mathfrak O}_{\bar{\gamma}s}^{\vee} $, the sum
\begin{equation}
\widetilde{\beta}^{g_{1}}+\dots +\widetilde{\beta}^{g_{k}}
\notag\end{equation}
is a well-defined element of $ {\mathcal Q}_{dk}{\mathfrak O}_{\bar{\gamma}s}^{\vee} $, thus defines a mapping $ S_{5} \to {\mathfrak Q}_{dk,M} $.

Taking a closure again, one obtains a mapping $ \beta_{6}\colon S_{6} \to {\mathfrak Q}_{dk,M} $ with $ S_{6} $
being a projective variety with a mapping $ S_{6} \to S $ which is a finite
Galois extension over an open subset of $ S $. Moreover, the action of Galois
group extends to $ S_{6} $ and the mapping $ \beta_{6} $ is Galois-invariant, thus induces
a mapping $ \beta_{7}\colon S_{7}=S_{6}/G \to {\mathfrak Q}_{dk,M} $. Now $ S_{7} $ is a blow-up of $ S $, and the family
on $ S_{7} $ induced by $ \beta_{7} $ coincides with $ d $ times the family induced by $ S $.
Dividing by $ d $, we prove the theorem in the case of smooth $ F_{\bullet} $ and $ \operatorname{char}{\mathbb K}=0 $.

If $ \operatorname{char}{\mathbb K}=p\not=0 $, suppose that degree of $ S_{4}/S $ in generic point is $ d_{1}p^{r} $
with $ \left(d_{1},p\right)=1 $. Instead of taking an arithmetic mean, split $ S_{4}/S $ into two
extensions, $ S_{4}/S_{4}' $ of degree $ d_{1} $, $ S_{4}'/S $ of degree $ p^{r} $. For the extension of
degree $ d_{1} $ one can consider the arithmetic mean, as above, thus what
remains is to consider the case of an extension of degree $ p^{r} $.

In the case of such an extension consider the product of different
$ \beta^{g_{l}} $ instead of the sum. Now the above arguments show that the family $ F_{\bullet}^{p^{r}} $
is induced by $ \beta_{7} $.

\begin{lemma} The mapping $ F \mapsto F^{p} $ induces a mapping $ {\mathfrak Q}_{k,m} \xrightarrow[]{\rho} {\mathfrak Q}_{pk,m} $ with a
closed image, and is a bijection with this image. \end{lemma}

\begin{proof} Note that $ \operatorname{Supp} F_{s}^{p^{r}}=\operatorname{Supp} F $ for any $ s $, thus it is enough to show
this in the case of the mapping $ {\mathcal Q}_{k}V \to {\mathcal Q}_{pk}V $ for a finite-dimensional
vector space $ V $. In this case the lemma is obvious. \end{proof}

The lemma shows that $ F_{s} $ (for fixed any $ s\in S_{7} $) may be written as a
rational expression of appropriate $ E_{z_{i}\left(s\right)}^{\bullet} $. Moreover, the $ p^{r} $-th degree of
coefficients $ c_{l} $ of these expressions are regular functions of $ s\in S_{7} $. What
remains is to show that the coefficients $ c_{l} $ themselves are regular
functions of $ s\in S_{7} $. On the other hand, fix sufficiently many test
functions $ \xi_{i}\in{\mathcal M}\left(Z\right) $, consider meromorphic functions $ t_{i}=F_{s}\left(\xi_{i}\right) $ of $ s\in S_{7} $. Each
function $ t_{i} $ is a fraction of two linear combinations of $ c_{l} $, coefficients
of these linear combinations are meromorphic functions on $ S_{7} $. These
provides a system of homogeneous linear equations on $ c_{l} $ with coefficients
being meromorphic functions on $ S_{7} $. Obviously, it is possible to chose
sufficiently many $ \xi_{i} $ so that this system determines $ c_{l} $ up to
multiplication by a common meromorphic function on $ S_{7} $. Since for an
appropriate value of this meromorphic function $ c_{l}^{p^{r}} $ become regular, so
are $ c_{l} $. This finishes the proof for the smooth case.

Consider now the general case, when $ F_{\bullet} $ is meromorphic. Substituting
$ S $ by an appropriate blow-up, one may assume that $ F_{\bullet}=\varphi\cdot G_{\bullet} $, with $ \varphi $ being a
scalar meromorphic function on $ S $, and $ G_{\bullet} $ is induced by $ \beta\colon S \to {\mathfrak Q}_{k,m} $. The
mapping $ \left(\varphi,\beta\right) $ maps on open subset of $ S $ to $ {\mathbb P}^{1}\times{\mathfrak Q}_{k,m} $, thus, replacing $ S $ by
an appropriate blow-up, one can assume that $ \left(\varphi,\beta\right)\colon S \to {\mathbb P}^{1}\times{\mathfrak Q}_{k,m} $ is
regular. On the other hand, for any vector space $ V $ there is a natural
multiplication mapping
\begin{equation}
{\mathbb P}^{1}\times{\mathcal Q}_{k}V\smallsetminus\delta \to {\mathcal Q}_{k}V \colon \left(\left(\lambda_{1},\lambda_{2}\right),\left(n,d\right)\right) \mapsto \left(\lambda_{1}n,\lambda_{2}d\right),
\notag\end{equation}
here $ \delta=\left(\left\{0\right\}\times{\mathcal Q}_{k}^{\infty}V\right)\cup\left(\left\{\infty\right\}\times{\mathcal Q}_{k}^{0}V\right) $. Since $ \varphi\not\equiv 0 $, $ \varphi\not\equiv\infty $, one can see that there is a
mapping from an open subset $ S_{0} $ of $ S $ to $ {\mathfrak Q}_{k,m} $. Now the same arguments
with taking a closure and a blow-up finish the proof of the theorem.\end{proof}

\begin{remark} Note that the values of $ k $ and $ m $ in the theorem are not uniquely
determined. Moreover, for any $ k_{2}\geq k_{1} $, $ m_{2}>m_{1} $ there is a mapping from an
appropriate blow-up of $ {\mathfrak Q}_{k_{1},m_{1}} $ to $ {\mathfrak Q}_{k_{2},m_{2}} $. Indeed, fix a point $ z_{0}\in Z $. For
$ Z_{0}\in\left(Z\smallsetminus\left\{z_{0}\right\}\right)^{\left[m\right]} $ the scheme $ Z_{0}\cup\left\{z_{0}\right\} $ is finite of length $ m+1 $, and evaluation
at $ z_{0} $ gives a linear functional $ E_{z_{0}} $ on $ {\mathfrak O}_{Z_{0}\cup\left\{z_{0}\right\}} $. Thus the mapping
\begin{equation}
\left(Z_{0},n:d\right) \to \left(Z_{0}\cup\left\{z_{0}\right\},nE_{z_{0}}^{l}:dE_{z_{0}}^{l}\right)
\notag\end{equation}
maps an open subset of $ {\mathfrak Q}_{k,m} $ to $ {\mathfrak Q}_{k+l,m+1} $ if $ l\geq0 $. This mapping extends to a
mapping of an appropriate blowup of $ {\mathfrak Q}_{k,m} $ to $ {\mathfrak Q}_{k+l,m+1} $. \end{remark}

\subsection{$ V $-conformal operators } The description of families of elements of
$ {\mathfrak M}{\mathcal M}\left(Z\right) $ from the previous sections allows a partial description of
conformal operators. Recall that to each element of $ {\mathfrak M}{\mathcal M}\left(Z\right) $ one may
associate the support, which is a finite subscheme of $ Z $. Forgetting all
the information about support except multiplicities of
points in $ Z $, one obtains a {\em support\/} $ 0 $-{\em cycle}, which is an element of $ Z^{m}/{\mathfrak S}_{m} $
for an appropriate $ m\in{\mathbb N} $.

\begin{theorem} \label{th13.42}\myLabel{th13.42}\relax  Suppose that $ V $ is an affine open subset of a
quasiprojective variety $ Z $, and $ {\mathbit C} $ is a $ V $-conformal operator on $ {\mathfrak M}{\mathcal M}\left(Z\right) $ which
dominates an open dense subset $ U\subset Z $. Then there is a function $ M\colon {\mathbb N}^{2} \to {\mathbb N} $,
blow-ups $ Z^{\left(m,k\right)} $ of $ Z^{m}/{\mathfrak S}_{m} $, and regular mappings
\begin{equation}
\alpha_{m,k}\colon Z^{\left(m,k\right)} \to Z^{M\left(m,k\right)}/{\mathfrak S}_{M\left(m,k\right)}.\qquad m,k\geq1,
\notag\end{equation}
such that for any element $ F $ of $ {\mathfrak M}{\mathcal M}\left(Z\right) $ with the support $ 0 $-cycle $ \sum a_{i}z_{i} $ there
is $ k>0 $ such that the support $ 0 $-cycle of $ {\mathbit C}\left(F\right) $ is smaller than $ \alpha_{m,k}\left(\sum a_{i}z_{i}\right) $
(here $ m=\sum a_{i} $) if $ {\mathbit C}\left(F\right) $ is defined.

Here $ \alpha_{m,k}\left(\sum a_{i}z_{i}\right) $ means $ \alpha_{m,k}\left(\eta\right) $, $ \eta $ being an arbitrary preimage of
$ \sum a_{i}z_{i} $ in $ Z^{\left(m,k\right)} $. \end{theorem}

\begin{proof} Apply $ {\mathbit C} $ to the family $ {\mathfrak F}_{k,m} $ restricted to a big open subset
$ S\subset{\mathfrak Q}_{k,m} $. Let $ T=U^{\left[m\right]} $, and $ S=\pi^{-1}T $, here $ \pi $ is the projection of $ {\mathfrak Q}_{k,m} $ to $ Z^{\left[m\right]} $.
Then $ {\mathbit C}\left({\mathfrak F}_{k,m}\left(s\right)\right) $ is defined for any $ s\in S $, and $ {\mathbit C}\left({\mathfrak F}_{k,m}\left(s\right)\right) $ is a smooth family
in $ {\mathfrak M}{\mathcal M}\left(V\right) $.

By Theorem~\ref{th13.32}, one can take a blow-up $ \widetilde{S} $ of a finite covering
of $ S $ and a mapping $ \alpha\colon \widetilde{S} \to {\mathfrak Q}_{k',m'} $ which induces the family $ {\mathbit C}\left({\mathfrak F}_{k,m}\left(s\right)\right) $.
Compose $ \alpha $ with the projection $ {\mathfrak Q}_{k',m'} \to Z^{\left[m'\right]} $ and the Chow mapping $ Z^{\left[m'\right]}
\to Z^{m'}/{\mathfrak S}_{m'} $. By the above condition and Definitions~\ref{def11.25} and~%
\ref{def11.37}, one can pick up $ \alpha $ (decreasing $ m' $ if necessary) in such a way
that this composition $ \beta $ goes to an affine variety $ V^{m'}/{\mathfrak S}_{m'} $.

Consider fibers of the projection from $ \widetilde{S} $ to $ S $. They are projective
varieties, thus each of them maps to a finite subset of $ V^{m'}/{\mathfrak S}_{m'} $. Since a
generic fiber of $ \widetilde{S} $ over $ S $ is a point, $ \beta $ induces a mapping from the
normalization of $ S $ to $ V^{m'}/{\mathfrak S}_{m'} $. Though fibers of $ \pi\colon S \to U^{\left[m\right]} $ are not
projective varieties, they are isomorphic to $ {\mathbb P}^{2l-1}\smallsetminus{\mathbb P}^{l-1} $, $ l>1 $, thus
contain many projective spaces. This means that any such fiber goes into
a point when mapped to an arbitrary affine variety. Hence $ \beta $ induces a
mapping from a normalization of $ T=U^{\left[m\right]} $ to $ V^{m'}/{\mathfrak S}_{m'} $. Finally, fibers of
projection from $ U^{\left[m\right]} $ to $ U^{m}/{\mathfrak S}_{m} $ are projective, and a generic fiber is a
point, thus $ \beta $ induces a mapping from normalization of $ U^{m}/{\mathfrak S}_{m} $ to $ V^{m'}/{\mathfrak S}_{m'} $.
In other words, this is $ {\mathfrak S}_{n} $-invariant mapping from normalization of $ U^{m} $ to
$ V^{m'}/{\mathfrak S}_{m'} $, thus a mapping from $ \widetilde{U}^{m}/{\mathfrak S}_{m} $, $ \widetilde{U} $ being the normalization of $ U $.

Moreover, $ \widetilde{S} $ and the mappings $ \alpha $ and $ \beta $ can be completed in such a way
that $ \widetilde{S} $ maps onto $ {\mathfrak Q}_{k,m} $ and to $ {\mathfrak Q}_{k',m'} $. This completion induces a mapping
from a blow-up of $ \widetilde{Z}^{m}/{\mathfrak S}_{m} $ to $ Z^{m'}/{\mathfrak S}_{m'} $, here $ \widetilde{Z} $ is the normalization of $ Z $. On
the other hand, $ \widetilde{Z} $ is a blow-up of $ Z $, thus $ \widetilde{Z}^{m}/{\mathfrak S}_{m} $ is a blow-up of $ Z^{m}/{\mathfrak S}_{m} $.

Since any element of $ {\mathfrak M}{\mathcal M}\left(Z\right) $ is in an appropriate $ {\mathfrak Q}_{k,m} $, this finishes
the proof. \end{proof}

Though the following statement does not add a lot of contents to
Theorem~\ref{th13.42}, it will be very useful for study of codominant
families. For two effective $ 0 $-cycles $ {\mathcal A}=\sum a_{i}z_{i} $ and $ {\mathcal A}'=\sum a'_{j}z_{j} $ let $ {\mathcal A}\smallsetminus{\mathcal A}' $ be
the sum of terms $ a_{i}z_{i} $ of $ {\mathcal A} $ such that $ z_{i} $ does not appear in $ {\mathcal A}' $.

\begin{amplification} \label{amp13.42}\myLabel{amp13.42}\relax  Suppose that $ V $ is an affine open subset of $ Z $, and
$ {\mathbit C} $ is a $ V $-conformal operator on $ {\mathfrak M}{\mathcal M}\left(Z\right) $ which dominates $ U\subset Z $. Then there are
functions $ M,L\colon {\mathbb N}^{2} \to {\mathbb N} $, blow-ups $ Z^{\left(m,k\right)} \xrightarrow[]{\pi_{m,k}} Z^{m}/{\mathfrak S}_{m} $, and regular mappings
\begin{equation}
\alpha_{m,k}\colon Z^{m,k} \to Z^{M\left(m,k\right)}/{\mathfrak S}_{M\left(m,k\right)}.\qquad m,k\geq1,
\notag\end{equation}
such that for any element $ F $ of $ {\mathfrak M}{\mathcal M}\left(Z\right) $ with the support $ 0 $-cycle $ \sum a_{i}z_{i} $ there
is $ k>0 $ such that the support $ 0 $-cycle of $ {\mathbit C}\left(F\right) $ is smaller than
\begin{equation}
\alpha_{m,k}\left(\sum a_{i}z_{i}\right) + L\left(m,k\right)\cdot\left(\sum a_{i}z_{i}\right)
\notag\end{equation}
(here $ m=\sum a_{i} $) if $ \left\{z_{i}\right\}\subset U $. Here $ \alpha_{m,k}\left(\sum a_{i}z_{i}\right) $ means $ \alpha_{m,k}\left(\eta\right) $, $ \eta $ being an
arbitrary preimage of $ \sum a_{i}z_{i} $ in $ Z^{\left(m,k\right)} $.

Moreover, let $ U^{m,k} $ be the part of $ Z^{\left(m,k\right)} $ which lies over $ U^{m}/{\mathfrak S}_{m} $. Then
for any $ u\in U^{m,k} $
\begin{equation}
\alpha_{m,k}\left(u\right) \smallsetminus \pi_{m,k}\left(u\right) \in V^{L}/{\mathfrak S}_{L}\qquad \text{for a suitable }L\geq0,
\notag\end{equation}
and for generic $ u $ supports of $ \alpha_{m,k}\left(u\right) $ and $ \pi_{m,k}\left(u\right) $ do not intersect. \end{amplification}

\begin{proof} For a cycle $ {\mathcal A}=\sum a_{i}z_{i} $ denote $ \sum a_{i} $ by $ |{\mathcal A}| $. One may assume that $ Z $ is
irreducible, since an appropriate blow-up of $ Z $ is. Obviously,

\begin{lemma} Suppose that $ Z $ is quasiprojective. Consider two regular mappings
\begin{equation}
\alpha\colon T \to Z^{M}/{\mathfrak S}_{M},\qquad \alpha'\colon T \to Z^{N}/{\mathfrak S}_{N}
\notag\end{equation}
with $ T $ being an irreducible variety. Then $ |\alpha\left(t\right)\smallsetminus\alpha'\left(t\right)| $ does not depend on
$ t\in U $ for an appropriate open dense subset $ U\subset T $. Moreover,
\begin{equation}
t \mapsto \alpha\left(t\right)\smallsetminus\alpha'\left(t\right)
\notag\end{equation}
gives a regular mapping $ \beta\colon U \to Z^{L}/{\mathfrak S}_{L} $ for an appropriate $ L $, thus for an
appropriate blow-up $ \widetilde{T} \xrightarrow[]{\pi}T $ this mapping may be extended to $ \beta\colon \widetilde{T} \to Z^{L}/{\mathfrak S}_{L} $.

Moreover, for any $ t\in\widetilde{T} $ one has $ \beta\left(t\right)\leq\alpha\left(\pi\left(t\right)\right) $. \end{lemma}

Consider $ \alpha $ from Theorem~\ref{th13.42}, let $ \alpha' $ be the projection of $ Z^{\left(m,k\right)} $
to $ Z^{m}/{\mathfrak S}_{m} $. Application of the previous lemma proves the first statement of
the amplification. To prove the second statement, note that $ V $-conformal
operator sends an element of $ {\mathfrak M}{\mathcal M}\left(Z\right) $ with support at $ A=\left\{z_{1},\dots ,z_{l}\right\}\subset Z $ to an
element with support in $ A\cup V $. This and Definition~\ref{def11.25} show that one can
assume that $ \operatorname{Supp} \alpha_{m,k}\left(u\right) $ is contained in $ \operatorname{Supp}\pi_{m,k}\left(u\right)\cup V $ (possibly after
decreasing $ L $). \end{proof}

\subsection{Codominant families } While the above description of $ V $-conformal
operators gives a lot of information, the description of codominant
families is yet more rigid:

\begin{theorem} \label{th13.52}\myLabel{th13.52}\relax  Consider quasiprojective varieties $ Z $ and $ S $, and a
codominant conformal family $ {\mathbit C}_{s} $ on $ S $ acting in $ {\mathfrak M}{\mathcal M}\left(Z\right) $. Then for any $ l\in{\mathbb N} $
there are $ M_{l},N_{l}\in{\mathbb N} $ and a mapping $ \alpha_{l}\colon S \to Z^{M_{l}}/{\mathfrak S}_{M_{l}} $ such that the mapping
$ \alpha_{m,k} $ from Theorem~\ref{th13.42} can be taken to be
\begin{equation}
\sum a_{i}z_{i} \mapsto N_{l}\sum a_{i}z_{i} + \alpha_{k}\left(s\right).
\notag\end{equation}
with $ \alpha_{k} $ depending only on $ s\in S $, but not on $ z_{i} $. \end{theorem}

\begin{proof} Suppose that $ {\mathbit C}_{\bullet} $ is defined on $ U\subset S\times Z $. Define $ U^{\left\{m\right\}}\subset S\times Z^{m}/{\mathfrak S}_{m} $ by the
condition that the fiber of $ U^{\left\{m\right\}} $ over any $ s\in S $ is the $ m $-th Chow power of
the fiber of $ U $ over $ s $:
\begin{equation}
U^{\left\{m\right\}}\cap\left(\left\{s\right\}\times Z^{m}/{\mathfrak S}_{m}\right) = \left\{s\right\}\times\left(\left(U\cap\left(\left\{s\right\}\times Z\right)\right)^{m}/{\mathfrak S}_{m}\right).
\notag\end{equation}
Obviously, Amplification~\ref{amp13.42} can be generalized to families of
conformal operators. This generalization provides a blow-up $ \widetilde{U} \xrightarrow[]{\pi} U^{\left\{m\right\}} $ and
a mapping $ \beta $ from $ \widetilde{U} $ to $ Z^{M}/{\mathfrak S}_{M} $ for an appropriate $ M>0 $. For a generic
$ u\in\widetilde{U} $ the supports of $ \beta\left(u\right) $ and $ \pi\left(u\right) $ do not intersect, for any $ u $ and any
open $ V\subset Z $ the support of $ \beta\left(u\right)\smallsetminus\pi\left(u\right) $ is in $ V $ as far as the $ S $-projection of
$ u $ is in an appropriate non-empty open dense subset $ S_{V}\subset S $.

We want to show that $ \beta\left(u\right) $ depends on $ S $-projection of $ u\in S\times Z^{m}/{\mathfrak S}_{m} $ only.
To show this in the case $ \dim  Z\geq2 $, it is enough to prove

\begin{proposition} Suppose that $ \dim  Z\geq2 $. For any open dense $ U\subset Z $ one can find
an open dense $ V\subset Z $ such that if a regular mapping $ b $ from any blow-up $ \widehat{U} \xrightarrow[]{p}
U^{m}/{\mathfrak S}_{m} $ to $ Z^{M}/{\mathfrak S}_{M} $ satisfies
\begin{equation}
\operatorname{Supp}\left(b\left(u\right)\smallsetminus p\left(u\right)\right)\subset V\qquad \text{for any }u\in\widehat{U}
\notag\end{equation}
and
\begin{equation}
\operatorname{Supp}\left(b\left(u\right)\right) \cap \operatorname{Supp} \left(p\left(u\right)\right)=\varnothing\qquad \text{for a generic }u\in\widehat{U},
\notag\end{equation}
then $ b $ is a constant mapping to $ V^{M}/{\mathfrak S}_{M} $. \end{proposition}

\begin{proof} One may suppose that $ Z $ is irreducible and affine, $ Z\subset{\mathbb A}^{L} $, $ \dim 
Z=d $, $ U=Z $. Note that for any fixed $ Z_{0}\in Z^{m-1}/{\mathfrak S}_{m-1} $ and any curve $ \Gamma\subset Z $ the
mapping $ z \mapsto z+Z_{0} $, $ z\in\Gamma $, sends $ \Gamma $ to $ Z^{m}/{\mathfrak S}_{m} $. This mappings lifts to any
blow-up of $ Z^{m}/{\mathfrak S}_{m} $. Denote the corresponding curve in such a blowup $ \Gamma\vee Z_{0} $.

Choose any $ d-1 $ coordinates $ \eta_{i_{1}},\dots ,\eta_{i_{d-1}} $ on $ {\mathbb A}^{L} $, consider the
corresponding projection $ Z \xrightarrow[]{{\text H}} {\mathbb A}^{d-1} $. Suppose that $ {\text H} $ is dominant. For an
open subset $ W\subset{\mathbb A}^{d-1} $ the preimages $ \Gamma_{w}={\text H}^{-1}\left(w\right) $, $ w\in W $, are curves with the same
number $ N $ of points ``at infinity'', i.e., points in $ \bar{\Gamma}_{w}\smallsetminus\Gamma_{w} $, $ \bar{\Gamma}_{w} $ being the
(maximal) completion of $ \Gamma_{w} $. Consider the mapping $ \Gamma_{w} \to Z^{M}/{\mathfrak S}_{M} $ induced by
the mapping $ b|_{\Gamma_{w}\vee Z_{0}} $. Note that $ \widetilde{\Gamma} $ has no more than $ NM $ points at infinity.
In particular, any regular mapping from $ \widetilde{\Gamma} $ to $ {\mathbb A}^{1}\smallsetminus\sigma $ is constant if $ \sigma $
contains at least $ NM $ different points.

Suppose than no coordinate function $ \eta_{j} $ is constant on $ Z $. Fix an
arbitrary subset $ \sigma $ of $ NM+m $ points on $ {\mathbb A}^{1} $, let $ V=Z\smallsetminus\bigcup_{j}\eta_{j}^{-1}\sigma $. Then any
mapping from $ \widetilde{\Gamma}_{w} $ to $ V\cup Z' $ is constant if any $ \eta_{j} $-projection of $ Z' $ contains
at most $ m $ points. Consider the function $ \varepsilon_{k}=\eta_{i_{k}}\circ b' $ restricted to $ \widetilde{\Gamma}_{w} $,
$ 1\leq k\leq d-1 $. The first condition of the proposition shows that $ \varepsilon_{k} $ sends $ \widetilde{\Gamma}_{w} $ to
$ \left({\mathbb A}^{1}\smallsetminus\sigma\right)\cup\varepsilon_{k}\left(\widetilde{\Gamma}_{w}\right)\cup\varepsilon_{k}\left(Z_{0}\right) $. Since $ \eta_{i_{k}}|_{\widetilde{\Gamma}_{w}} $ is constant, $ \varepsilon_{k}\left(\widetilde{\Gamma}_{w}\right)\cup\varepsilon_{k}\left(Z_{0}\right) $ contains at
most $ m $ points, thus $ \varepsilon_{k} $ is constant on $ \widetilde{\Gamma}_{w} $.

This shows that the composition of $ b $ with the projection of $ Z^{M}/{\mathfrak S}_{M} $ to
$ \left({\mathbb A}^{d-1}\right)^{M}/{\mathfrak S}_{M} $ is constant on any $ \Gamma_{w}\vee Z_{0} $ for any $ Z_{0} $ and $ w\in W $. Note that one can
do the same with all $ \left(d-1\right) $-tuples of independent coordinates on $ Z $ (one
may need to increase $ N $ to do it).

For a generic $ Z_{0} $ and $ z $ the preimage $ p^{-1}\left(z+Z_{0}\right) $ contains one point
only. The above observation applied to all independent $ \left(d-1\right) $-tuples of
coordinates shows that for any coordinate $ \eta_{j} $ on $ {\mathbb A}^{L} $ the composition
$ \widehat{\eta}_{j}\circ b\circ p^{-1}\left(z+Z_{0}\right) $ (here $ \widehat{\eta}_{j} $ is the extension of $ \eta_{j} $-projection to a mapping of
$ Z^{M}/{\mathfrak S}_{M} $ to $ \left({\mathbb A}^{1}\right)^{M}/{\mathfrak S}_{M} $) depends only on $ \eta_{j}\left(z\right) $ and $ Z_{0} $ for generic $ z $ and $ Z_{0} $.
Denote by $ \zeta_{j} $ (or $ \zeta_{j,Z_{0}} $) the corresponding mapping
\begin{equation}
{\mathbb A}^{1}\supset W_{j} \to \left({\mathbb A}^{1}\right)^{M}/{\mathfrak S}_{M}\colon \eta_{j}\left(z\right) \mapsto \eta_{j}\left(b\left(z+Z_{0}\right)\right)
\notag\end{equation}
defined on an appropriate open $ W_{j}=W_{j,Z_{0}} $.

Let $ \eta_{j} $, $ \eta_{j'} $ be two independent coordinates on $ Z\subset{\mathbb A}^{L} $ (this is where we
use the condition that $ \dim  Z>1 $). We may assume that $ \eta_{j}+\eta_{j'} $ coincides with
another coordinate $ \eta_{j''} $ on $ {\mathbb A}^{L} $ when restricted on $ Z $. If we consider $ \zeta_{j} $,
$ \zeta_{j'} $, $ \zeta_{j''} $ as multivalued functions on open subsets of $ {\mathbb A}^{1} $, the above
remark shows that $ \zeta_{j}\left(a\right)+\zeta_{j'}\left(b\right) $ is a (multivalued) function of $ a+b $. As in
the case of usual functions, this shows that all the branches of $ \zeta_{j} $ and
$ \zeta_{j'} $ are parallel translations of each other:

\begin{lemma} Consider curves $ \Gamma,\Gamma'\subset{\mathbb A}^{2} $. Denote coordinates on $ {\mathbb A}^{2}\times{\mathbb A}^{2} $ by
$ x_{1},x_{2},y_{1},y_{2} $. Consider the projection $ \pi\colon {\mathbb A}^{2}\times{\mathbb A}^{2} \to {\mathbb A}^{2} $ given by
$ \left(x_{1}+y_{1},x_{2}+y_{2}\right) $. If $ \pi|_{\Gamma_{1}\times\Gamma_{2}} $ is not dominant, there is an additive mapping
$ \rho\colon {\mathbb A}^{1} \to {\mathbb A}^{1} $ such that $ \Gamma $ and $ \Gamma' $ are unions of graphs of $ \rho+\operatorname{const} $. \end{lemma}

\begin{proof} All the tangent lines to $ \Gamma $ and $ \Gamma' $ have the same slope $ t $. Doing
a coordinate change $ x_{2}'=x_{2}-tx_{1} $, one may assume that this slope is 0.
Hence the equations of $ \Gamma $ can be written as $ F\left(x_{1}^{p},x_{2}\right)=0 $, $ p=\operatorname{char}{\mathbb K} $. If $ p\not=0 $,
and $ F\left(s,t\right) $ depends on $ s $, do a coordinate change $ x_{1}'=x_{1}^{p} $ and repeat the
argument enough times. \end{proof}

In general there may be several components of the graphs $ \Gamma $ and $ \Gamma' $ of
$ \zeta_{i} $ and $ \zeta_{j} $, so we may need to apply the lemma to some components of $ \Gamma\times\Gamma' $
only. Thus $ \zeta_{i} $ and $ \zeta_{j} $ may be unions of several parts, each of these parts
being described by the lemma. Thus $ \zeta_{i} $, $ \zeta_{j} $ may be unions of graphs of
$ \rho_{t}+\operatorname{const} $, with the same collection of additive functions $ \rho_{t} $ for both $ \zeta_{i} $
and $ \zeta_{j} $.

Let $ k=2 $ if $ \operatorname{char}{\mathbb K}\not=2 $, otherwise $ k=3 $. One may suppose that $ \eta_{j}^{k} $
coincides on $ Z $ with another coordinate $ \eta_{j'''} $ on $ {\mathbb A}^{L} $ when restricted on $ Z $.
Applying the lemma to $ \eta_{j'''} $ and $ \eta_{j'} $ shows that the additive functions $ \rho_{t} $
which correspond to the pair $ \left(\eta_{j},\eta_{j'}\right) $ commute with taking $ k $-th power.
Thus for a fixed $ t $ either $ \rho_{t}=0 $, or $ \rho_{t}\left(X\right)=X^{l} $, here $ l=p^{a} $, $ a\geq0 $, if
$ \operatorname{char}{\mathbb K}=p\not=0 $, and $ l=1 $ if $ \operatorname{char}{\mathbb K}=0 $.

If $ \operatorname{char}{\mathbb K}\not=0 $, and $ {\mathbb K} $ contains a transcendental $ \lambda $, then, as above, one
can conclude that $ \rho_{t} $ commute with multiplication by $ \lambda $. Thus $ a=0 $ in the
above case, or otherwise $ {\mathbb K} $ is the algebraic closure of $ {\mathbb F}_{p} $, and $ \rho $ is a
power of Frobenius morphism $ \Phi $.

If $ \rho=0 $, then the corresponding part of $ \zeta $ is a constant mapping $ {\mathbb A}^{1} \to
\left({\mathbb A}^{1}\right)^{M_{1}}/{\mathfrak S}_{M_{1}} $. If $ \rho=\operatorname{id} $, then the corresponding part of $ \zeta $ is constant, thus
$ bp^{-1}\left(z+Z_{0}\right) $ contains a multiple of $ z $, what contradicts the second
condition of the proposition, since $ Z_{0} $ is arbitrary. One can see that in
general $ bp^{-1}\left(z+Z_{0}\right) $ (for a fixed $ Z_{0} $) is some constant term plus
$ \sum_{a=0}^{A}c_{a}\Phi^{a}\left(z\right) $.

The constant term in this expression may depend on $ Z_{0} $, breaking $ Z_{0} $
into a sum $ z'+Z_{1} $ and repeating above arguments one can see that $ b $ is
\begin{equation}
\sum_{a=0}^{A}c_{a}\Phi^{a}p+C,\qquad c_{a}\geq0,\quad C\in Z^{M'}/{\mathfrak S}_{M'}\quad \sum_{a}c_{a}+M'=M.
\notag\end{equation}
Use $ \dim  Z\geq2 $ again. Taking a projection on $ {\mathbb A}^{2} $, it is enough to consider
the case $ Z\subset{\mathbb A}^{2} $. One may suppose that $ Z\cap\left({\mathbb A}^{1}\times\left\{0\right\}\right)\not=\varnothing $. Let $ V=Z\smallsetminus\left({\mathbb A}^{1}\times\left\{0\right\}\right) $.
Obviously, if the above $ b $ satisfies the conditions of proposition with
this $ V $, then $ c_{a}=0 $ for any $ a $. This finishes the proof of the proposition.
\end{proof}

Consider the case $ \dim  Z=1 $. Decreasing $ U $, one may assume that $ Z $ is
projective and normal. Fix $ Z_{0} $ again, restrict $ \beta $ to elements of the form
$ \left(s,z+Z_{0}\right) $, $ z\in Z $. This defines a mapping $ \beta_{Z_{0}} $ from an open subset of $ S\times Z $ to
$ Z^{M}/{\mathfrak S}_{M} $, replacing $ S $ by an open subset, one may extend it to a mapping of
$ S\times Z $. Note that to a mapping $ Z \xrightarrow[]{\gamma} Z^{M}/{\mathfrak S}_{M} $ one may associate a curve $ \Gamma\subset Z\times Z $
given by the condition that $ \gamma\left(z_{0}\right) $ is the intersection of $ \left\{z_{0}\right\}\times Z $ with $ \Gamma $
(it is possible that one may need to assign multiplicities to components
of $ \Gamma $). This gives a family $ \Gamma_{s} $ of curves in $ Z\times Z $ parameterized by $ S $. All
these curves have the same degree of projection to the second component,
call it $ L $. To show that $ \beta $ depends on $ s $ only, it is enough to show that
$ L=0 $.

Suppose that $ L>0 $. Taking a finite mapping of $ Z $ to $ {\mathbb P}^{1} $, and
multiplying $ M $ by degree, one can see that it is enough to consider $ Z={\mathbb P}^{1} $.
Fix a generic $ s_{0}\in S $, consider $ \Gamma=\Gamma_{s_{0}}\subset{\mathbb P}^{1}\times{\mathbb P}^{1} $. We know that $ \Gamma $ does not contain
the diagonal as an irreducible component. Consider two projections $ \pi_{1,2} $
of $ \Gamma $ to $ {\mathbb P}^{1} $ of degrees $ M $ and $ L $ as meromorphic functions on $ \Gamma $. Then $ \pi_{1}-\pi_{2} $
is a correctly defined meromorphic function on $ \Gamma $ which is not zero on any
irreducible component of $ \Gamma $, thus has no more than $ L+M $ zeros. Suppose $ Z\smallsetminus U $
contains $ d $ point, let $ V_{0} $ be $ {\mathbb P}^{1} $ without any collection of $ N $ points,
$ N>L+M+Md. $ Obviously, $ \Gamma\cap U\times{\mathbb P}^{1} $ cannot be contained in $ \left({\mathbb P}^{1}\times V_{0}\right)\cup\Delta $, $ \Delta $ being a
diagonal in $ {\mathbb P}^{1}\times{\mathbb P}^{1} $.

This shows that a family for which $ L>0 $ cannot be $ V_{0} $-conformal on any
open subset of $ S $, which proves the theorem. \end{proof}

\begin{definition} Say that the family $ {\mathbit C}_{s} $ is {\em clean\/} if Theorem~\ref{th13.52} remains
true with one mapping $ \alpha\colon S \to Z^{M}/{\mathfrak S}_{M} $ instead of a sequence of mappings $ \alpha_{l}:
S \to Z^{M_{l}}/{\mathfrak S}_{M_{l}} $. \end{definition}

\begin{remark} \label{rem13.41}\myLabel{rem13.41}\relax  One should expect that an appropriate modification of
definition of a conformal field will make all codominant conformal fields
to be clean. With the definitions used here, any infinite sum
$ {\mathbit C}^{\left(1\right)}+{\mathbit C}^{\left(2\right)}+\dots $ of conformal families on $ S $ acting on $ {\mathfrak M}{\mathcal M}\left(Z\right) $ is a conformal
family again, as far as all the summands but a finite number vanish on $ X_{{\mathcal L}} $
for any given line bundle $ {\mathcal L} $ on $ M $. It is the acceptability of such
infinite sums which makes it easy to construct examples of non-clean
codominant conformal families. \end{remark}

\begin{example} Consider the family $ E_{z} $, $ z\in Z $. Then one can consider
$ \alpha\left(z\right)=z $, thus this is a clean family.

A composition $ {\mathbit C}_{s}{\mathbit C}'_{t} $ of two clean families is clean again. A
subfamily of a clean family is clean, similarly, all the residues of a
clean family are clean. \end{example}

\begin{remark} Note that, inversely, any conformal family may be written as a
sum of clean families satisfying conditions of Remark~\ref{rem13.41}. \end{remark}

\begin{remark} Say that a conformal operator $ {\mathfrak M}{\mathcal M}\left(U\right) \xrightarrow[]{{\mathbit c}} {\mathfrak M}{\mathcal M}\left(Z\right) $ is {\em geometric\/} if it
is a homomorphism of algebras, say that it is {\em semigeometric\/} if there is a
homomorphism of algebras $ {\mathfrak M}{\mathcal M}\left(U\right) \xrightarrow[]{{\mathbit c}'} {\mathfrak M}{\mathcal M}\left(Z\right) $ such that $ {\mathbit c} $ is a corresponding
homomorphism of modules, i.e., $ {\mathbit c}\left(ab\right)={\mathbit c}'\left(a\right){\mathbit c}\left(b\right) $. Since $ {\mathfrak M}{\mathcal M}\left(Z\right) $ is considered
as a field of functions on $ {\mathcal M}\left(Z\right) $, one should consider a geometric operator
as a rational mapping $ {\mathcal M}\left(Z\right) \xrightarrow[]{f} {\mathcal M}\left(Z\right) $, and a semigeometric operator as a
rational fiber-wise mapping of line bundles over $ {\mathcal M}\left(Z\right) $.

To describe a geometric mapping $ {\mathbit c} $ it is enough to describe the
family $ {\mathbit c}E_{z} $, $ z\in Z $, of elements of $ {\mathfrak M}{\mathcal M}\left(Z\right) $. In particular, to any such mapping
one may associate a mapping $ Z \to Z^{m}/{\mathfrak S}_{m} $ for an appropriate $ m>0 $, or, what
is the same, a mapping $ \beta $ from a finite extension $ \widetilde{Z} $ of $ Z $ to $ Z $. As before,
one can show that a geometric mapping may be included into a codominant
family only if $ \beta $ is constant on $ \widetilde{Z}\smallsetminus\bar{Z} $, and either $ \bar{Z}\subset\widetilde{Z} $ is empty, or is an
irreducible component of $ \widetilde{Z} $ which projects isomorphically on $ Z $, and $ \beta|_{\bar{Z}} $
coincides with this projection.Note that in this case the element $ \alpha $ of
$ Z^{M}/{\mathfrak S}_{M} $ from Theorem~\ref{th13.52} may be taken to be $ \beta\left(\widetilde{Z}\smallsetminus\bar{Z}\right) $.

Note that a semigeometric mapping $ {\mathbit c} $ uniquely determines the
corresponding geometric mapping $ {\mathbit c}' $, since $ {\mathbit c}'\left(a\right)={\mathbit c}\left(a\right)/{\mathbit c}\left(1\right) $. One can write
$ {\mathbit c}\left(1\right) = \Phi\left(E_{z_{1}},\dots ,E_{z_{m}}\right) $ with an appropriate rational function $ \Phi $ (possibly
including some derivatives of $ E_{z_{k}} $ as well), thus $ {\mathbit c}=\Phi\left({\mathbit E}_{z_{1}},\dots ,{\mathbit E}_{z_{m}}\right){\mathbit c}' $. The
mapping $ {\mathbit c} $ may be included into a codominant family iff $ {\mathbit c}' $ may. The
corresponding element $ \alpha $ of $ Z^{M}/{\mathfrak S}_{M} $ is the union of the corresponding
element for $ {\mathbit c}' $ and of $ \left\{z_{1},\dots ,z_{m}\right\} $. One may see that smooth families of
semigeometric mappings are clean, and Theorem~\ref{th13.52} gives estimates
for the positions of $ z_{i} $.

Note that all the examples of conformal fields considered in this
paper are semigeometric. \end{remark}

\subsection{Conformal algebras }\label{s13.50}\myLabel{s13.50}\relax  Consider conformal fields $ {\mathbit C}\left(s\right) $ and $ {\mathbit C}'\left(t\right) $,
$ s\in S $, $ t\in T $, acting on $ {\mathfrak M}{\mathcal M}\left(Z\right) $. Suppose that $ {\mathbit C} $ is defined on $ U\subset S\times Z $, $ {\mathbit C}' $ is a
clean codominant family with the associated mapping $ \alpha\colon T \to Z^{m}/{\mathfrak S}_{m} $. Define
the {\em exceptional subset\/} $ \sigma_{{\mathbit C}{\mathbit C}'}\subset S\times T $ by
\begin{equation}
\left(s,t\right)\in\sigma_{{\mathbit C}{\mathbit C}'}\text{ if }\alpha\left(t\right)\cap R_{s}\not=\varnothing,
\notag\end{equation}
here $ R_{s_{0}}\subset Z $ is the subset of $ Z $ where $ {\mathbit C}\left(s_{0}\right) $ is not defined, i.e.,
$ \left(S\times Z\smallsetminus U\right)\cap\left(\left\{s_{0}\right\}\times Z\right) $.

The following obvious theorem shows that a product of several clean
codominant conformal fields has singularity uniquely determined by
singularities of pairwise products. In particular, triple product has no
fusions which do not correspond to fusions of pairwise products. This
sets a foundation stone for the theory of conformal algebras, since it
allows one to investigate relations on fusions implied by associativity
of multiplication.

\begin{theorem} \label{th13.22}\myLabel{th13.22}\relax  Consider smooth clean codominant families $ {\mathbit C}_{1}\left(s_{1}\right) $,
$ {\mathbit C}_{2}\left(s_{2}\right) $, \dots , $ {\mathbit C}_{k}\left(s_{k}\right) $, $ s_{i}\in S_{i} $, $ i=1,\dots ,k $, acting on $ {\mathfrak M}{\mathcal M}\left(Z\right) $. Then the product
family
\begin{equation}
{\mathbit C}_{1}\left(s_{1}\right){\mathbit C}_{2}\left(s_{2}\right)\dots {\mathbit C}_{k}\left(s_{k}\right)
\notag\end{equation}
is smooth at $ \left(s_{1},s_{2},\dots ,s_{k}\right) $ if $ \left(s_{i},s_{j}\right)\notin\sigma_{{\mathbit C}_{i}{\mathbit C}_{j}} $ for any $ 1\leq i<j\leq k $. \end{theorem}

\begin{remark} Note that by complicating notations by having a sequence of
subschemes $ \sigma_{{\mathbit C}_{i}{\mathbit C}_{j}}^{\left(l\right)} $ instead of one subscheme $ \sigma_{{\mathbit C}_{i}{\mathbit C}_{j}} $ one can remove
the requirement that the families are clean. Similarly, one may remove
smoothness requirement by introducing subschemes $ \sigma_{{\mathbit C}_{i}}\subset S_{i} $ where $ {\mathbit C}_{i} $ is not
smooth, and stating that
\begin{equation}
{\mathbit C}_{1}\left(s_{1}\right){\mathbit C}_{2}\left(s_{2}\right)\dots {\mathbit C}_{k}\left(s_{k}\right)
\notag\end{equation}
is smooth at $ \left(s_{1},s_{2},\dots ,s_{k}\right) $ if $ s_{i}\notin\sigma_{{\mathbit C}_{i}} $ for any $ 1\leq i\leq k $, and $ \left(s_{i},s_{j}\right)\notin\sigma_{{\mathbit C}_{i}{\mathbit C}_{j}} $
for any $ 1\leq i<j\leq k $.

Thus the only important requirement in the theorem~\ref{th13.22} is the
codominantness. Recall that the codominantness ensures that the family
may be composed on the left with an arbitrary conformal operator. \end{remark}

\begin{remark} \label{rem13.57}\myLabel{rem13.57}\relax  Yet another reference to classification of meromorphic
families in $ {\mathfrak M}{\mathcal M}\left(Z\right) $ shows that for any $ k,m\geq0 $ there is a number $ L=L\left(k,m\right) $
such that for a ``big enough'' family $ F_{t} $, $ t\in T $, the family
\begin{equation}
\Phi\left(s_{1},s_{2},\dots ,s_{k}\right){\mathbit C}_{1}\left(s_{1}\right){\mathbit C}_{2}\left(s_{2}\right)\dots {\mathbit C}_{k}\left(s_{k}\right)F_{t}
\notag\end{equation}
is smooth in $ \left(s_{1},s_{2},\dots ,s_{k}\right) $ (in the sense that it is smooth on a subset
of $ S_{1}\times\dots \times S_{k}\times T $ which projects onto $ S_{1}\times\dots \times S_{k} $) as far as $ \Phi $ has a zero of
order $ L $ on the subsets given by $ s_{i}\in\sigma_{{\mathbit C}_{i}} $ for any $ 1\leq i\leq k $, and by
$ \left(s_{i},s_{j}\right)\in\sigma_{{\mathbit C}_{i}{\mathbit C}_{j}} $ for any $ 1\leq i<j\leq k $. \end{remark}

\begin{remark} Theorem~\ref{th13.22} shows that in the settings of this paper
the conformal associativity relations are corollaries of geometry of
punctual Hilbert schemes. Note the similarity with \cite{Kon97Def}, where the
associativity relations are obtained as corollaries of geometry of
moduli spaces. \end{remark}

The same arguments show

\begin{proposition} \label{prop13.23}\myLabel{prop13.23}\relax  Consider clean codominant families $ {\mathbit C}\left(s\right) $ and $ {\mathbit C}'\left(t\right) $,
$ s\in S $, $ t\in T $, let $ {\mathbit D} $ be a fusion of $ {\mathbit C} $ and $ {\mathbit C}' $ on an irreducible subvariety
$ P\subset S\times T $. If $ P $ projects onto $ S $ and onto $ T $, then
\begin{enumerate}
\item
$ {\mathbit D} $ is clean codominant;
\item
if $ {\mathbit D}'\left(r\right) $, $ r\in R $, is another clean codominant family, then
\begin{equation}
\sigma_{{\mathbit D}{\mathbit D}'}\subset\left(\left(\sigma_{{\mathbit C}{\mathbit D}'}\times T\right)\cup\left(S\times\sigma_{{\mathbit C}'{\mathbit D}'}\right)\right)\cap\left(P\times R\right)\subset S\times T\times R.
\notag\end{equation}

\end{enumerate}
\end{proposition}

\subsection{Fractions of formal series }\label{s13.60}\myLabel{s13.60}\relax  Consider a meromorphic function $ F $ on
$ Z $ which has no pole at $ z_{0}\in Z $. Then one can consider $ F $ as a regular
function on a formal completion $ \widehat{Z} $ of $ Z $ at $ z_{0} $, i.e., as a formal function
near $ z_{0} $. Choosing a coordinate system $ \left(\eta_{1},\dots ,\eta_{d}\right) $ near a smooth point $ z_{0} $
makes $ F $ into a formal series in $ \eta_{1},\dots ,\eta_{d} $.

If one lifts the restriction that $ F $ has no pole at $ z_{0} $, then some
extension of the algebra $ {\mathbb K}\left[\left[\eta_{1},\dots ,\eta_{d}\right]\right] $ of formal series in $ \eta_{1},\dots ,\eta_{d} $
should be introduced. A natural choice is the field of fractions of
$ {\mathbb K}\left[\left[\eta_{1},\dots ,\eta_{d}\right]\right] $, but it is a pain to work with if $ d>1 $, since
$ {\mathbb K}\left[\left[\eta_{1},\dots ,\eta_{d}\right]\right] $ has too many minimal prime ideals. One way to construct a
field which contains formal functions {\em and\/} meromorphic functions is to
consider an iterated construction of residue from Section~\ref{s11.8}.

Consider a {\em flag\/} $ Z=Z_{d}\supset\dots \supset Z_{0}=\left\{z_{0}\right\} $ of irreducible subvarieties of $ Z $
(or $ \infty $-jets of subvarieties of $ Z $ near $ z_{0} $, i.e., subvarieties of $ \widehat{Z} $). A
{\em normalized flag\/} is a flag with a lifting $ \iota_{k}\colon \widetilde{Z}_{k} \hookrightarrow \widetilde{Z}_{k+1} $ of the inclusion
$ Z_{k} \hookrightarrow Z_{k+1} $ into the normalizations $ \widetilde{Z}_{k} $, $ \widetilde{Z}_{k+1} $ of $ Z_{k} $ and $ Z_{k+1} $, for every
$ 0\leq k<d $. Fix an {\em equation\/} $ \eta_{k+1} $ of $ Z_{k} $ in $ Z_{k+1} $, i.e., a meromorphic function
on $ Z_{k+1} $ which has a zero of the first order at $ \iota_{k}\widetilde{Z}_{k}\subset\widetilde{Z}_{k+1} $. We will also
fix meromorphic vector fields $ v_{k+1} $ on $ \widetilde{Z}_{k+1} $ such that $ \left< v_{k+1},d\eta_{k+1} \right>=0 $,
and the restriction of $ v_{k+1} $ to $ \iota_{k}\widetilde{Z}_{k} $ is defined and not tangent to $ \iota_{k}\widetilde{Z}_{k} $.
Such data allows one taking residues of meromorphic functions on $ Z_{k+1} $
along $ Z_{k} $ (compare Section~\ref{s11.8}). Note that if all $ Z_{k} $ are normal,
$ \left(\eta_{1},\dots ,\eta_{d}\right) $ may be extended to become a {\em regular sequence\/} of functions on
an open subset of $ Z $ (see \cite{Har77Alg}), in smooth case one may assume
instead that $ \eta_{\bullet} $ gives a coordinate system near $ z_{0} $.

Any meromorphic function $ F $ on $ Z $ (or an element of the field of
fractions of $ {\mathbb K}\left[\left[\eta_{1},\dots ,\eta_{d}\right]\right] $) is uniquely determined by its residues on
$ Z_{d-1} $ of all orders, which are in turn meromorphic functions on $ Z_{d-1} $ (or
elements of the corresponding field of fractions). Now one can take
residues of these meromorphic functions on $ Z_{d-2} $, and so on.

As a result, one associates to $ F $ a {\em formal Laurent series\/}\footnote{Note that though $ \eta_{i} $ are defined on subvarieties $ Z_{i} $ only, one can extend
them to (at least) the formal neighborhood of $ Z_{i} $ in $ Z $ using vector fields
$ v_{i} $.}
\begin{equation}
\sum_{k_{1},\dots ,k_{d}}\phi_{k_{1},\dots ,k_{d}}\eta_{1}^{k_{1}}\dots \eta_{d}^{k_{d}},
\label{equ13.10}\end{equation}\myLabel{equ13.10,}\relax 
and the coefficients $ \phi_{k_{1},\dots ,k_{d}} $ vanish if
\begin{equation}
k_{1}<-K_{1}\text{, }k_{2}<-K_{2}\left(k_{1}\right)\text{, \dots , }k_{d}<-K_{d}\left(k_{1},\dots ,k_{d-1}\right).
\notag\end{equation}
Note that the series as the above one may be added, multiplied, {\em and\/}
divided. Denote the field of such formal series by
$ {\mathbb K}\left(\left(\eta_{d}\right)\right)\left(\left(\eta_{d-1}\right)\right)\dots \left(\left(\eta_{1}\right)\right) $. We consider elements of this field as functions
on the {\em formal completion\/} $ \widehat{Z}' $ {\em of\/} $ Y $ {\em at the normalized flag\/} $ Z=Z_{d}\supset\dots \supset Z_{0}=\left\{z_{0}\right\} $.

If a different sequence $ \left(\eta_{1}',\dots ,\eta_{d}'\right) $ results in the same flag, then
one can easily establish a canonical isomorphism of
$ {\mathbb K}\left(\left(\eta_{d}\right)\right)\left(\left(\eta_{d-1}\right)\right)\dots \left(\left(\eta_{1}\right)\right) $ and $ {\mathbb K}\left(\left(\eta'_{d}\right)\right)\left(\left(\eta'_{d-1}\right)\right)\dots \left(\left(\eta'_{1}\right)\right) $. However,
completions at different normalized flags are not related, as the
following example shows:

Consider a transposition $ \sigma\in{\mathfrak S}_{d} $, suppose that $ \eta_{i} $ and $ v_{i} $ are extended
to $ Z $, and $ \left(\eta_{\sigma_{1}},\dots ,\eta_{\sigma_{d}}\right) $ is a regular sequence too. Then it defines a
different normalized flag, and a different completion $ \widehat{Z}'' $ of $ Z $. A
question arises whether one is able to compare functions on $ \widehat{Z}' $ and $ \widehat{Z}'' $
for equality.

Suppose that a function $ F $ in $ {\mathbb K}\left(\left(\eta_{d}\right)\right)\left(\left(\eta_{d-1}\right)\right)\dots \left(\left(\eta_{1}\right)\right) $ becomes an
element $ F_{1} $ in $ {\mathbb K}\left[\left[\eta_{1},\dots ,\eta_{d}\right]\right] $ after multiplication by an element
$ F_{2}\in{\mathbb K}\left[\left[\eta_{1},\dots ,\eta_{d}\right]\right] $, i.e., $ F=F_{1}/F_{2} $. Obviously, $ F $ can be represented as an
element of $ {\mathbb K}\left(\left(\eta_{\sigma_{d}}\right)\right)\left(\left(\eta_{\sigma_{d-1}}\right)\right)\dots \left(\left(\eta_{\sigma_{1}}\right)\right) $ as well.

\begin{example} Consider $ F=\sum_{k\geq0}\eta_{1}^{k}\eta_{2}^{-k}\in{\mathbb K}\left(\left(\eta_{2}\right)\right)\left(\left(\eta_{1}\right)\right) $. Then $ F=\frac{\eta_{2}}{\eta_{2}-\eta_{1}} $, thus
can be written as $ F=-\sum_{k\geq1}\eta_{1}^{-k}\eta_{2}^{k}\in{\mathbb K}\left(\left(\eta_{1}\right)\right)\left(\left(\eta_{2}\right)\right) $. \end{example}

Similarly, one can define a relation of {\em being an analytic
continuation\/} of each other for completions of $ Z $ at arbitrary pair of
normalized flags, but as with usual analytic continuation, it can be made
$ 1 $-to-1 on very small subsets only.

\begin{remark} Consider how the Laurent coefficient at a normalized flag are
related to the ``usual'' Laurent coefficients of a function of two complex
variables. Given an analytic function $ \xi\left(z_{1},z_{2}\right) $ defined on an open subset
$ U $ of $ {\mathbb C}^{2} $ one may define Laurent coefficients as integrals
\begin{equation}
\xi_{kl}=\frac{1}{\left(2\pi i\right)^{2}}\iint\Sb|z_{1}|=\varepsilon_{1} \\ |z_{2}|=\varepsilon_{2}\endSb\frac{\xi\left(z_{1},z_{2}\right)}{z_{1}^{k_{1}+1}z_{2}^{k_{2}+1}}dz_{1}dz_{2},
\notag\end{equation}
here $ \varepsilon_{1} $, $ \varepsilon_{2} $ are appropriate numbers such that the torus $ {\mathbb T}_{\varepsilon_{1}\varepsilon_{2}}^{2}
=\left\{|z_{1}|=\varepsilon_{1},|z_{2}|=\varepsilon_{2}\right\} $ is in $ U $. The coefficients do not change if one deforms
$ \varepsilon_{1} $, $ \varepsilon_{2} $ in such a way that $ {\mathbb T}_{\varepsilon_{1}\varepsilon_{2}}^{2} $ remains in $ U $ during the deformation.

Consider the flag defined by the sequence $ \left(z_{1},z_{2}\right) $. Laurent
coefficients near this flag correspond to considering the residues at
$ z_{1}=0 $, and residues of these residues at $ z_{2}=0 $. In other words, one takes
the limit $ \varepsilon_{1} \to $ 0, and then the limit $ \varepsilon_{2} \to $ 0. On the other hand, if $ F $ is
meromorphic, then the possibility of description of roots of denominator
by Puiseux series shows that there exist $ \varepsilon>0 $ and $ K\geq0 $ such that $ F\left(z_{1},z_{2}\right) $
is defined on
\begin{equation}
0<|z_{1}|<|z_{2}|^{K}\text{, }0<|z_{2}|<\varepsilon.
\notag\end{equation}
Obviously, the above double limit coincides with the ``usual'' Laurent
coefficients for any $ 0<\varepsilon_{1}<\varepsilon_{2}^{K} $, $ 0<\varepsilon_{2}<\varepsilon $.

Similarly, the flag given by the sequence $ \left(z_{2},z_{1}\right) $ corresponds to
taking $ 0<\varepsilon_{2}<\varepsilon_{1}^{K} $, $ 0<\varepsilon_{1}<\varepsilon $. When one deforms the contour of integration from
one region to another, it may hit the divisor of poles of $ F $, each time
adding residues at the corresponding irreducible component of this
divisor. \end{remark}

\subsection{Laurent coefficients of conformal fields } Similarly to what is done in
the previous section, given a meromorphic function $ F $ on $ S\times X $, and a
normalized flag $ \left\{S_{i}\right\} $ on $ S $ equipped with equations $ \eta_{i} $ and vector fields
$ v_{i} $, one can consider $ F $ as an element of $ {\mathcal M}\left(X\right)\left(\left(\eta_{d}\right)\right)\left(\left(\eta_{d-1}\right)\right)\dots \left(\left(\eta_{1}\right)\right) $.
Moreover, if $ F_{s} $ is a meromorphic function on $ {\mathcal M}\left(Z\right) $ parameterized by $ s\in S $,
one can consider $ F $ as an element of $ {\mathfrak M}{\mathcal M}\left(Z\right)\left(\left(\eta_{d}\right)\right)\left(\left(\eta_{d-1}\right)\right)\dots \left(\left(\eta_{1}\right)\right) $.

Apply the defined above completion to a conformal field $ {\mathbit C} $ on $ S $
acting in $ {\mathfrak M}{\mathcal M}\left(Z\right) $. Suppose that $ {\mathbit C} $ dominates $ U\subset Z $. Given a function $ F\in{\mathfrak M}{\mathcal M}\left(U\right) $,
$ {\mathbit C}\left(s\right)F $ is a meromorphic function on $ S\times{\mathcal M}\left(Z\right) $, thus one can consider it as a
function on $ \widehat{S}\times{\mathcal M}\left(Z\right) $, $ \widehat{S} $ being the formal completion of $ S $ at the normalized
flag $ \left\{S_{i}\right\} $. Consider Laurent coefficients $ \phi_{k_{1}\dots k_{D}}\in{\mathfrak M}{\mathcal M}\left(Z\right) $ of $ {\mathbit C}\left(s\right)F $, $ D=\dim 
S $. Each of these coefficient is obtained by a sequence of taking a fine
residue (defined in Section~\ref{s11.85}) at subvarieties of the normalized
flag $ \left\{S_{i}\right\} $ (next subvariety is a hypersurface in the previous one), thus
there are operators $ {\mathbit c}_{k_{1}\dots k_{D}} $ such that
\begin{equation}
\phi_{k_{1}\dots k_{D}} = {\mathbit c}_{k_{1}\dots k_{D}}F
\notag\end{equation}
and $ {\mathbit c}_{k_{1}\dots k_{D}} $ gives a linear mapping $ {\mathfrak M}{\mathcal M}\left(U\right) \to {\mathfrak M}{\mathcal M}\left(Z\right) $ (which is not
necessarily a conformal operator, since it may fail to send smooth
families to smooth families). Call these operators {\em Laurent coefficients\/}
of $ {\mathbit C} $ at the normalized flag $ \left\{S_{i}\right\} $. (Strictly speaking, one should
associate Laurent coefficients to the sequence of equations $ \eta_{i} $ and vector
fields $ v_{i} $, but as it was noted above, different systems of equations lead
to isomorphic fields of functions).

For any family $ F_{t} $, $ t\in T $, of functions $ F $ the functions $ {\mathbit c}_{k_{1}\dots k_{D}}F_{t} $
vanish for sufficiently big negative $ \left(k_{1},\dots ,k_{D}\right) $ (big in the same sense
as in~\eqref{equ13.10}). As we have seen in Remark~\ref{rem11.315}, the {\em operators\/}
$ {\mathbit c}_{k_{1}\dots k_{D}} $ do not necessarily vanish for big negative $ \left(k_{1},\dots ,k_{D}\right) $.

\begin{proposition} \label{prop13.69}\myLabel{prop13.69}\relax  Consider a clean codominant conformal field $ {\mathbit C}\left(s\right) $ on
$ S $ acting on $ {\mathfrak M}{\mathcal M}\left(Z\right) $ and a flag $ S_{\bullet} $ on $ S $. Suppose that $ {\mathbit C} $ dominates $ U\subset Z $. Then
there is a finite subset $ Z_{0}\subset Z $ such that for any multiindex $ \left(k_{1},\dots ,k_{D}\right) $
operator $ {\mathbit c}_{k_{1}\dots k_{D}} $ sends $ {\mathfrak M}{\mathcal M}\left(U_{1}\right) \to {\mathfrak M}{\mathcal M}\left(U_{2}\right) $ if $ U_{1}\subset U $, $ U_{2}\supset U_{1}\cup Z_{0} $. \end{proposition}

\begin{proof} This is an immediate corollary of Theorem~\ref{th13.52} and the
following obvious

\begin{lemma} Consider a family $ F_{s} $ in $ {\mathfrak M}{\mathcal M}\left(Z\right) $ parameterized by $ S $. Suppose that
this family corresponds to a mapping $ \beta\colon S \to {\mathfrak Q}_{k,m} $. Let $ \pi $ be the
projection $ {\mathfrak Q}_{k,m} \to Z^{\left[m\right]} $. Then for any flag $ S_{0}\subset\dots \subset S_{d} $ in $ S $ the Laurent
coefficients $ F_{k_{1}\dots k_{d}} $ of $ F $ at $ S_{\bullet} $ are in $ {\mathfrak M}{\mathcal M}\left(U\right) $ if $ U\supset\pi\left(\beta\left(S_{0}\right)\right) $. \end{lemma}\end{proof}

\subsection{Laurent coefficients of a composition }\label{s13.81}\myLabel{s13.81}\relax  Consider two conformal
fields, $ {\mathbit C}\left(s\right) $ on $ S $, and $ {\mathbit C}'\left(t\right) $ on $ T $, both acting on $ {\mathfrak M}{\mathcal M}\left(Z\right) $. Suppose that $ {\mathbit C}' $
is codominant, thus the composition $ {\mathbit C}\left(s\right){\mathbit C}'\left(t\right) $ is defined. Let $ {\mathcal S}=\left\{S_{i}\right\} $,
$ i=0,\dots ,D $, be a normalized flag on $ S $, $ {\mathcal T}=\left\{T_{j}\right\} $, $ j=0,\dots ,D' $, be a normalized
flag on $ T $, define a normalized flag $ {\mathcal S}*{\mathcal T} $ on $ S\times T $ to be
\begin{equation}
S_{D}\times T_{D'}\supset\dots \supset S_{D}\times T_{1}\supset S_{D}\times T_{0}\supset\dots \supset S_{0}\times T_{0}.
\notag\end{equation}
Note that similarly defined $ {\mathcal T}*{\mathcal S} $ can be considered as a flag on $ S\times T\simeq T\times S $
too, and that additional data $ \eta_{i} $, $ v_{i} $ for $ {\mathcal S} $ and $ {\mathcal T} $ gives similar data for
$ {\mathcal S}*{\mathcal T} $ and $ {\mathcal T}*{\mathcal S} $.

\begin{theorem} \label{th13.80}\myLabel{th13.80}\relax  Consider Laurent coefficients $ {\mathbit c}_{\bullet} $, $ {\mathbit c}'_{\bullet} $ of $ {\mathbit C} $ and $ {\mathbit C}' $ at $ {\mathcal S} $
and $ {\mathcal T} $. Suppose that $ {\mathbit C}' $ is codominant, let $ Z_{0} $ is the corresponding finite
subset of $ Z $, as in Proposition~\ref{prop13.69}. Suppose that $ {\mathbit C} $ dominates $ U\supset Z_{0} $,
$ {\mathbit C}' $ dominates $ U' $, so the composition $ {\mathbit c}_{K}\circ{\mathbit c}'_{K'} $ is well-defined on $ {\mathfrak M}{\mathcal M}\left(U\cap U'\right) $.
Then $ {\mathbit c}_{K}\circ{\mathbit c}'_{K'} $ restricted to $ {\mathfrak M}{\mathcal M}\left(U\cap U'\right) $ coincides with the $ \left(K,K'\right) $-th Laurent
coefficient of $ {\mathbit C}\circ{\mathbit C}' $ at $ {\mathcal S}*{\mathcal T} $. \end{theorem}

\begin{proof} Fix a function $ F\in{\mathfrak M}{\mathcal M}\left(U\cap U'\right) $, then $ {\mathbit C}'\left(t\right)F $ and $ {\mathbit C}\left(s\right)\left({\mathbit C}'\left(t\right)F\right) $ are
well-defined meromorphic families. Apply the classification of
meromorphic families of Theorem~\ref{th13.32} to $ {\mathbit C}'\left(t\right)F $. This leads to a
blow-up $ T' $ of $ T $.

\begin{lemma} Consider a normalized flag $ T_{\bullet} $ in $ T $, and a blow-up $ T' \to T $. Then
$ T_{\bullet} $ may be naturally lifted to a normalized flag $ T'_{\bullet} $ in $ T' $, each $ T'_{k} $ being
a blow-up of $ T_{k} $. Given equations $ \eta_{\bullet} $ of $ T_{\bullet} $ and vector fields $ v_{\bullet} $ (as in
Section~\ref{s13.60}), there is a natural lifting of these data to $ T' $. \end{lemma}

\begin{proof} Let $ d=\dim  T $. It is enough to consider $ T_{d-1} $, and then proceed by
induction. The normalization $ \widetilde{T}' $ of $ T' $ is a blow-up of the normalization $ \widetilde{T} $
of $ T $ in a closed subscheme $ Z_{0}\subset\widetilde{T} $ of codimension greater than 2. In
particular, $ Z_{0} $ does not contain $ \iota_{d-1}\widetilde{T}_{d-1} $, thus there is a canonically
defined blowup $ \widetilde{T}'_{d-1} $ of $ \widetilde{T}_{d-1} $ and a mapping $ \widetilde{T}'_{d-1} \to \widetilde{T}'=\widetilde{T}'_{d} $. Obviously,
$ \widetilde{T}'_{d-1} $ is a normalization of its image in $ T' $.

The vector field $ v_{d} $ and the meromorphic function $ \eta_{d} $ on $ \widetilde{T}=\widetilde{T}_{d} $ may be
lifted to $ \widetilde{T}' $ in a canonical way. What we need to prove is that the
lifting of $ \eta_{d} $ has a zero at $ \widetilde{T}'_{d-1} $, and the lifting of $ v_{d} $ is defined at
$ \widetilde{T}'_{d-1} $, and is not tangent to $ \widetilde{T}'_{d-1} $. Both parts follow from $ \widetilde{T}'_{d-1} $ being
isomorphic to $ \widetilde{T}_{d-1} $ at generic point. \end{proof}

Call such a lifting a {\em strict preimage\/} of the flag $ T_{\bullet} $. Note that one
can also define a strict preimage for {\em incomplete\/} normalized flags
$ T_{k}\subset T_{k+1}\subset\dots \subset T_{d}=T $ as far as each subvariety is a hypersurface in the next
subvariety.

Obviously, for a given meromorphic family $ F_{t} $ in $ {\mathfrak M}{\mathcal M}\left(Z\right) $, $ t\in T $, and a
normalized flag $ T_{\bullet} $ on $ T $ the Laurent coefficients of $ F_{t} $ at $ T_{\bullet} $ do not
change after any blowup of $ T $ (with the corresponding lifting of $ \eta_{i} $ and
$ v_{i} $). Moreover, multiplication of $ {\mathbit C}\left(s\right) $ and $ {\mathbit C}'\left(t\right) $ by scalar meromorphic
functions of $ s $ and $ t $ does not change the validity of the theorem. In
particular, one may assume that both $ {\mathbit C} $ and $ {\mathbit C}' $ are smooth, and in
conditions of Theorem~\ref{th13.32} the family $ F_{t}={\mathbit C}'\left(t\right)F $ leads to no blow-up
of $ T $.

This eliminates the field $ {\mathbit C}'\left(t\right) $, instead one may consider Laurent
coefficients of $ {\mathbit C}\left(s\right)F_{t} $ for an arbitrary smooth family $ F_{t} $ in $ {\mathfrak M}{\mathcal M}\left(U\right) $.

Note that for an incomplete normalized flag $ R_{k}\subset R_{k+1}\subset\dots \subset R_{d}=R $ with
each subvariety being a hypersurface in the next subvariety (with
associated data $ \eta_{i} $ and $ v_{i} $), and for a meromorphic function on $ R $ one can
define {\em partial Laurent coefficients\/} $ F_{l_{k}\dots l_{d}} $ of $ F $ which are meromorphic
functions on $ R_{k} $.

\begin{proposition} \label{prop13.330}\myLabel{prop13.330}\relax  Suppose that a conformal field $ {\mathbit C}\left(s\right) $ on $ S $ acting in
$ {\mathfrak M}{\mathcal M}\left(Z\right) $ dominates $ U\subset S $, and $ F_{t} $, $ t\in T $, is a smooth family in $ {\mathfrak M}{\mathcal M}\left(U\right) $. Fix a
normalized flag $ {\mathcal T} $ in $ T $. Then $ \left(l_{1},\dots ,l_{D'}\right) $-partial Laurent coefficient of
$ {\mathbit C}\left(s\right)F_{t} $ on the partial flag $ S\times{\mathcal T}\buildrel{\text{def}}\over{=}\left\{S\times T_{k-D}\right\}_{k=D,\dots ,D+D'} $ coincides with $ {\mathbit C}\left(s\right)
F_{l_{1}\dots l_{D'}} $, $ F_{l_{1}\dots l_{D'}} $ being $ \left(l_{1},\dots ,l_{D'}\right) $-Laurent coefficient of $ F $ at the
flag $ {\mathcal T} $. \end{proposition}

\begin{proof} Again, $ {\mathbit C}\left(s\right)F_{t} $ is a family on $ R=S\times T $, thus Theorem~\ref{th13.32} leads
to a blow-up $ R' $ of $ R $. There is a mapping $ \beta\colon R' \to {\mathfrak Q}_{k,m} $ which induces the
family $ {\mathbit C}\left(s\right)F_{t} $.

Suppose that $ R' $ is the blow-up of $ R $ in $ \rho\subset R=S\times T $. For $ s $ in an open
subset $ V\subset S $ the intersection $ \rho\cap\left(\left\{s\right\}\times T\right) $ depends smoothly on $ s $ (i.e.,
$ \rho\cap\left(V\times T\right) $ is flat over $ V $). Decreasing $ V $, one may assume that $ \rho\cap\left(V\times T_{k}\right) $ is
also flat over $ V $ for any $ k\leq\dim  T $. This shows that for $ s\in V $ the preimage
$ T_{\bullet}^{\left(s\right)} $ in $ R' $ of the flag $ \left\{s\right\}\times T_{\bullet} $ in $ \left\{s\right\}\times T\subset R $ is well-defined, and depends
smoothly on $ s $. In particular, for any $ s\in V $ one can canonically define a
lifting of $ T_{0} $ to $ R' $, i.e., there is a mapping $ \tau_{0}\colon V \to R' $. Let $ V' $ be the
image of this mapping.

On the other hand, consider an incomplete flag $ {\mathcal R}=\left\{{\mathcal R}_{k}\right\}=S\times{\mathcal T} $ on $ S\times T $. It
has a strict preimage $ {\mathcal R}' $ in $ R' $. Obviously, on the preimage of $ V\times T $ the
subvariety $ {\mathcal R}'_{D} $ coincides with the image of $ \tau_{0} $.

Denote the smallest variety of the flag $ {\mathcal T} $ by $ \left\{t_{0}\right\} $. Since $ {\mathbit C}\left(s\right) $
dominates $ U $, $ {\mathbit C}\left(s\right)F_{t_{0}} $ is smooth on an open subset of $ S $. Decrease $ V $ so that
$ {\mathbit C}\left(s\right)F_{t_{0}} $ is smooth if $ s\in V $. Consider $ \beta_{1}=\beta\circ\tau_{0} $, it sends $ V $ into $ \mathring{{\mathfrak Q}}_{k,m} $. Thus
for any $ s_{0}\in V $ there is an open subset $ W^{\left(s_{0}\right)}\subset T^{\left(s_{0}\right)} $ such that $ W^{\left(s_{0}\right)}\ni\tau_{0}\left(s_{0}\right) $
and $ \beta $ sends $ W^{\left(s_{0}\right)} $ into $ \mathring{{\mathfrak Q}}_{k,m} $, hence the family $ {\mathbit C}\left(s_{0}\right)F_{t} $ becomes smooth when
restricted/lifted to an open subset $ W^{\left(s_{0}\right)} $ of a blowup of $ T $.

\begin{lemma} Consider a conformal operator $ {\mathbit C}\colon {\mathfrak M}{\mathcal M}\left(U\right) \to {\mathfrak M}{\mathcal M}\left(Z\right) $ and a smooth
family $ F_{t} $ of elements of $ {\mathfrak M}{\mathcal M}\left(U\right) $ parameterized by $ T $. Fix a flag $ T_{\bullet} $ in $ T $.
Then the Laurent ($ = $Taylor) coefficients of $ {\mathbit C}F_{t} $ in $ t $ at $ T_{\bullet} $ are images by
$ {\mathbit C} $ of Laurent coefficients of $ F_{t} $ in $ t $. \end{lemma}

\begin{proof} This obvious for $ \left(0,\dots ,0\right) $-th Laurent coefficient, since $ {\mathbit C} $
sends smooth families to smooth families. Note that $ {\mathbb K} $-linearity of $ {\mathbit C} $ and
the property of Laurent coefficients we know already imply that
\begin{enumerate}
\item
taking a blowup of $ T $ does not change validity of the lemma;
\item
if $ \varphi $ is a regular scalar function on $ T $, then $ {\mathbit C}\varphi\left(t\right)F_{t}=\varphi\left(t\right){\mathbit C}F_{t} $;
\item
thus if $ F_{t} $ has a zero of $ k $-th order at a hypersurface $ T_{d-1}\subset T $, $ {\mathbit C}F_{t} $
has a zero of $ k $-th order at the same hypersurface;
\item
If $ F_{t} $ does not depend on $ t $, $ {\mathbit C}F_{t} $ does not depend on $ t $.
\item
If $ F_{t} $ depends on $ t $ as $ \varphi\left(t\right)F^{\left(0\right)} $, $ {\mathbit C}F_{t} $ is $ \varphi\left(t\right){\mathbit C}F^{\left(0\right)} $.
\end{enumerate}

Consider the case $ \dim  T=1 $ first. The lemma follows from the above
properties alone. Indeed, the first property shows that it is enough to
consider normal curves. The third property insures that the $ k $-th Laurent
coefficient of $ {\mathbit C}F_{t} $ depends only on $ l $-th coefficients of $ {\mathbit C}_{t} $ with $ l\leq k $, and
the last one shows that $ l $-th coefficient of $ {\mathbit C}_{t} $ influences $ l $-th
coefficient of $ {\mathbit C}F_{t} $ only. (Note that in case $ \dim  T>1 $ there may be no
projection $ T \to T_{d-1} $, which breaks the last argument.)

Now proceed by induction in $ d=\dim  T $. One may assume that the
statement is true when the family is restricted to $ T_{d-1} $. It is enough to
show that residues of $ F_{t} $ and $ {\mathbit C}F_{t} $ on $ T_{d-1} $ are related by $ {\mathbit C} $. On the other
hand, since the lemma is already proven if $ d=1 $, one may apply it to
curves in $ T_{d} $ which are transversal to $ T_{d-1} $, which proves the statement
about the residues. \end{proof}

This finishes the proof of the proposition, since a blow-up does not
change Laurent coefficients. \end{proof}

By definition, Laurent coefficients at $ {\mathcal S}*{\mathcal T} $ are Laurent coefficients
at $ {\mathcal S}\times T_{0} $ of partial Laurent coefficients at $ S\times{\mathcal T} $. Since partial Laurent
coefficients of $ {\mathbit C}\left(s\right)F_{t} $ at $ S\times{\mathcal T} $ can be written as $ {\mathbit C}\left(s\right)F_{l_{1}\dots l_{D'}} $, and
Laurent coefficients of $ {\mathbit C}\left(s\right)G $ are (by definition) Laurent coefficients of
$ {\mathbit C} $ applied to $ F_{l_{1}\dots l_{D'}} $, this finishes the proof of the theorem. \end{proof}

\begin{remark} Note that since a blow-up of a manifold is not necessarily a
manifold, use of non-smooth flags is necessary (as far as one follows the
same arguments). I do not know whether one can omit consideration of
non-normal flags. \end{remark}

\begin{remark} Suppose that conformal fields $ {\mathbit C}\left(s\right) $ and $ {\mathbit C}'\left(t\right) $ commute:
\begin{equation}
{\mathbit C}\left(s\right){\mathbit C}'\left(t\right) = {\mathbit C}'\left(t\right){\mathbit C}\left(s\right).
\notag\end{equation}
However, a fundamental asymmetry between $ {\mathcal S}*{\mathcal T} $ and $ {\mathcal T}*{\mathcal S} $ in Theorem~\ref{th13.80}
shows that Laurent coefficients of $ {\mathbit C} $ and of $ {\mathbit C}' $ {\em do not necessarily
commute\/}! As we have seen it in Section~\ref{s13.60}, the non-commutativity
corresponds to fusions of $ {\mathbit C} $ and $ {\mathbit C}' $ on subvarieties passing through $ s_{0}\times t_{0} $,
here $ \left\{s_{0}\right\} $ and $ \left\{t_{0}\right\} $ are smallest elements of flags on $ S $ and $ T $
correspondingly.

This observation is the core of the {\em vertex commutation relations},
considered in the following section. \end{remark}

\begin{remark} Note that Example~\ref{ex14.30} will provide some pseudo-explanation
why it is better to pick up the flag $ {\mathcal S}*{\mathcal T} $ than $ {\mathcal T}*{\mathcal S} $. Indeed, if Theorem~%
\ref{th13.80} would contain $ {\mathcal T}*{\mathcal S} $, then~\eqref{equ55.11} would have another sign, thus~%
\eqref{equ6.12} would give a wrong sign for commutators of Laurent coefficients
of $ \psi $ and $ \psi^{+} $. \end{remark}


\section{Chiral algebras and commutation relations }\label{h14}\myLabel{h14}\relax 

Recall that in Section~\ref{s13.81} it was shown that Laurent
coefficients of commuting conformal fields do not necessarily commute,
and the commutator may be calculated basing on fusions of these conformal
fields. Moreover, in Section~\ref{s13.50} it was shown that products of
several codominant conformal fields have fusions at subvarieties which
are determined by pairwise fusions only.

In this section we investigate how the properties of conformal
fields studied in Section~\ref{h13} lead to the notion of a chiral algebra,
study vertex commutation relations in the case of a chiral algebra. We
also show that Laurent coefficients of the operator $ \psi $, $ \psi^{+} $ (constructed in
Section~\ref{s2.40}) satisfy the same (anti)commutation relations as fermion
destruction/creation operators, and that Laurent coefficients of
operators $ {\mathcal G}\left(r,s\right) $ (constructed in Section~\ref{s12.30}) (and their smooth
parts) generate Lie algebras $ {\mathfrak g}{\mathfrak l}_{\infty} $ and its central extension $ {\mathfrak g}\left(\infty\right) $.

Note that the calculation in the homology space of a vector space
with several removed hyperspaces coincide literally with the calculation
done in for conformal fields studied in physics-related literature, and
in the approach of formal power series. However, they have a different
meaning, since we perform them with actual correctly defined families of
operators, instead of some symbols which depend on a parameter in a
formal sense only.

Consider conformal fields $ {\mathbit C}_{1},\dots ,{\mathbit C}_{k} $ on the same variety $ S $ acting in
the same space $ {\mathfrak M}{\mathcal M}\left(Z\right) $. Say that these conformal fields are {\em conformally
based on\/} $ \left(S,S^{\left(1\right)},S^{\left(2\right)}\right) $, $ S^{\left(1\right)}\subset S $, $ S^{\left(2\right)}\subset S\times S $, if all $ {\mathbit C}_{i} $ are smooth in $ S\smallsetminus S^{\left(1\right)} $,
and pairwise products of $ {\mathbit C}_{i} $ and $ {\mathbit C}_{j} $ have singularities at components of
$ S^{\left(2\right)} $ only. If one can take $ S^{\left(1\right)}=\varnothing $, $ S^{\left(2\right)}=\varnothing $, call the fields {\em conformally
trivial}.

We start with the simplest case when $ S^{\left(2\right)} $ is the diagonal $ \Delta=\Delta_{S} $ of $ S $.

\subsection{Chiral algebras }\label{s5.24}\myLabel{s5.24}\relax  Let $ S $ be a curve. Consider clean codominant
conformal fields $ {\mathbit C}\left(s\right) $, $ {\mathbit C}'\left(s\right) $, $ s\in S $, acting on $ {\mathfrak M}{\mathcal M}\left(Z\right) $ which are conformally
based on $ \left(S,\varnothing,\Delta\right) $. By definition one can write\footnote{To avoid defining what $ s_{i}-s_{j} $ means, one may assume that in this
expression $ {\mathbb K}={\mathbb C} $, and it is written in local analytic coordinate system on
$ S $, or that the expression is written in a {\em formal\/} coordinate system on a
completion of $ S $ in a particular point.
\endgraf
However, since $ {\mathbit D}^{\left(l\right)} $ depends on $ l $-jet of such coordinate systems
only, one can easily define it the algebraic case too.}
\begin{equation}
{\mathbit C}\left(r\right){\mathbit C}'\left(s\right) = \sum_{l\geq1}^{L}\frac{{\mathbit D}^{\left(l\right)}\left(r\right)}{\left(s-r\right)^{l}} + {\mathbit D}\left(r,s\right),
\label{equ4.20}\end{equation}\myLabel{equ4.20,}\relax 
here $ {\mathbit D}\left(\bullet,\bullet\right) $ has no pole at $ x\approx y $, and the conformal families $ {\mathbit D}^{\left(l\right)} $ are
smooth (otherwise $ {\mathbit C} $ and $ {\mathbit D} $ would have fusions on subvarieties of $ S\times S $ of
codimension 2, so they would not be based on $ \left(S,\bullet,\Delta\right) $). Note that we do
not exclude the case $ L=\infty $, but in this case $ {\mathbit D}^{\left(l\right)}F $ with $ l\gg 1 $ vanishes for
any fixed $ F $ (or $ F $ changing in a finite-dimensional family).

As Proposition~\ref{prop13.23} shows, the fields $ {\mathbit D}^{\left(l\right)} $ are also
codominant. As Remark~\ref{rem13.57} shows, for a generic fixed $ F $ (or $ F $
changing in a generic finite-dimensional family) and $ k\geq2 $ there is an
appropriate $ M $ such that
\begin{equation}
\left(\prod_{i<j}\left(s_{i}-s_{j}\right)\right)^{M}{\mathbit A}_{1}\left(s_{1}\right){\mathbit A}_{2}\left(s_{2}\right)\dots {\mathbit A}_{k}\left(s_{k}\right)F
\notag\end{equation}
is smooth in $ s_{i} $, here $ {\mathbit A}_{i} $ is $ {\mathbit C} $ or $ {\mathbit C}' $. In particular, the only non-zero
residues of $ {\mathbit A}_{1}\left(s_{1}\right){\mathbit A}_{2}\left(s_{2}\right)\dots {\mathbit A}_{k}\left(s_{k}\right) $ are on subvarieties given by
several equations $ s_{i_{m}}=s_{j_{m}} $.

Therefore, fusions of $ {\mathbit D}^{\left(l\right)}\left(r\right) $ with $ {\mathbit C}\left(s\right) $, $ {\mathbit C}'\left(s\right) $, or with other
$ {\mathbit D}^{\left(m\right)}\left(s\right) $ are all smooth, and defined on the diagonal $ r=s $. Same is true
for fusions of these fusions, and so on. One obtains a huge family of
conformal fields conformally based on $ \left(S,\varnothing,\Delta\right) $. Following \cite{BeiDriChiral}
and \cite{HuaLep96Mod}, introduce chiral algebras.

\begin{definition} A {\em chiral algebra\/} on $ S $ is a vector space $ {\mathcal C} $ of
conformal families on $ S $ acting in $ {\mathfrak M}{\mathcal M}\left(Z\right) $ such that
\begin{enumerate}
\item
If $ {\mathbit C}\left(s\right)\in{\mathcal C} $ and $ \Phi\left(s\right) $ is a holomorphic function on $ S $, then $ \Phi\left(s\right){\mathbit C}\left(s\right)\in{\mathcal C} $;
\item
any fusion of any two families $ {\mathbit C}_{1,2}\in{\mathcal C} $ is defined on $ \Delta_{S}\simeq S $ and is an
element of $ {\mathcal C}. $\footnote{The standard definition of the chiral algebra has an additional
structure: a linear mapping from vector fields on $ X $ to smooth families in
$ {\mathcal C} $.
\endgraf
Thus to any vector field on $ X $ one can associate a smooth family
from the vector space.}
\end{enumerate}
\end{definition}

It is possible to codify this definition by defining an {\em abstract
chiral algebra\/} in terms of operations of taking fusion only, without any
reference to the action in $ {\mathfrak M}{\mathcal M}\left(Z\right) $. After such a formalization one can
consider the object in the above definition as a {\em representation\/} of an
abstract chiral algebra in the lattice of vector spaces $ {\mathfrak M}{\mathcal M}\left(U\right) $, $ U\subset Z $,
compare \cite{Bor97Ver,BeiDriChiral}.

Given a chiral algebra spanned by families $ {\mathbit C}_{i} $, $ i\in I $, one can interpret~%
\eqref{equ4.20} as {\em chiral commutation relations\/}
\begin{equation}
{\mathbit C}_{i}\left(r\right){\mathbit C}_{j}\left(s\right) = \sum\Sb1\leq l\leq L \\ m\endSb \Gamma_{ij}^{lm}\left(r\right)\frac{{\mathbit C}_{m}\left(r\right)}{\left(s-r\right)^{l}} + \left[\text{ a smooth family near }r=s\,\right].
\notag\end{equation}
\subsection{Conformal commutation relations } Show how standard calculations with
chiral algebras look in our context.

Consider a curve $ S $ and two conformal fields $ {\mathbit C}\left(s\right) $, $ {\mathbit C}'\left(s\right) $ based on
$ \left(S,\varnothing,\Delta\right) $. Suppose that
\begin{equation}
{\mathbit C}\left(s\right){\mathbit C}'\left(s'\right) = q{\mathbit C}'\left(s'\right){\mathbit C}\left(s\right),\qquad q\in{\mathbb K},
\notag\end{equation}
and that $ {\mathbit C} $, $ {\mathbit C}' $ have chiral commutation relations
\begin{equation}
{\mathbit C}\left(r\right){\mathbit C}'\left(s\right) = \sum_{m\geq1}^{L}\Gamma^{\left(m\right)}\left(r\right)\frac{{\mathbit C}^{\left(m\right)}\left(r\right)}{\left(s-r\right)^{m}} + \left[\text{ a smooth family near }r=s\,\right].
\notag\end{equation}
Consider a smooth point $ s_{0}\in S $ and Laurent coefficients $ {\mathbit C}_{k} $, $ {\mathbit C}'_{k'} $ of these
fields near $ s_{0} $. Denote by $ {\mathbit C}_{k}^{\left(m\right)} $ the Laurent coefficients of the families
$ {\mathbit C}^{\left(m\right)}\left(s\right) $, by $ \Gamma_{k}^{\left(m\right)} $ the Laurent coefficients of the functions $ \Gamma^{\left(m\right)} $.

\begin{proposition} These Laurent coefficients have commutation relations
\begin{equation}
{\mathbit C}_{k}{\mathbit C}'_{k'}-q{\mathbit C}_{k'}'{\mathbit C}_{k} = -\sum\Sb1\leq m\leq L \\ t\in{\mathbb Z}\endSb \Gamma_{t}^{\left(m\right)}\binom{-k'-1}{m-1}{\mathbit C}_{\left(k+k'+m-t\right)}^{\left(m\right)},
\label{equ55.11}\end{equation}\myLabel{equ55.11,}\relax 
\end{proposition}

\begin{proof} By Theorem~\ref{th13.80}, $ {\mathbit C}_{k}{\mathbit C}'_{k'} $ is $ \left(k,k'\right) $-Laurent coefficient of
$ {\mathbit C}\left(r\right){\mathbit C}'\left(s\right) $ w.r.t.~the flag $ {\mathcal F}_{1}=\left\{\left(s_{0},s_{0}\right)\right\}\subset S\times\left\{s_{0}\right\}\subset S\times S\ni\left(r,s\right) $, $ {\mathbit C}'_{k'}{\mathbit C}_{k'} $ is
$ \left(k',k\right) $-Laurent coefficient of $ {\mathbit C}'\left(s\right){\mathbit C}\left(r\right) $ w.r.t.~the flag
$ {\mathcal F}_{2}=\left\{\left(s_{0},s_{0}\right)\right\}\subset\left\{s_{0}\right\}\times S\subset S\times S\ni\left(r,s\right) $.

Thus the difference in the left-hand side of~\eqref{equ55.11} is the
difference of Laurent coefficients of $ {\mathbit C}\left(r\right){\mathbit C}\left(s\right) $ w.r.t.~the flags $ {\mathcal F}_{1} $ and
$ {\mathcal F}_{2} $. As it was seen in Section~\ref{s13.60}, if $ {\mathbb K}={\mathbb C} $, one can write the
Laurent coefficients as integrals along the torus $ |s|=|r|^{K} $ in the case of
$ {\mathcal F}_{1} $, and $ |r|=|s|^{K} $ in the case of $ {\mathcal F}_{2} $, in both cases with $ K\gg 1 $, $ |s|=\varepsilon\ll 1 $.

Note that if $ \gamma>\alpha>\beta $ and $ \delta<\alpha $
\begin{equation}
\int_{|r|=\alpha}\int_{|s|=\beta}-\int_{|r|=\alpha}\int_{|s|=\gamma}=\int_{|r|=\alpha}\left(\int_{|s|=\beta}-\int_{|s|=\gamma}\right)=-\int_{|r|=\alpha}\int_{|s-r|=\delta}
\notag\end{equation}
in homology of $ U={\mathbb C}^{2}\smallsetminus\left(\left\{r=0\right\}\cup\left\{s=0\right\}\cup\left\{r=s\right\}\right) $. Thus for any holomorphic
function $ \Phi\left(r,s\right) $ on $ U $ the difference of Laurent coefficients of $ \Phi $ w.r.t.~
$ {\mathcal F}_{1} $ and $ {\mathcal F}_{2} $ is
\begin{equation}
-\frac{1}{\left(2\pi i\right)^{2}}\int_{|r|=\alpha}\int_{|s-r|=\delta}\Phi\left(r,s\right)\frac{dr}{r^{k+1}}\frac{ds}{s^{k'+1}}.
\notag\end{equation}
Since $ \delta $ may be arbitrarily small, one should be able to express the
integral $ \int_{|y-x|=\delta} $ via the singular part of $ \Phi $ on $ r=s $ only. Indeed, the part
$ \frac{dr}{r^{k+1}}\frac{ds}{s^{k'+1}} $ is smooth near $ r=s $, $ r\not=0 $, so only the poles of $ \Phi\left(r,s\right) $
along $ r=s $ matter. Suppose that $ \Phi\left(r,s\right)=\sum_{1\leq l\leq L}\frac{\Phi^{\left[l\right]}\left(r\right)}{\left(s-r\right)^{l}} $, with
Laurent coefficients of $ \Phi^{\left[l\right]}\left(r\right) $ being $ \Phi_{k}^{\left[l\right]} $.

Writing the holomorphic on $ \left\{\left(r,s\right) \mid r\approx s, r\not=0\right\} $ differential form
$ \frac{dr}{r^{k+1}}\frac{ds}{s^{k'+1}} $ as a Taylor series in $ \left(s-r\right) $ for a fixed $ r $, one gets
\begin{equation}
\frac{dr}{r^{k+1}}\frac{ds}{s^{k'+1}} =\sum_{m\geq0}\binom{-k'-1}{m}\frac{dr}{r^{k+k'+2+m}}\left(s-r\right)^{m}d\left(s-r\right)
\notag\end{equation}
thus the difference of Laurent coefficients of $ \Phi $ is
\begin{equation}
-\frac{1}{2\pi i}\int_{|r|=\alpha}\sum_{m\geq0}\binom{-k'-1}{m}\Phi^{\left[m+1\right]}\left(r\right)\frac{dr}{r^{k+k'+2+m}}=
-\sum_{m\geq0}\binom{-k'-1}{m}\Phi_{k+k'+1+m}^{\left[m+1\right]}
\notag\end{equation}
Applying this to $ \Phi\left(r,s\right)={\mathbit C}\left(r\right){\mathbit C}'\left(s\right)F $, one obtains
\begin{equation}
{\mathbit C}_{k}{\mathbit C}'_{k'}-q{\mathbit C}_{k'}'{\mathbit C}_{k} = -\sum\Sb0\leq m\leq L-1 \\ t\in{\mathbb Z}\endSb \Gamma_{t}^{\left(m+1\right)}\binom{-k'-1}{m}{\mathbit C}_{\left(k+k'+m-t+1\right)}^{\left(m+1\right)}.
\notag\end{equation}

In the case of generic $ {\mathbb K} $ one gets the same answer, since homology of
$ S\times S\smallsetminus\left(\left\{s_{0}\right\}\times S\cup S\times\left\{s_{0}\right\}\cup\Delta\right) $ does not depend on the field. \end{proof}

If the fields $ {\mathbit C}^{\left(m\right)}\left(s\right) $ are independent, one can incorporate
coefficients $ \Gamma^{\left(m\right)}\left(s\right) $ into the fields, and assume $ \Gamma^{\left(m\right)}\left(s\right) =1 $. Then
\begin{equation}
{\mathbit C}_{k}{\mathbit C}'_{k'}-q{\mathbit C}_{k'}'{\mathbit C}_{k} = -\sum_{1\leq m\leq L} \binom{-k'-1}{m-1}{\mathbit C}_{\left(k+k'+m\right)}^{\left(m\right)}.
\notag\end{equation}
\begin{example} \label{ex14.30}\myLabel{ex14.30}\relax  Consider fields $ \psi\left(s\right) $, $ \psi^{+}\left(s\right) $. Denote the Laurent
coefficients $ \psi_{k} $, $ \psi_{k}^{+} $. One concludes that the anticommutators satisfy
\begin{equation}
\left[\psi_{k},\psi_{k'}\right]_{+}=0,\qquad \left[\psi_{k}^{+},\psi_{k'}^{+}\right]_{+}=0,\qquad \left[\psi_{k},\psi_{k'}^{+}\right]_{+}=-\delta_{k+k'+1},
\label{equ6.12}\end{equation}\myLabel{equ6.12,}\relax 
since the only fusion of these fields is {\bf1}, which does not depend on $ s $,
thus has all the Laurent coefficients but the $ 0 $-th one being 0.

One can easily check the sign in this identity, since by Corollary~%
\ref{cor12.40}
\begin{equation}
\psi_{-1}\left(1\right)=0,\quad \psi_{0}\left(1\right)=T,\quad \psi_{0}^{+}\left(1\right)=T^{-1},\quad \psi_{-1}\left(T^{-1}\right)=-T^{-1}\psi_{0}\left(1\right)=-1.
\notag\end{equation}
\end{example}

\subsection{The family $ {\protect \mathcal G} $ }\label{s6.30}\myLabel{s6.30}\relax  Consider a conformal family $ {\mathcal G}\left(r,s\right) $ on $ {\mathbb A}^{2} $ defined in
Section~\ref{s12.30}:
\begin{equation}
\left({\mathcal G}\left(r,s\right)F\right)\left(\xi\right) = \frac{1}{s-r}\frac{\xi\left(r\right)}{\xi\left(s\right)}F\left(\frac{z-r}{z-s}\xi\left(z\right)\right).
\label{equ5.35}\end{equation}\myLabel{equ5.35,}\relax 
As calculations in Section~\ref{s2.40} show, $ {\mathcal G}\left(r,s\right)=\psi\left(r\right)\psi^{+}\left(s\right) $.

As it was seen in Section~\ref{s12.30}, the family $ {\mathcal G}\left(r,s\right) $ is
based on
$ \left({\mathbb A}^{2}, D_{0}, D_{1}\right) $ with $ D_{0}=\Delta_{{\mathbb A}^{1}}\subset{\mathbb A}^{2} $, $ D_{1}=\left\{\left(r_{1},s_{1},r_{2},s_{2}\right) \mid s_{1}\not=r_{2}, s_{2}\not=r_{1}\right\} $. The
residue of $ {\mathcal G}\left(r,s\right) $ on $ r=s $ is {\bf1}. In fact the definition of $ {\mathcal G}\left(r,s\right) $ shows
that in fact $ {\mathcal G}\left(r,s\right) $ can be considered as a family based on
$ \left(\left({\mathbb P}^{1}\right)^{2},\Delta_{{\mathbb P}^{1}},\bar{D}_{1}\right) $, here $ \bar{D}_{1}=\left\{\left(r_{1},s_{1},r_{2},s_{2}\right) \mid s_{1}\not=r_{2}, s_{2}\not=r_{1}\right\}\subset\left({\mathbb P}^{1}\right)^{4} $.

Consider Laurent coefficients of $ {\mathcal G}\left(r,s\right) $ for $ r,s\approx0 $. Consider two
different cases, one corresponds to the flag $ {\mathcal F}^{<}=\left\{\left(0,0\right)\right\}\subset{\mathbb A}^{1}\times\left\{0\right\}\subset{\mathbb A}^{2} $ (i.e.,
taking a contour with $ |r|<|s| $), another to the flag $ {\mathcal F}^{>}=\left\{\left(0,0\right)\right\}\subset\left\{0\right\}\times{\mathbb A}^{1}\subset{\mathbb A}^{2} $
(or to $ |r|>|s| $). Denote the corresponding Laurent coefficients $ {\mathcal G}_{nm}^{<} $ and
$ {\mathcal G}_{nm}^{>} $ (assume that $ n $ gives the power of $ r $, $ m $ the power of $ s $). In the
same way as we did it with $ \psi $ and $ \psi^{+} $, one gets commutation relations for
$ {\mathcal G}_{nm}^{<,>} $ (a more complicated case of $ {\mathcal G}_{nm} $ is done in details below):
\begin{align} \left[{\mathcal G}_{nm}^{<}, {\mathcal G}_{kl}^{<}\right] & = {\mathcal G}_{nl}^{<}\delta_{m+k+1} - {\mathcal G}_{km}^{<}\delta_{n+l+1},
\notag\\
\left[{\mathcal G}_{nm}^{>}, {\mathcal G}_{kl}^{>}\right] & = {\mathcal G}_{nl}^{>}\delta_{m+k+1} - {\mathcal G}_{km}^{>}\delta_{n+l+1}.
\notag\end{align}

On the other hand, since $ {\mathcal G}\left(r,s\right) $ has a simple pole at $ r\approx s $ with
residue $ -\boldsymbol1 $, one can write the integrals for $ {\mathcal G}_{nm}^{<} - {\mathcal G}_{nm}^{>} $ as
\begin{equation}
\int_{|r|=\alpha}\int_{|s|=\beta} - \int_{|r|=\gamma}\int_{|s|=\beta} = \int_{|r-s|=\alpha'}\int_{|s|=\beta},\qquad \alpha,\alpha'<\beta<\gamma,
\notag\end{equation}
thus $ {\mathcal G}_{nm}^{<} - {\mathcal G}_{nm}^{>} = \delta_{n+m+1}\boldsymbol1 $. To understand the relationship between $ {\mathcal G}_{nm}^{<} $
and $ {\mathcal G}_{nm}^{>} $ better, consider the smooth part $ {\mathcal G}^{\text{smooth}}\left(r,s\right) $ of $ {\mathcal G}\left(r,s\right) $:
\begin{equation}
{\mathcal G}^{\text{smooth}}\left(r,s\right) = {\mathcal G}\left(r,s\right) - \frac{1}{s-r}\boldsymbol1.
\notag\end{equation}
Denote by $ {\mathcal G}_{nm} $ the Laurent coefficients of the smooth function $ {\mathcal G}^{\text{smooth}} $
(they do not depend on the choice of the flag). Writing Laurent series
for $ \frac{1}{r-s} $ near flags $ {\mathcal F}_{1} $ and $ {\mathcal F}_{2} $ (i.e., in regions $ |r|<|s| $ and $ |r|>|s| $),
one gets
\begin{equation}
{\mathcal G}_{nm} = {\mathcal G}_{nm}^{<} - \delta_{n+m+1}\chi_{n\geq0} = {\mathcal G}_{nm}^{>} + \delta_{n+m+1}\chi_{n<0}.
\notag\end{equation}
Here $ \chi_{S} $ is the characteristic function of a set $ S $. Using these
relationship, one can easily get the commutation relations for $ {\mathcal G}_{nm} $:
\begin{equation}
\left[{\mathcal G}_{nm}, {\mathcal G}_{kl}\right] = {\mathcal G}_{nl}\delta_{m+k+1} - {\mathcal G}_{km}\delta_{n+l+1} + 2\delta_{n+l+1}\delta_{m+k+1}\left(\chi_{n\geq0}-\chi_{k\geq0}\right).
\notag\end{equation}

Indeed, the divisor of poles for $ {\mathcal G}\left(a,b\right){\mathcal G}\left(c,d\right)=-\psi\left(a\right)\psi\left(c\right)\psi^{+}\left(b\right)\psi^{+}\left(d\right) $ is
\begin{equation}
\frac{1}{b-a}\psi\left(c\right)\psi^{+}\left(d\right) - \frac{1}{b-c}\psi\left(a\right)\psi^{+}\left(d\right) -\frac{1}{d-a}\psi\left(c\right)\psi^{+}\left(b\right)+ \frac{1}{d-c}\psi\left(a\right)\psi^{+}\left(b\right),
\notag\end{equation}
since $ {\mathcal G}^{\text{smooth}}\left(a,b\right){\mathcal G}^{\text{smooth}}\left(c,d\right) $ coincides with
\begin{equation}
{\mathcal G}\left(a,b\right){\mathcal G}\left(c,d\right)-\frac{1}{b-a}{\mathcal G}^{\text{smooth}}\left(c,d\right) - \frac{1}{d-c}{\mathcal G}^{\text{smooth}}\left(a,b\right)-\frac{1}{\left(b-a\right)\left(d-c\right)}
\notag\end{equation}
the poles coincide with poles of
\begin{equation}
\frac{1}{\left(b-a\right)\left(d-c\right)} -\frac{2}{\left(b-c\right)\left(d-a\right)}-\frac{1}{b-c}{\mathcal G}^{\text{smooth}}\left(a,d\right)-\frac{1}{d-a}{\mathcal G}^{\text{smooth}}\left(c,b\right).
\notag\end{equation}
To calculate the commutation relations for Laurent coefficients one should
take difference of integrals along $ |a|\approx|b|<|c|\approx|d| $ and $ |a|\approx|b|>|c|\approx|d| $.
Take contours $ C_{1} $ with $ |a|<|b|<|c|<|d| $ and $ C_{2} $ with $ |c|<|d|<|a|<|b| $.

When pulling $ C_{1} $ to $ C_{2} $ preserving $ |a|<|b| $ and $ |c|<|d| $, we cross poles
at $ b=c $ and at $ a=d $, thus the first summand out of four gives no
contribution at all. During intersection of $ b=c $ the fourth one gives no
contribution, and during intersection of $ a=d $ the third gives no
contribution.

The pole at $ b=c $ has residue $ \frac{2}{d-a}+{\mathcal G}^{\text{smooth}}\left(a,d\right) $ w.r.t.~$ c-b $, thus
its contribution to $ \left[{\mathcal G}_{nm}, {\mathcal G}_{kl}\right] $ is
\begin{equation}
\frac{1}{\left(2\pi i\right)^{2}}\int_{|a|<|d|}\frac{da}{a^{n+1}}\frac{dd}{d^{l+1}}\left(\frac{2}{d-a}+{\mathcal G}^{\text{smooth}}\left(a,d\right)\right)\delta_{k+m+1},
\notag\end{equation}
i.e., it is $ \delta_{k+m+1} $ times the Laurent coefficient of
$ \frac{2}{d-a}+{\mathcal G}^{\text{smooth}}\left(a,d\right) $. Note that during this intersection $ |a|<|d| $,
thus the contribution is
\begin{equation}
+2\chi_{n\geq0}\delta_{k+m+1}\delta_{n+l+1}+\delta_{k+m+1}{\mathcal G}_{nl}.
\notag\end{equation}
Similarly, the pole at $ a=d $ gives contribution of $ \delta_{n+l+1} $ times the Laurent
coefficient of $ -\frac{2}{b-c}-{\mathcal G}^{\text{smooth}}\left(c,b\right) $, and during this intersection
$ |b|>|c| $. Thus this pole gives
\begin{equation}
-2\chi_{k\geq0}\delta_{k+m+1}\delta_{n+l+1}-\delta_{n+l+1}{\mathcal G}_{km}.
\notag\end{equation}
This proves commutation relations for $ {\mathcal G}_{ab} $.

\subsection{Generators $ X_{ij} $ }\label{s6.31}\myLabel{s6.31}\relax  Let $ X_{ij}={\mathcal G}_{i\,\,-j-1} $, same with $ X_{ij}^{<} $ and $ X_{ij}^{>} $. It is
obvious now that $ X_{ij}^{<} $ commute as elementary matrices of infinite size
(same with $ X_{ij}^{>} $), thus generate $ {\mathfrak g}{\mathfrak l}_{\infty} $, and $ X_{ij} $ commute as generators of central
extension $ {\mathfrak g}\left(\infty\right) $ of $ {\mathfrak g}{\mathfrak l}_{\infty} $.

Fix $ F\in{\mathfrak M}{\mathcal M}\left(Z\right) $. Note that $ X_{ij}^{<}F $ is a Laurent coefficient of the family
$ {\mathcal G}\left(r,s\right)F $. This function is smooth if $ r,s\notin\operatorname{Supp} F $ and $ r\not=s $. Since $ \operatorname{Supp} F $ may
contain 0, one can conclude that $ X_{ij}^{<}F $ is a Laurent coefficient of a
function which is smooth in $ 0<|r|<|s|\ll 1 $. Similarly, $ X_{ij}F $ is a Laurent
coefficient of a function which is smooth in $ 0<|r|,|s|\ll 1 $.

Consider a collection $ a_{ij} $, $ i,j\in{\mathbb Z} $, such that the sum $ \sum a_{ij}\phi_{ij} $ contains
only a finite number of non-zero summands if $ \phi_{ij} $ are Laurent coefficients
of a meromorphic function $ \phi $ which is smooth in $ 0<|r|<|s|\ll 1 $. Then $ a_{ij}=0 $
if $ j>J\left(-i\right)-i $ or $ i>I\left(-i-j\right) $, here $ J $ and $ I $ are appropriate non-decreasing
functions. If the same is true if $ \phi $ is smooth in $ 0<|r|<|s|\ll 1 $, then $ a_{ij}=0 $
if $ j>J\left(-i\right) $ or $ i>I\left(-j\right) $. This is compatible with the facts that $ {\mathfrak g}\left(\infty\right) $
contains arbitrary sums $ \sum a_{ij}X_{ij} $, with $ a_{ij}=0 $ if $ j<J\left(-i\right) $ or $ i>I\left(j\right) $, and
that cocycle for $ {\mathfrak g}\left(\infty\right) $ is a coboundary on sums $ \sum a_{ij}X_{ij} $, with $ a_{ij}=0 $ if
$ j<-J\left(-i\right)+i $ or $ i>I\left(j-i\right) $.

One can construct geometric realizations of some particular infinite
sums $ \sum a_{ij}X_{ij} $. Consider
\begin{equation}
Y_{kl} = \sum_{i\in{\mathbb Z}}\binom{i}{k}X_{i\,\,i-k-l}=\sum\binom{i}{k}{\mathcal G}_{i\,\,-i+k+l-1}^{\text{smooth}},\qquad k\geq0.
\notag\end{equation}
Taking into account that $ {\mathcal G}_{ab} $ are Laurent coefficients for $ {\mathcal G}\left(r,s\right) $, one can
see that
\begin{equation}
Y_{kl} = \frac{1}{\left(2\pi i\right)^{2}}\sum_{i}\binom{i}{k}\iint{\mathcal G}\left(r,s\right)r^{-i-1}s^{i-k-l}dr\,ds
= \frac{1}{2\pi i}\int{\mathcal G}^{\text{smooth}\,\,\left(k,0\right)}\left(s,s\right)\frac{ds}{s^{l}},
\notag\end{equation}
so $ Y_{kl} $ is the $ \left(l-1\right) $-th Laurent coefficient of $ -k $-th residue of $ \psi\left(r\right)\psi^{+}\left(s\right) $
at $ r=s $ (the vector field $ v $ being $ \partial/\partial r $).

Let us calculate Laurent coefficients of $ 0 $-th residue $ B\left(s\right) $ of
$ \psi\left(r\right)\psi^{+}\left(s\right) $ at $ r=s $. By~\eqref{equ5.35}
\begin{equation}
\left(B\left(s\right)F\right)\xi = -\frac{\xi'\left(s\right)}{\xi\left(s\right)}F\left(\xi\left(z\right)\right) - \frac{dF\left(\left(1+\frac{\tau}{z-s}\right)\xi\left(z\right)\right)}{d\tau}|_{\tau=0}.
\notag\end{equation}
In coordinates
\begin{equation}
\xi\left(z\right)=Te^{\sum_{k}t_{k}z^{k}}
\notag\end{equation}
one can write $ B $ as
\begin{equation}
B\left(s\right)F\left(T,t\right)= -\left(\sum_{l\geq1}lt_{l}s^{l-1}\right)F\left(T,t\right)+\sum_{l\geq1}\frac{1}{s^{l+1}}\partial_{t_{l}}F\left(T,t\right)+\frac{1}{s}\frac{\partial_{T}}{T}F\left(T,t\right).
\notag\end{equation}
We conclude that $ Y_{0l}=-lt_{l} $ if $ l>0 $, $ Y_{0l}=\partial_{t_{-l}} $ if $ l<0 $, and $ Y_{00}=\frac{\partial_{T}}{T} $.

Taking the ``holomorphic in $ |s|\leq1 $ part'' of the function $ B\left(s\right) $, one
gets
the the standard boson creation operator in the Fock space $ {\mathbb K}\left[t_{1},t_{2},\dots \right] $,
and the ``holomorphic in $ |s|\geq1 $ part'' gives the boson destruction operator.
Moreover, one can see that $ B\left(s\right) $ lies in the completion of algebra
generated by $ \psi\left(s\right) $ and $ \psi^{+}\left(s\right) $. On the other hand, the formulae for vertex
operators show that visa versa, $ \psi\left(s\right) $ and $ \psi^{+}\left(s\right) $ are in the completion of
algebra generated by Laurent coefficients of $ B\left(s\right) $.

Note that Laurent coefficients of $ \psi\left(s\right) $, $ \psi^{+}\left(s\right) $ commute as fermion
destruction/creation operators (see~\eqref{equ6.12}), and Laurent
coefficients
of $ B\left(s\right) $ commute as $ t_{k} $, $ \frac{\partial}{\partial t_{k}} $, i.e., as boson creation/destruction
operators. This the base of so called {\em boson-fermion correspondence}, which
is studied in the following section.

\section{Boson-fermion correspondence }\label{h15}\myLabel{h15}\relax 

In this section we establish a {\em geometric boson-fermion
correspondence\/} (in dimension $ 1+1 $), show that this correspondence provides
heuristic for considering $ {\mathcal M}\left({\mathcal M}\left({\mathbb P}^{1}\right)\right) $ as a configuration space for a system
of fermions close to Dirac's vacuum state, show that the Fock space
coincides with the localization of $ {\mathcal M}\left({\mathcal M}\left({\mathbb P}^{1}\right)\right) $ near $ 0\in{\mathbb P}^{1} $, and investigate
the relation of the geometric boson-fermion correspondence with the
usual, algebraic boson-fermion correspondence.

The mechanical description of a physical system consists of two
parts: kinematics describes the set of possible {\em states\/} of the system, and
dynamics examines how the states change when time goes. As far as we are
concerned here, dynamics is not a question, so we restrict our attention
to kinematics only.

The quantum mechanics describes {\em states\/} of a system as
elements of a Hilbert space called {\em(quantum) configuration space}. In
the case of $ d $ space dimensions, this
space may be the Hilbert space of square-integrable functions on $ {\mathbb R}^{d} $.
A particle ``concentrated'' near a point $ Y\in{\mathbb R}^{d} $ is described by a function
with a support in a neighborhood of $ Y $. Two proportional vectors represent
the same state of the system.

In general the relation to the classical mechanics is that the
quantum configuration space is the vector space of functions on {\em classical
configuration space\/} ($ {\mathbb R}^{d} $ above), and functions with support in a
small neighborhood of the given point represent quantum states which are
``almost classical''.

In what follows we ignore the Hilbert structure completely (though
it can be added back later at some moment). This would allow one to
consider more general functions than square-integrable one, say, one can
consider functions of the form
\begin{equation}
\sum_{k=1}^{K}a_{k}\delta\left(Y-Y_{k}\right),
\label{equ3.10}\end{equation}\myLabel{equ3.10,}\relax 
vectors $ Y_{i} $ being some points in $ {\mathbb R}^{d} $. If $ K=1 $, the above state describes a
particle concentrated {\em exactly\/} at $ Y_{1} $, thus a classical state.

\subsection{Bosonic space } What is described above is the approach of {\em primary
quantization}, when the number of particles is fixed, thus we can model $ m $
particles in $ d $-dimensional space by one particle in $ md $-dimensional space.
To allow the number of particle to be variable, one needs some
modifications.

In what follows we assume $ d=1 $.

If the quantum system is a union of two parts, each one described by
a (Hilbert) space $ H_{1,2} $, then the system in whole is described by the
tensor product $ H_{1}\otimes H_{2} $, vector $ v_{1}\otimes v_{2} $ describing the state of the system
when the first part is in the state $ v_{1} $, the second in the state $ v_{2} $. Thus
the system of $ n $ particles having the space of states $ V $ each is described
by
\begin{equation}
\underbrace{V\otimes\dots \otimes V}_{n\text{ times}}.
\notag\end{equation}
If the particles are indistinguishable one needs to require that $ v_{1}\otimes v_{2} $ and
$ v_{2}\otimes v_{1} $ describe the {\em same\/} state of the system, thus in this case the space
of states of the system is some quotient of the above space. In this
quotient the images of $ v_{1}\otimes v_{2} $ and $ v_{2}\otimes v_{1} $ should be proportional.

A little bit of representation theory of the symmetric group $ {\mathfrak S}_{n} $ says
that there are only two such quotients: the symmetric power $ S^{n}V $ of $ V $ and
the skewsymmetric one $ \Lambda^{n}V $. Which one is relevant is determined by the
physical properties of the particle in question. Particles for which $ S^{n} $
should be taken are called {\em bosons}, the other ones {\em fermions}.

Summing up, we conclude that the configuration space of $ n $ bosons on a
line $ {\mathbb R} $ consists of symmetric functions on $ {\mathbb R}^{n} $. If one consider states of
one boson of the form~\eqref{equ3.10}, one should consider the space $ B_{n} $ which
consists of symmetric generalized functions\footnote{I.e., generalized functions which are invariant w.r.t.~permutations of
coordinates.} on $ {\mathbb R}^{n} $ which are sums of
$ \delta $-functions.

To consider variable number of particles, one should consider the
direct sum of spaces of the above form for all possible values of $ n $,
$ 0\leq n<\infty $:
\begin{equation}
{\mathbit B}=\sum_{k\geq0}B_{k}.
\notag\end{equation}
\subsubsection{The dual space } The above space is too convoluted to work with, thus
in fact we are going to describe {\em some\/} dual space to this one. Note that
there is a canonically defined non-degenerate pairing between the above
space of generalized functions on $ {\mathbb R}^{n} $ and the space of symmetric polynomials
on $ {\mathbb R}^{n} $. This pairing is the calculation of value
\begin{equation}
\left< P\left(Y\right), \delta\left(Y-Y_{0}\right) \right> = P\left(Y_{0}\right),\qquad Y\in{\mathbb R}^{n}.
\notag\end{equation}
\subsubsection{Quotient space } One can identify the space of symmetric polynomial
with the space of polynomials in new variables $ \left(\sigma_{1},\dots ,\sigma_{n}\right) $, which are
elementary symmetric functions in the old variables $ z_{i} $. Associate with
the collection $ \left(\sigma_{1},\dots ,\sigma_{n}\right) $ the polynomial
\begin{equation}
\pi\left(z\right) = z^{n}-\sigma_{1}z^{n-1}+\sigma_{2}z^{n-2}-\dots \pm\sigma_{n} = \left(z-z_{1}\right)\left(z-z_{2}\right)\dots \left(z-z_{n}\right).
\notag\end{equation}

Denote the manifold of polynomials of degree $ n $ with the leading
coefficient 1 by $ {\mathcal P}_{n} $. One concludes that the dual space to the bosonic
configuration space can be described as the space of polynomials on the
affine variety $ {\mathcal P}_{n} $. Similarly, the bosonic space itself can be identified
with the space spanned by $ \delta $-functions on $ {\mathcal P}_{n} $.

To consider variable number of particles, one should consider the
variety
\begin{equation}
{\mathcal P} = \bigcup_{n}{\mathcal P}_{n}.
\notag\end{equation}
\subsection{Boson creation operator } Fix a point $ z_{0}\in{\mathbb R} $. Consider the
transformation of adding a boson at the point $ z_{0} $. Obviously, if the
initial state of the system was described by $ \delta\left(Y-Y_{0}\right) $, $ Y,Y_{0}\in{\mathbb R}^{n} $, the new
state is $ \delta\left(\bar{Y}-\bar{Y}_{0}\right) $, $ \bar{Y},\bar{Y}_{0}\in{\mathbb R}^{n+1} $, $ \bar{Y}_{0}=\left(z_{0}, Y_{0}\right) $. In other words, it is
described by the mapping
\begin{equation}
{\mathbb R}^{n} \to {\mathbb R}^{n+1}\colon \left(z_{1},\dots ,z_{n}\right) \mapsto \left(z_{0},z_{1},\dots ,z_{n}\right)
\notag\end{equation}
(followed by symmetrization, of course).

Similarly, in the language of quotient space, this operator is
\begin{equation}
{\mathcal P}_{n} \to {\mathcal P}_{n+1}\colon \pi\left(z\right) \mapsto \left(z-z_{0}\right)\pi\left(z\right).
\notag\end{equation}
Multiplication law for polynomials shows that this mapping is compatible
with the affine structures on the spaces $ {\mathcal P}_{\bullet} $.

In the dual language, the dual operator sends symmetric polynomial on
$ {\mathbb R}^{n+1} $ to symmetric polynomials on $ {\mathbb R}^{n} $:
\begin{equation}
P\left(z, z_{1},\dots ,z_{n}\right) \mapsto Q\left(z_{1},\dots ,z_{n}\right),\quad Q\left(z_{1},\dots ,z_{n}\right) = P\left(z_{0},z_{1},\dots ,z_{n}\right).
\notag\end{equation}

In the dual language for the quotient, when one considers
polynomials $ f\left(\pi\right) $ on $ {\mathcal P} $, this mapping is written as
\begin{equation}
\alpha_{z_{0}}\colon f \mapsto \alpha_{z_{0}}f,\qquad \left(\alpha_{z_{0}}f\right)\left(\pi\right)=f\left(\left(z-z_{0}\right)\pi\right).
\label{equ3.25}\end{equation}\myLabel{equ3.25,}\relax 
\begin{remark} Note that this formula for the boson creation operator
coincides with the formula for the operator $ {\mathbit M}_{\xi} $ from Section~\ref{s11.31} in
the case $ \xi\left(z\right)=z-z_{0} $ and the rational function $ F $ being the polynomial $ \pi $.
Thus it is the restriction of operator $ a^{-1}\left(s\right) $ considered in Section~\ref{h0.2}. \end{remark}

\subsection{Naive boson-fermion correspondence }\label{s15.30}\myLabel{s15.30}\relax  Consider the configuration
space for the system of $ n $ fermions on the line $ {\mathbb R} $. As described above, it
is the space of skew-symmetric functions on $ {\mathbb R}^{n} $. To use the dual language,
note that there is a canonical pairing of the space of skew-symmetric
generalized functions on $ {\mathbb R}^{n} $ (with compact support) and the space of
skew-symmetric polynomials.

Thus we consider the space of skew-symmetric polynomials on $ {\mathbb R}^{n} $ as
the dual space to the configuration space. On the other hand, any
skew-symmetric polynomial $ q $ in $ z_{i} $ can be written as
\begin{equation}
q\left(z_{i}\right) = \operatorname{Disc}\left(z_{i}\right)p\left(z_{i}\right),\qquad \operatorname{Disc}\left(z_{i}\right) = \prod_{i<j}\left(z_{i}-z_{j}\right),
\notag\end{equation}
with $ p $ being a symmetric polynomial.

Call the correspondence $ p \iff q $ {\em the (geometric) boson-fermion
correspondence}. This correspondence makes it possible to describe states
of fermions in terms of states of bosons, via an explicit isomorphism of
configuration spaces.

\subsection{Fermion creation operator } Consider a transformation of configuration
space which corresponds to adding a fermion at the point $ z_{0}\in{\mathbb R} $. As in
the case of bosons, it is described by the mapping
\begin{equation}
{\mathbb R}^{n} \to {\mathbb R}^{n+1}\colon \left(z_{1},\dots ,z_{n}\right) \mapsto \left(z_{0},z_{1},\dots ,z_{n}\right),
\notag\end{equation}
the only difference being that it should be followed by
skewsymmetrization (instead of symmetrization). In the dual language, the
formula is the same as for bosons:
\begin{equation}
P\left(z, z_{1},\dots ,z_{n}\right) \mapsto Q\left(z_{1},\dots ,z_{n}\right),\quad Q\left(z_{1},\dots ,z_{n}\right) = P\left(z_{0},z_{1},\dots ,z_{n}\right),
\notag\end{equation}
but it is applied to skew-symmetric functions.
Let us rewrite this operator in terms of bosonic configuration
space using the boson-fermion correspondence. Let $ P_{\text{bos}} $, $ Q_{\text{bos}} $ be the
images of $ P $ and $ Q $ under boson-fermion correspondence:
\begin{equation}
P\left(z,z_{i}\right)=\operatorname{Disc}\left(z,z_{i}\right)\cdot P_{\text{bos}}\left(z,z_{i}\right),\qquad Q\left(z_{i}\right)=\operatorname{Disc}\left(z_{i}\right)\cdot Q_{\text{bos}}\left(z_{i}\right).
\notag\end{equation}
One can see that in terms of bosons this operator is
\begin{align} P_{\text{bos}}\left(z, z_{1},\dots ,z_{n}\right) & \mapsto Q_{\text{bos}}\left(z_{1},\dots ,z_{n}\right),
\notag\\
Q_{\text{bos}}\left(z_{1},\dots ,z_{n}\right) & = \frac{\operatorname{Disc}\left(z_{0},z_{i}\right)}{\operatorname{Disc}\left(z_{i}\right)} P_{\text{bos}}\left(z_{0},z_{1},\dots ,z_{n}\right)
\notag\\
& =\left(\prod_{k=1}^{n}\left(z_{0}-z_{k}\right)\right) P_{\text{bos}}\left(z_{0},z_{1},\dots ,z_{n}\right).
\notag\end{align}
On the other hand, in terms of dual language for the quotient, the factor
$ \left(\prod_{k=1}^{n}\left(z_{0}-z_{k}\right)\right) $ is just $ \pi\left(z_{0}\right) $, thus the fermion creation operator is
\begin{equation}
\psi_{z_{0}}^{+}\colon f \mapsto \psi_{z_{0}}^{+}f,\qquad \left(\psi_{z_{0}}^{+}f\right)\left(\pi\right)=\pi\left(z_{0}\right)f\left(\left(z-z_{0}\right)\pi\right).
\label{equ3.30}\end{equation}\myLabel{equ3.30,}\relax 
\begin{remark} Note that the last formula coincides with Equation~\eqref{equ12.31} of
Section~\ref{s2.40} for $ \psi^{+} $ in the case of the rational function $ \xi $ being the
polynomial $ \pi $. \end{remark}

\subsection{Fock space I}\label{s16.50}\myLabel{s16.50}\relax  The above description of the configuration space in
the fermion case is extremely naive, since it deals with collections of a
finite number of fermions. However, the physically meaningful picture
deals with collections close to the Dirac vacuum state.

This state corresponds to having infinitely many fermions, all the
fermions with ``negative energy'' being present, all with ``positive energy''
being absent. (This is not relevant to the discussion in this paper, but
we state it nevertheless: in terms of functions on $ {\mathbb R} $, we say that the
state $ f\left(z\right) $, $ z\in{\mathbb R} $, has negative energy if its Fourier transform $ \widehat{f}\left(\xi\right) $
vanishes for positive $ \xi $, similarly for positive energy. However, we use
``complex logarithmic'' coordinate system $ z=e^{i\widetilde{z}} $, thus will use Laurent
transform instead of Fourier transform.)

What follows may have no mathematical meaning, but it provides a
heuristic for relation of conformal field theory to the secondary
quantization in $ 1 $-dimensional case. Note that this heuristic is yet
incomplete, compare with Remark~\ref{rem3.60}.

Fix a finite collection $ S $ of points on the line $ {\mathbb R} $. The classical
states which correspond to bosons sitting at the points of $ S $ {\em and\/} at some
other points are described by the subset
\begin{equation}
p_{S}\cdot{\mathcal P}\subset{\mathcal P},\qquad p_{S}\left(z\right)=\prod_{s\in S}\left(z-s\right).
\notag\end{equation}
Clearly, the same is true for fermion states after applying boson-fermion
correspondence.

If we consider a classical state corresponding to $ S $ as a new {\em ground
state},\footnote{Say, one can suppose that the Hamiltonian for some physical system has
this collection as a minimum.} then the parameter change $ \pi=p_{S}\cdot\pi_{S} $ describes a new coordinate $ \pi_{S} $ on
the set of admissible states such that the ground state corresponds to
$ \pi_{S}=1 $.

Suppose that the collection $ S $ is very big, but one want to consider
classical states which are {\em close\/} to $ S $, say, add only a few additional
bosons or fermions. One finishes with the subset of classical
configuration space:
\begin{equation}
p_{S}\cdot{\mathcal P}_{\leq K}\subset{\mathcal P},\qquad {\mathcal P}_{\leq K}=\bigcup_{k\leq K}{\mathcal P}_{k}.
\notag\end{equation}
If we allow taking a few particles from the collection $ S $ as well as
adding a few particles, the set of classical states becomes
\begin{equation}
\left(p_{S}\cdot{\mathcal Q}_{K}\right) \cap {\mathcal P} = p_{S}\cdot{\mathcal Q}_{K}^{\left(S\right)} \subset {\mathcal P}.
\notag\end{equation}
Here $ {\mathcal Q}_{K} $ is the set of rational functions of the form $ \frac{p_{1}}{p_{2}} $, $ p_{1,2}\in{\mathcal P}_{\leq k} $,
and $ {\mathcal Q}_{K}^{\left(S\right)}={\mathcal Q}_{K}\cap\left(p_{S}^{-1}\cdot{\mathcal P}\right) $.

We conclude that if we consider classical states {\em close to\/} the ground
state $ S $ (in the above sense), then the corresponding subset of
configurations space can be identified with the set $ {\mathcal Q}_{K}^{\left(S\right)} $ for some small
value of $ K $. When we increase the set $ S $, the set $ {\mathcal Q}_{K}^{\left(S\right)} $ increases and
becomes closer and closer to $ {\mathcal Q}_{K} $. In other words, the set $ {\mathcal Q}_{K} $ describes
{\em relative\/} classical states of a configuration of particles with respect to
a {\em very big\/} collection $ S $, or, in other words, a very populated ground
state. The corresponding quantum states are described by functions on $ {\mathcal Q}_{K} $.

Let us describe boson and fermion creation operators in terms of
representation by functions on $ {\mathcal Q}_{K} $. Obviously, the boson creation operator $ \alpha_{z_{0}} $
is described by the same formula~\eqref{equ3.25}, thus the expression for this
operator does not depend on $ S $. To write down the fermion
creation operator $ \psi_{z_{0}}^{+} $, one need to fix the family $ S $. Writing down $ \psi_{z_{0}}^{+} $ in
terms of $ \pi=p_{S}\cdot\pi_{S} $, $ \pi_{S}\in{\mathcal Q}_{K}^{\left(S\right)} $, one gets
\begin{equation}
\psi_{z_{0}}^{+S}\colon f \mapsto \psi_{z_{0}}^{+S}f,\qquad \left(\psi_{z_{0}}^{+S}f\right)\left(\pi_{S}\right)=p_{S}\left(z_{0}\right)\pi_{S}\left(z_{0}\right)f\left(\left(z-z_{0}\right)\pi_{S}\right),\quad \pi_{S}\in{\mathcal Q}_{K}^{\left(S\right)},
\notag\end{equation}
or $ \psi_{z_{0}}^{+S} = p_{S}\left(z_{0}\right)\psi_{z_{0}}^{+} $. Since the operators $ \psi_{z_{0}}^{+} $ and $ \psi_{z_{0}}^{+S} $ differ by a
constant $ p_{S}\left(z_{0}\right) $, and vectors which differ by a constant describe the same
quantum state, one can safely use $ \psi_{z_{0}}^{+} $ instead of $ \psi_{z_{0}}^{+S} $.

It is obvious that the (extended as above) operator $ \alpha_{z_{0}} $ sends
polynomial functions on $ {\mathcal Q}_{K} $ to polynomial functions on $ {\mathcal Q}_{K+1} $. This is not
true for $ \psi_{z_{0}}^{+} $, since $ \pi_{S} \mapsto \pi_{S}\left(z_{0}\right) $ is not a {\em polynomial\/} function of
$ \pi_{S}\in{\mathcal Q}_{K}^{\left(S\right)} $. However, $ \pi_{S} \mapsto \pi_{S}\left(z_{0}\right) $ is a {\em rational\/} function of $ \pi_{S}\in{\mathcal Q}_{K}^{\left(S\right)} $, thus
the extended operator $ \psi_{z_{0}}^{+} $ sends rational functions on $ {\mathcal Q}_{K} $ to rational
functions on $ {\mathcal Q}_{K+1} $.

Obviously, the operators $ \alpha_{z_{0}} $, $ \psi_{z_{0}}^{+} $ are compatible with natural
inclusions $ {\mathcal Q}_{K} \hookrightarrow {\mathcal Q}_{K+L} $, $ L\geq0 $, thus one can extend them to rational
functions on
\begin{equation}
{\mathcal Q}= \lim _{\to} {\mathcal Q}_{K}.
\notag\end{equation}
Note that the description of $ {\mathcal Q} $ as of inductive limit leads to the same
definition of the vector space of rational functions on $ {\mathcal Q} $ as we had in
the Section~\ref{s11.10}, taking into account $ {\mathcal Q}={\mathcal M}\left({\mathbb P}^{1}\right) $.

\begin{nwthrmiii} Extend operators $ \alpha_{z_{0}} $ and $ \psi_{z_{0}}^{+} $ to the space of rational
functions on $ {\mathcal Q} $. These operators describe addition of a boson or a fermion
(correspondingly) to a collection of particles {\em relative\/} to a very
populated ground state. \end{nwthrmiii}

\begin{remark} \label{rem3.60}\myLabel{rem3.60}\relax  Note that the given heuristic does not describe why the
ground state which is a limit of collections of particles with more and
more of them should correspond to the {\em Fock ground state}. However, as we
will see later, it is possible to add {\em fermion destruction operators\/} $ \psi_{z_{0}} $
which have correct commutations relations with the operators $ \psi_{z_{0}}^{+} $, and
which annihilate the quantum state corresponding to the classical state
$ \pi=1 $, i.e., the classical ground state.

Some additional arguments are required to find a heuristic for this
coincidence. \end{remark}

\subsection{Fock space II }\label{s16.60}\myLabel{s16.60}\relax  Now we can show that the operators constructed in
Section~\ref{s16.50} can be indeed considered as fermion creation operator.
One needs to construct a subspace of the space of meromorphic functions
on $ {\mathcal Q} $ (the {\em space of states\/}) and to show that there are fermion destruction
operators $ \psi\left(s\right) $ such that
\begin{enumerate}
\item
they satisfy correct (anti)commutation relations with operators
$ \psi^{+}\left(s\right) $;
\item
the negative Laurent coefficients of $ \psi\left(s\right) $ kill the vacuum state;
\item
the negative Laurent coefficients of $ \psi^{+}\left(s\right) $ kill the vacuum state;
\item
all the states may be obtained as combinations of Laurent
coefficients of $ \psi\left(s\right) $ and $ \psi^{+}\left(s\right) $ applied to the vacuum state.
\end{enumerate}

Recall that in the boson case the classical vacuum state corresponds
to having no bosons, thus to the polynomial 1 which has no roots, so the
quantum vacuum state corresponded to the function on $ \cup_{k\geq0}{\mathcal P}_{k} $ which
vanishes on $ {\mathcal P}_{\geq1} $, and is 1 on the one-point set $ {\mathcal P}_{0} $. As seen in Section~%
\ref{s16.50}, introduction of Fock space corresponds to multiplication by an
arbitrary polynomial $ \pi_{S} $, and taking the limit when the degree of this
polynomial goes to infinity. These multiplications (for different $ \pi_{S} $ of
fixed degree $ k $) spread the point $ {\mathcal P}_{0} $ over the whole space $ {\mathcal P}_{k} $. Thus it is
natural to consider a constant function 1 on $ {\mathcal Q}={\mathcal M}\left({\mathbb P}^{1}\right) $ as a natural
analogue of the above boson ground state.

Since this function 1 is translation-invariant, one can see that
$ {\mathbit M}_{\xi}1=1 $, thus $ \psi\left(s\right)1=E_{s} $. On the other hand, $ E_{s} $ depends smoothly on $ s $, thus
negative Laurent coefficients of $ \psi\left(s\right) $ kill 1. This is obviously so in
$ {\mathfrak M}{\mathcal M}\left({\mathbb P}^{1}\right) $.

Similarly, consider the operators $ \psi^{+}\left(s\right) $ acting on meromorphic
functions on $ {\mathcal Q} $. Clearly, $ \psi^{+}\left(s\right)1 =E_{s}^{-1} $. This expression depends smoothly
on $ s $ considered as an element of $ {\mathfrak M}{\mathcal M}\left({\mathbb P}^{1}\right) $.

The (anti)commutation relations for Laurent coefficients of $ \psi\left(s\right) $ and
$ \psi^{+}\left(s\right) $ were already proved in~\eqref{equ6.12}.

Consider now the completeness part. Note that by Corollary~\ref{cor12.40}
the action of Laurent coefficients of $ \psi\left(s\right) $ and $ \psi^{+}\left(s\right) $ preserves the
subspace of $ {\mathfrak M}{\mathcal M}\left({\mathbb P}^{1}\right) $ consisting of polynomials in $ T $, $ T^{-1} $ and $ t_{k} $. This means
that the space of states is contained in the space of these polynomials.

On the other hand, as seen in Section~\ref{s6.31}, the operators $ T\frac{\partial}{\partial T} $,
$ \frac{\partial}{\partial t_{k}} $ and the operators of multiplication by $ t_{k} $ may be expressed via
Laurent coefficients of $ \psi\left(s\right) $ and $ \psi^{+}\left(s\right) $. All these operators preserve the
degree in $ T $, but $ \psi\left(s\right) $ and $ \psi^{+}\left(s\right) $ do not, thus {\em all\/} the polynomials in $ T $,
$ T^{-1} $ and $ t_{k} $ may be obtained by application of combinations of Laurent
coefficients of $ \psi\left(s\right) $ and $ \psi^{+}\left(s\right) $.

This would finish the description of the Fock space as a subspace of
the space of meromorphic functions on $ {\mathcal Q} $, but note that above we had shown
that $ E_{s} $ and $ E_{s}^{-1} $ depend smoothly on $ s $ as elements of $ {\mathfrak M}{\mathcal M}\left({\mathbb P}^{1}\right) $, not as
elements of the space of states. However, the formulae of Corollary~%
\ref{cor12.40} show that this is true in the space of states too.

We conclude that the subspace of $ {\mathfrak M}{\mathcal M}\left({\mathbb P}^{1}\right) $ which corresponds to the
Fock space is a subspace of $ {\mathfrak M}{\mathcal M}\left(\widehat{Z}\right) $, here $ \widehat{Z} $ is the localization of $ {\mathbb P}^{1} $ near
$ 0\in{\mathbb P}^{1} $. Recall that $ {\mathfrak M}{\mathcal M}\left(\widehat{Z}\right) $ consists of rational expressions of $ E_{0}^{\frac{d^{k}}{dz^{k}}} $,
$ k\geq0 $, and the subspace consists of polynomials in $ E_{0}^{\frac{d^{k}}{dz^{k}}} $, $ k\geq0 $, and in
$ E_{0}^{-1} $.

\subsection{Bosons: geometric $ vs $ algebraic } Our approach introduced a geometric
connection between a boson with the creations operators given by Laurent
coefficients of the conformal field $ \alpha\left(s\right)=a^{-1}\left(s\right) $, and a fermion with the
creation operators given by Laurent coefficients of the conformal field
$ \psi^{+}\left(s\right) $. However, usually people consider a different correspondence, which
associate to this fermion a boson with the creation operator given by
Laurent coefficients $ B_{+}\left(s\right) $ and destruction operator given by Laurent
coefficients of $ B_{-}\left(s\right) $, here $ B_{\pm} $ are positive and negative parts of the
field $ B\left(s\right) $ from Section~\ref{s6.31}:
\begin{align} \left(B_{+}\left(s\right)F\right)\xi & = -\frac{\xi'\left(s\right)}{\xi\left(s\right)}F\left(\xi\left(z\right)\right),
\notag\\
\left(B_{-}\left(s\right)F\right)\xi & = - \frac{dF\left(\left(1+\frac{\tau}{z-s}\right)\xi\left(z\right)\right)}{d\tau}|_{\tau=0},
\notag\end{align}
with Laurent coefficients
\begin{align} B_{+}\left(s\right)F\left(T,t\right) & = -\left(\sum_{l\geq1}lt_{l}s^{l-1}\right)F\left(T,t\right),
\notag\\
B_{-}\left(s\right)F\left(T,t\right) & = \sum_{l\geq1}\frac{1}{s^{l+1}}\partial_{t_{l}}F\left(T,t\right)+\frac{1}{s}\frac{\partial_{T}}{T}F\left(T,t\right).
\notag\end{align}
It is easy to see that $ B_{+}\left(s\right)=-{\mathbit E}_{s}^{\frac{d}{ds}}/{\mathbit E}_{s} $, $ B_{-}\left(s\right)=-\frac{d}{d\tau}{\mathbit M}_{1+\frac{\tau}{z-s}} $. It
follows that $ B_{+}\left(r\right) $ commutes with $ B_{+}\left(s\right) $ without fusions, similarly for
$ B_{-}\left(r\right) $ and $ B_{-}\left(s\right) $, $ \left[B_{+}\left(r\right),B_{-}\left(s\right)\right]=\frac{1}{\left(r-s\right)^{2}} $, and $ B_{+}\left(r\right)B_{-}\left(s\right) $ has no fusion,
while the fusions of $ B_{-}\left(s\right)B_{+}\left(r\right) $ are residues of $ \frac{1}{\left(r-s\right)^{2}} $.

Any pair of fields $ B'_{+}\left(s\right) $, $ B'_{-}\left(s\right) $ with the same commutation and
fusion relations may be considered as fields of creation/destruction of a
boson. Since $ \alpha\left(s\right)={\mathbit M}_{z-s} $, $ \alpha\left(s\right) $ and $ \alpha\left(r\right) $ commute for any $ r $ and $ s $ without
fusion, thus they definitely fulfill the share of properties required
from a boson creation field taken separately from boson destruction
field. On the other hand, $ \alpha\left(s\right)^{-1}\frac{d}{ds}\alpha\left(s\right)=B_{-}\left(s\right)! $ Heuristically, this
means that $ \alpha\left(s\right)=C \exp \int B_{-}\left(s\right)ds. $ Thus one should expect that
$ \left[B_{+}\left(r\right),\alpha\left(s\right)\right]=\frac{1}{r-s}\alpha\left(s\right) $, what is easy to check by a direct calculation.
Also, $ B_{+}\left(r\right)\alpha\left(s\right) $ has no fusion.

It is easy to see that if the integral $ \beta\left(s\right)=\int B_{-}\left(s\right)ds=\log \alpha\left(s\right) $
existed, then the pair $ \beta\left(s\right) $, $ \beta^{+}\left(s\right)=-\frac{d}{ds}B_{+}\left(r\right) $ would satisfy the
required properties for a boson destruction/creation fields. Note that
$ \beta^{+}\left(s\right)=\frac{\left(E_{s}^{\frac{d}{ds}}\right)^{2}-E_{s}E_{s}^{\frac{d}{ds}\frac{d}{ds}}}{E_{s}^{2}} $.

Thus we decompose the algebraic boson-fermion correspondence into
the following 3 steps
\begin{equation}
B\left(s\right) \xrightarrow[]{\int} \beta\left(s\right) \xrightarrow[]{\exp } \alpha\left(s\right) \xrightarrow[]{\text{Disc}} \psi^{+}\left(s\right),
\notag\end{equation}
here Disc is an ``extrapolation'' to Dirac vacuum state of the operation of
multiplication by the discriminate Disc (see Section~\ref{s15.30}). As we
have seen it above, the first step is not literally defined for operators
acting in $ {\mathfrak M}{\mathcal M}\left(Z\right) $, but formally $ \beta\left(s\right) $ is the field of destruction of a boson
(which is ``an integral'' of algebraic boson $ B\left(s\right) $). The exponent of this
field becomes a creation operator $ \alpha\left(s\right) $ for a geometric boson, and
geometric boson-fermion correspondence expresses a fermion creation
operator $ \psi^{+}\left(s\right) $ in terms of boson creation operator $ \alpha\left(s\right) $.

Since the first two operations correspond to integration of a
connection, the geometric boson may be considered as a connection, and
the algebraic boson as a differential form corresponding to a connection.
This fact is widely used in physical literature to quantify connections.
I do not know a heuristic which would explain why exponent of $ \beta $-boson
{\em destruction\/} operator appears a geometric boson {\em creation\/} operator.

\bibliography{ref,outref,mathsci}
\end{document}